\documentclass[12pt,onecolumn]{article}
\usepackage{color}
\usepackage{amsmath}
\usepackage{amssymb}
\usepackage{amsfonts}
\usepackage{geometry}
\usepackage{setspace}
\usepackage{amsmath,bm}
\usepackage{graphicx}
\usepackage[section]{placeins}
\usepackage{subfigure}
\usepackage[sort&compress]{natbib}
\usepackage[colorlinks, citecolor=blue]{hyperref}

\author{Shuwei Zhou$^{3,4,5}$, Xiaoying Zhuang$^{4,5}$, Timon Rabczuk$^{1,2*}$}
\title {Phase-field modeling of fluid-driven dynamic cracking in porous media}
\begin{document}

\bibliographystyle{elsarticle-num-names}
\setcitestyle{numbers,square,aysep={},yysep={,},citesep={,}}
\date{}
\maketitle

\spacing {1.2}
\noindent
1 Division of Computational Mechanics, Ton Duc Thang University, Ho Chi Minh City, Viet Nam \\
2 Faculty of Civil Engineering, Ton Duc Thang University, Ho Chi Minh City, Viet Nam \\
3 Institute of Structural Mechanics, Bauhaus-University Weimar, Weimar 99423, Germany\\
4 Department of Geotechnical Engineering, College of Civil Engineering, Tongji University, Shanghai 200092, P.R. China\\
5 Institute of Continuum Mechanics, Leibniz University Hannover, Hannover 30167, Germany.\\
* Corresponding author: timon.rabczuk@tdtu.edu.vn\\

\section*{Research highlight}
\begin{spacing}{1.0}
\begin{itemize}
\item A phase-field modeling of dynamic fracture propagation in porous media is proposed.
\item The phase field method for dynamic cracks in a single-phasic solid is extended for fluid-driven dynamic cracks.
\item The crack propagation and branching is driven by elastic energy.
\item The presented results agree well with existing analytical results.
\item Examples of dynamic crack branching and its interaction with pre-existing natural fractures are presented.
\end{itemize}
\end{spacing}

\begin{abstract}
\noindent A phase field model for fluid-driven dynamic crack propagation in poroelastic media is proposed. Therefore, classical Biot poroelasticity theory is applied in the porous medium while arbitrary crack growth is naturally captured by the phase field model. We also account for the transition of the fluid property from the intact medium to the fully broken one by employing indicator functions. We employ a staggered scheme and implement our approach into the software package COMSOL Multiphysics. Our approach is first verified through three classical benchmark problems which are compared to analytical solutions for dynamic consolidation and pressure distribution in a single crack and in a specimen with two sets of joints. Subsequently, we present several 2D and 3D examples of dynamic crack branching and their interaction with pre-existing natural fractures. All presented examples  demonstrate the capability of the proposed approach of handling dynamic crack propagation, branching and coalescence of fluid-driven fracture. 
\end{abstract}

\noindent Keywords: Phase field, Dynamic crack, Hydraulic fractures, Poroelasticity, COMSOL

\section {Introduction}\label{Introduction}

Hydraulic fracture (HF) is an effective technology for extracting petroleum and natural gas (e.g. shale gas) from reservoirs with low-permeability. The fractures driven by the pressurized fluid form artificial channels, which connect the wellbores with expected resources in the rock matrix. HF promises huge economic benefits because of feasible extraction of vast amounts of resources that have been unexploitable in the past. Hence, fracture propagation in porous media has received increasing attention in mechanical, energy and environmental engineering \citep{mikelic2013phase}. However, hydraulic fracturing sometimes unintentionally creates many extra channels, which facilitates the fracturing fluid or gas to contaminate the groundwater. Therefore, predicting the fracture patterns in porous media is crucial for the application of best HF practices. This requires accurate mathematical models and proper numerical simulation tools that can describe and predict complex fracture behaviors in porous media, such as branching and merging.

Fracture propagation in porous media is a hydro-mechanical coupling process \citep{yang2017hydraulic} and the Biot's theory is often used as a base of the coupling. The flow in the porous media is usually assumed to be laminar \citep{yang2017hydraulic} and modeled as Darcy type \citep{mikelic2013phase} or equivalent Darcy type from Poiseuille flow \citep{miehe2015minimization,miehe2016phase}. 

So far, a large number of discrete, continuous, and hybrid approaches have been developed to model fractures in a solid. In the discrete approach, discontinuities in the displacement field are introduced. Typical discrete approaches exploit remeshing techniques as in \citep{ingraffea1985numerical, areias2016damage, areias2016novel, areias2017effective}, extended finite element method (XFEM) \citep{moes2002extended, chen2012extended},  phantom-node method \citep{song2006method, chau2012phantom, rabczuk2008new}, cohesive element method \citep{zhou2004dynamic, nguyen2001cohesive} or element-erosion method \citep{belytschko1987three, johnson1987eroding}. In the framework of mesh-free method \citep{fu2018boundary}, typical discrete approaches are cracking-particle method \citep{rabczuk2010simple,rabczuk2004cracking} or immersed particle method \citep{rabczuk2010immersed}.  Other discrete approaches are established within the framework of boundary element method (BEM) \citep{wu2015simultaneous, fu2018singular, fu2019semi}, discrete element method (DEM) \citep{SHIMIZU2011712}, and discontinuous deformation analysis (DDA) \citep{NAG:NAG2314}. BEM has relatively low efficiency in handling nonlinear and heterogeneous materials, and also in predicting reliably crack interactions. In DEM and DDA, the calculation domain is discretized into particles or blocks. Fracture forms when the bonds between the particles or blocks are removed. However, the fracture path is not arbitrary but depends on the initial arrangement of the particles and blocks. Furthermore, calibrating the material parameters is more complex.

Continuous approaches to fracture smear the crack over a certain region \citep{santillan2017phase} without introducing strong discontinuities in the displacement field. Typical continuous approaches include gradient damage models \citep{peerlings1996some}, screened-poisson models \citep{areias2016damage,areias2016novel,Areiasinpress}, and phase field models  \citep{miehe2010thermodynamically,miehe2010phase,borden2012phase,hofacker2012continuum, hofacker2013phase, areias2013finite, amiri2014phase, zhou2018phase, zhou2018phase2, zhou2018phase3}. Continuous approaches to fracture do not require complex track cracking algorithms and are much easier to implement compared to discrete approaches. Phase field models (PFM) to fracture can be traced back to \citet{bourdin2008variational} though they were first named by Miehe et al. \citep{miehe2010thermodynamically, miehe2010phase}. Hybrid approach aim to combine the advantages of discrete and continuous approaches. Common hybrid approaches include the DDA-FEM (FEM/DDA) \citep{choo2016hydraulic} and the finite discrete element method (FDEM) \citep{yan2016combined}. Both approaches are more effective than traditional discrete methods but the hybrid approach still cannot solve the intrinsic drawbacks of the discrete methods. For example, fractures from the FEM/DDA simulation only propagate along the element boundaries and the simulation is affected by mesh size and configuration.

For hydraulic fracturing, many discrete, continuous, and hybrid approaches have been successfully exploited and developed. For example, \citet{wu2015simultaneous} used BEM to model 2D fracture propagation. \citet{lecampion2009extended} employed the extended finite element method to simulate hydraulic fracturing by applying fluid pressure along a line fracture. However, some intrinsic difficulties that the hydraulic fracturing faces still exists such as prediction of complex fracture interaction such as branching and joining cracks.  The phase field model is an attractive choice for such problems as 1. the crack is a natural outcome of the simulation and thanks to the thermodynamic framework, it is particularly well suited for coupled problems. 

Some application of phase field models to fluid-driven fractures in porous media are reported for instance in \citep{bourdin2012variational,wheeler2014augmented,mikelic2015quasi, mikelic2015phase,heister2015primal, lee2016pressure,wick2016fluid, yoshioka2016variational,miehe2015minimization,miehe2016phase,ehlers2017phase,santillan2017phase}.  \citet{bourdin2012variational} assumed the material as an impermeable medium while \citet{wheeler2014augmented} extended the phase field model to porous media by introducing poroelastic terms into the energy functional.
The evolution of the fracture domains is treated as a moving boundary problem, and the implementation of their model is facilitated by using a global pressure field in terms of a so-called diffraction system, proposed in \citep{mikelic2015phase3}. In addition, only the mathematical analysis requires explicit knowledge about the fracture boundary.  Later, \citet{mikelic2015quasi, mikelic2015phase} modified the energy functional and fully coupled elasticity, phase field, and pressure. The fracture and intact medium have the same dimension and Biot equations are used for the fluid flow. Moreover, the permeability tensor was modified to form a higher permeability along the fracture. Subsequently, \citet{mikelic2015quasi, mikelic2015phase} suggested adaptive schemes for computational savings. \citet{wick2016fluid, yoshioka2016variational} coupled the phase field model to reservoir simulators. Miehe et al. \citep{miehe2015minimization, miehe2016phase} coupled Darcy-Biot-type flow in poroelastic media and the phase field model. The effective stress in the solid skeleton drives the evolution of the phase field and a stress threshold was set. Recently, \citet{ehlers2017phase} embedded a phase-field approach in the theory of porous media and \citet{santillan2017phase} proposed an immersed-fracture formulation for impermeable porous media. However, these recently developed approaches of phase field modeling in porous media rarely considered inertial effects.

This paper proposes a phase field approach for modeling dynamic fracture propagation in poroelastic media. First, the classical Biot poroelasticity theory is applied in the porous medium. We revisit the phase field method for dynamic fractures in a single-phasic solid \citep{borden2012phase}. Subsequently, we revise the energy functional by adding a fluid pressure-related term and work by external loads and then derive the governing equations in strong form. In addition, the phase field is used as an interpolation function to transit fluid property from the intact medium to the fully broken one. We use COMSOL Multiphysics to implement the proposed approach and adopt a staggered scheme where the displacement, pressure, and phase field are calculated independently. Three examples are performed to verify the feasibility and accuracy of the proposed approach before we present some 2D and 3D examples of dynamic crack branching and their interaction with pre-existing natural fractures. 

The content of this paper is outlined as follows. We present the mathematical models for dynamic fractures in Section 2. In Section 3, we show the numerical implementation of the proposed approach in COMSOL. In Section 4, we verify the numerical simulation by three examples. Section 5 presents 2D and 3D example of dynamic crack branching and Section 6 presents an example of interaction of dynamic hydraulic fracturing with natural cracks. Finally, we end with conclusions regarding our work in Section 7.

\section {Mathematical models of fracture in porous media}\label{Mathematical models of fracture in porous media}
\subsection {Theory of brittle fracture}\label{Theory of brittle fracture}

Let us consider an arbitrary bounded computational domain $\Omega\subset \mathbb R^d$ ($d\in \{2,3\} $) as illustrated in Fig. \ref{Sharp and diffusive crack shape}. The domain contains an internal crack boundary $\Gamma $ and is bounded by an external boundary $\partial \Omega$. We denote the displacement field of the body $\Omega$ at time $t\in[0,T]$ as $\bm u(\bm x,t)\subset \mathbb R^d$ with $\bm x $ being the position vector. The domain $\Omega$ is subjected to  time-dependent Dirichlet boundary conditions, $u_i(\bm x,t)=g_i(\bm x,t)$, on $\partial \Omega_{g_i} \in \Omega$, and the time-dependent von Neumann boundary $\partial \Omega_{h_i} \in \Omega$. The von Neumann conditions impose the traction $\bm f(\bm x,t)$ on $\partial \Omega_{h_i}$. In addition, a body force $\bm b(\bm x,t)\subset \mathbb R^d$ acts throughout the domain. We assume:

\begin{spacing}{1.0}
\begin{itemize}
\item The intrinsic length scale parameter of the phase field is large enough with respect to the pore size.
\item The porous media is linear elastic, homogeneous, and isotropic.
\item The fluid in the media is compressible and viscous.
\end{itemize}
\end{spacing}

For dynamic fracture, the energy functional $\Psi(\bm u,\Gamma)$ of a single-phase solid can be additively decomposed into the kinetic energy $\psi_{kin}$,  elastic energy $\psi_{\varepsilon}(\bm \varepsilon)$, fracture energy, and external work. In this paper, we use the Griffith's theory \citep{francfort1998revisiting} and assume the energy to create a fracture surface per unit area is equal to the critical energy release rate $G_c$. Thus, the energy functional $\Psi(\bm u,\Gamma)$ is written as
	\begin{equation}
	\Psi(\bm u,\Gamma) = \int_{\Omega} \psi_{kin} \mathrm{d}{\Omega}-\int_{\Omega}\psi_{\varepsilon}(\bm \varepsilon) \mathrm{d}{\Omega}-\int_{\Gamma}G_c \mathrm{d}S+\int_{\Omega} \bm b\cdot{\bm u}\mathrm{d}{\Omega}+\int_{\partial\Omega_{h_i}} \bm f\cdot{\bm u}\mathrm{d}S
	\label{functional1}
	\end{equation}

\noindent with the linear strain tensor $\bm\varepsilon = \bm\varepsilon(\bm u)$ given by
	\begin{equation}
	\varepsilon_{ij}=\frac 1 2 \left(\frac{\partial u_i}{\partial x_j}+\frac{\partial u_j}{\partial x_i}\right)
	\end{equation}

The kinetic energy $\psi_{kin}$ is evaluated by
	\begin{equation}
	\psi_{kin}=\psi_{kin}(\dot {\bm u})=\frac 1 2 \rho \dot{u_i} \dot{u_i}
	\end{equation}
\noindent with density $\rho$ and $\dot{\bm u}=\frac{\partial{\bm u}}{\partial t}$.

If the single-phase solid is isotropic and linear elastic, the elastic energy density $\psi_{\varepsilon}(\bm \varepsilon)$ is given by \citep{miehe2010phase}
	\begin{equation}
	\psi_{\varepsilon}(\bm \varepsilon) = \frac{1}{2}\lambda\varepsilon_{ii}\varepsilon_{jj}+\mu\varepsilon_{ij}\varepsilon_{ij}
	\end{equation}

\noindent where $\lambda,\mu>0$ are the Lam\'e constants.

For a porous solid filled with fluid, the effect of the fluid pressure on the energy functional must be considered, see the different forms of the additional pressure-related term in \citet{mikelic2015quasi}, \citet{mikelic2015phase}, and \citet{lee2016pressure}. We follow the formulation of \citet{lee2016pressure} and rewrite Eq. \eqref{functional1} as
	\begin{equation}
	\Psi(\bm u,\Gamma) = \int_{\Omega} \psi_{kin} \mathrm{d}{\Omega}-\int_{\Omega}\psi_{\varepsilon}(\bm \varepsilon) \mathrm{d}{\Omega}+\int_{\Omega}\alpha p \cdot (\nabla \cdot \bm u) \mathrm{d}{\Omega}-\int_{\Gamma}G_c \mathrm{d}S+\int_{\Omega} \bm b\cdot{\bm u}\mathrm{d}{\Omega} + \int_{\partial\Omega_{h_i}} \bm f\cdot{\bm u}\mathrm{d}S
	\label{functional2}
	\end{equation}

\noindent where $p:\Omega\times[0,T]\rightarrow\mathbb R$ is the fluid pressure and $\alpha\in [\epsilon_p,1]$ is the Biot coefficient with $\epsilon_p$ the porosity of the porous media.

\subsection{Phase filed approximation for fracture energy}\label{Phase filed approximation for fracture energy}

\begin{figure}[htbp]
	\centering
	\subfigure[]{\includegraphics[height = 5cm]{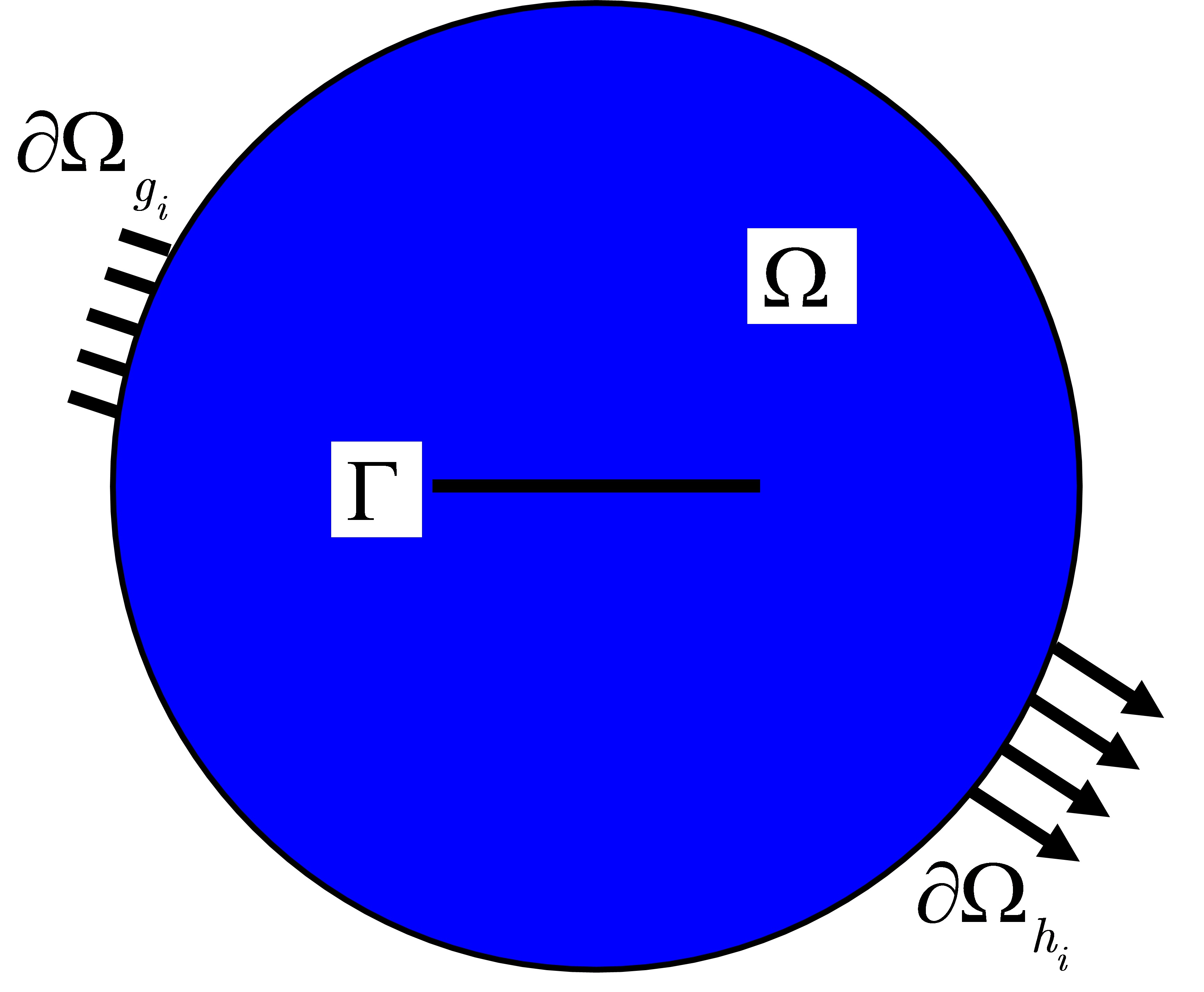}}
	\subfigure[]{\includegraphics[height = 5cm]{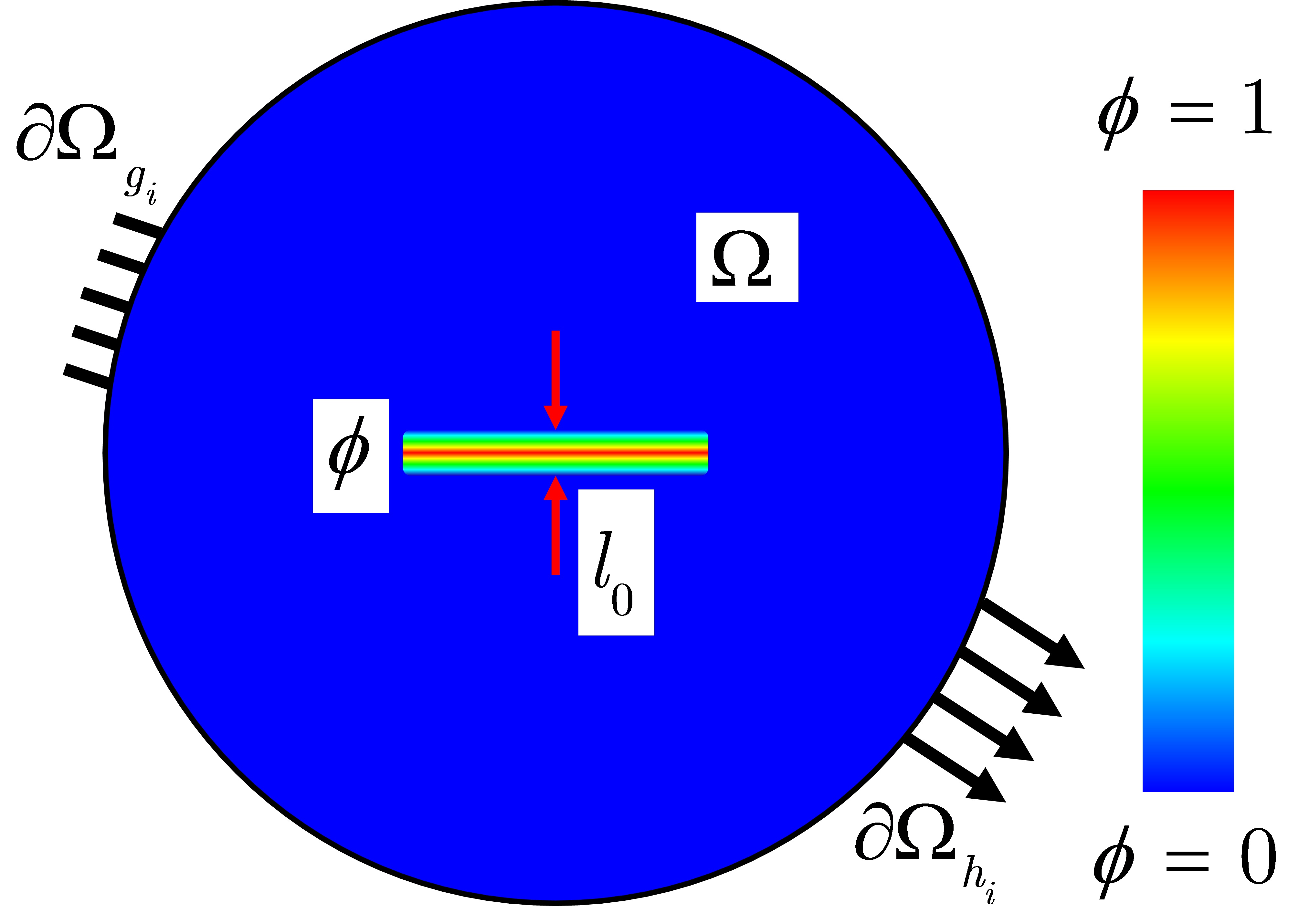}}
	\caption{Sharp and diffusive crack shape. (a) A sharp crack $\Gamma$ and (b) the diffusive phase field.}
	\label{Sharp and diffusive crack shape}
	\end{figure}

In the phase field model (PFM), a scalar field (phase field) is used to diffuse the sharp crack topology \citep{borden2012phase, miehe2010phase, miehe2010thermodynamically} over a certain domain which avoids complex crack tracking procedures and an explicit representation of the crack surface as in discrete crack approaches \cite{moes2002extended}. Therefore, a narrow transition band connects the fully fractured and intact domains with the displacement being still continuous. The phase field $\phi(\bm x,t)\in[0,1]$, which represents the crack as shown in Fig. \ref{Sharp and diffusive crack shape}b should satisfy the following conditions:
	\begin{equation}
	\phi = 
		\begin{cases}
		0,\hspace{1cm}\text{if material is intact}\\1,\hspace{1cm}\text{if material is cracked}
		\end{cases}
	\end{equation}

More details can be found in \citet{miehe2010phase}. A typical one dimensional phase field is approximated with the exponential function:
	\begin{equation}
	\phi(x)= e^{-|x|/l_0}
	\end{equation}

\noindent  with length scale parameter $l_0$ which controls the transition region between the fracture and intact material. 

For 2D and 3D problems, the crack surface density per unit volume of the solid is given by \cite{miehe2010phase}	
	\begin{equation}
	\gamma(\phi,\bigtriangledown\phi)=\frac{\phi^2}{2l_0}+\frac{l_0}2\frac{\partial\phi}{\partial x_i}\frac{\partial\phi}{\partial x_i}
	\label{phase field approximation}
	\end{equation}

Thus, exploiting Eq. \eqref{phase field approximation}, the fracture energy in Eq. \eqref{functional2} can be rewritten as
	\begin{equation}
	\int_{\Gamma}G_c \mathrm{d}S=\int_{\Omega}G_c\left[\frac{\phi^2}{2l_0}+\frac{l_0}2\frac{\partial\phi}{\partial x_i}\frac{\partial\phi}{\partial x_i}\right]\mathrm{d}{\Omega}
	\label{phase field approximation for the fracture energy}
	\end{equation}

\subsection{Governing equations for evolution of the phase field}\label{Governing equations for evolution of the phase field}

We use the variational approach to obtain the governing equations for dynamic crack propagation in porous media. In the PFM \citep{miehe2010thermodynamically}, the elastic energy drives crack propagation and can be transformed into fracture surface energy during cracking. Thus, energy decomposition is made to ensure cracking only under tension. We follow \citet{miehe2010thermodynamically} and the elastic energy is decomposed into positive and compressive parts. The phase field is only related to the positive elastic energy and evolution of the phase field under compression is not allowed. Therefore, the strain tensor $\bm\varepsilon$  is decomposed as $\bm\varepsilon=\bm\varepsilon_++\bm\varepsilon_-$ and
	\begin{equation}
	  \left\{
	   \begin{aligned}
	\bm\varepsilon_+=\sum_{a=1}^d \langle\varepsilon_a\rangle_+\bm n_a\otimes\bm n_a \\ \bm\varepsilon_-=\sum_{a=1}^d \langle \varepsilon_a\rangle_-\bm n_a\otimes\bm n_a
	   \end{aligned}\right.
	\end{equation}

\noindent where $\bm\varepsilon_+$  and $\bm\varepsilon_-$  are the tensile and compressive parts of the strain tensor, respectively; $\varepsilon_a$  and  $\bm n_a$ are the principal strains and their directions. The operators $\langle\centerdot\rangle_+$  and $\langle\centerdot\rangle_-$  are defined as \citep{miehe2010phase}: $\langle\centerdot\rangle_+=(\centerdot+|\centerdot|)/2$, and $\langle\centerdot\rangle_-=(\centerdot-|\centerdot|)/2$.

Using the decomposed strain tensor, the tensile and compressive parts of the elastic energy density are given by
	\begin{equation}
	  \left\{
	   \begin{aligned}
	\psi_{\varepsilon}^+(\bm \varepsilon) = \frac{\lambda}{2}\langle tr(\bm\varepsilon)\rangle_+^2+\mu tr \left(\bm\varepsilon_+^2\right) 
	\\ \psi_{\varepsilon}^-(\bm \varepsilon) = \frac{\lambda}{2}\langle tr(\bm\varepsilon)\rangle_-^2+\mu tr \left(\bm\varepsilon_-^2\right) 
	   \end{aligned}\right.
	\end{equation}

We follow \citet{borden2012phase} and assume that the phase field affects only the tensile part of the elastic energy density. A quadratic equation is used and the stiffness reduction is modeled by the following equation:
	\begin{equation}
	\psi_{\varepsilon}(\bm\varepsilon)=\left[(1-k)(1-\phi)^2+k\right]\psi_{\varepsilon}^+(\bm \varepsilon)+\psi_{\varepsilon}^-(\bm \varepsilon)
	\label{decomposition of the elastic energy}
	\end{equation}

\noindent where $0<k\ll1$ is a  parameter that prevents the tensile part of the elastic energy density from disappearing and avoids  the numerical singularity when the phase field $\phi$  tends to 1. Taking advantage of Eqs. \eqref{phase field approximation for the fracture energy}  and \eqref{decomposition of the elastic energy},  Eq. \eqref{functional2} can be rewritten as
	\begin{multline}
	L=\Psi(\bm u,\Gamma)=\frac 1 2 \int_{\Omega}\rho \dot{u_i} \dot{u_i} \mathrm{d}{\Omega}-\int_{\Omega}\left\{\left[(1-k)(1-\phi)^2+k\right]\psi_{\varepsilon}^+(\bm \varepsilon)+\psi_{\varepsilon}^-(\bm \varepsilon)\right\}\mathrm{d}{\Omega}+\\ \int_{\Omega}\alpha p \cdot (\nabla \cdot \bm u) \mathrm{d}{\Omega}-\int_{\Omega}G_c\left[\frac{\phi^2}{2l_0}+\frac{l_0}2\frac{\partial\phi}{\partial x_i}\frac{\partial\phi}{\partial x_i}\right]\mathrm{d}{\Omega}+ \int_{\Omega} b_iu_id{\Omega}+\int_{\partial\Omega_{h_i}} f_iu_i\mathrm{d}S
	\label{final functional}
	\end{multline}

Crack initiation, propagation and branching of the crack $\Gamma(\bm x,t)$ at  time $t\in[0,T]$ for $x\in\Omega$ occurs when the functional achieves an extreme value. Hence, we calculate the first variation of the functional $L$ and set it zero. After assembling all the items related to the variation of displacement and phase field, we obtain
	\begin{equation}
	  \left\{
	   \begin{aligned}
	\frac {\partial {\sigma_{ij}^{por}}}{\partial x_j}+b_i=\rho\ddot{u_i}, &\hspace{0.5cm} in \hspace{0.1cm} \Omega\times(0,T]
	\\ \left[\frac{2l_0(1-k)\psi_{\varepsilon}^+}{G_c}+1\right]\phi-l_0^2\frac{\partial^2 \phi}{\partial {x^2}}=\frac{2l_0(1-k)\psi_{\varepsilon}^+}{G_c}, &\hspace{0.5cm} in \hspace{0.1cm} \Omega\times(0,T]
	   \end{aligned}\right.
	\label{governing equations 0}
	\end{equation}

\noindent where $\sigma_{ij}^{por}$ are the components of the Cauchy stress tensor $\bm\sigma^{por}$ and 
	\begin{equation}
	\bm \sigma^{por}(\bm\varepsilon)=\bm \sigma(\bm\varepsilon)-\alpha p \bm I,\hspace{0.5cm} in \hspace{0.1cm} \Omega\times(0,T]
	\end{equation}
\noindent with $\bm I$ the identity tensor $\in \mathbb R^{d\times d}$ and $\sigma_{ij}$ component of the effective linear elastic stress tensor $\bm \sigma(\bm\varepsilon)$. The effective stress $\bm \sigma(\bm\varepsilon)$ is calculated by
	\begin{equation}
	\sigma_{ij}=\left [(1-k)(1-\phi)^2+k \right]\frac {\partial{\psi_\varepsilon^+}}{\partial {\varepsilon_{ij}}}+\frac {\partial{\psi_\varepsilon^-}}{\partial {\varepsilon_{ij}}}
	\end{equation}

	\begin{equation}
	\bm \sigma=\left [(1-k)(1-\phi)^2+k \right]\left[\lambda \langle tr(\bm\varepsilon)\rangle_+ \bm I+ 2\mu \bm\varepsilon_+ \right]+\lambda \langle tr(\bm\varepsilon)\rangle_- \bm I+ 2\mu \bm\varepsilon_-
	\end{equation}

The irreversibility condition $\Gamma(\bm x,s)\in\Gamma(\bm x,t)(s<t)$ is required for the phase field model, meaning cracks cannot be recovered to uncracked states. To ensure a monotonically increasing phase field, a strain-history field method \citep{miehe2010phase,miehe2010thermodynamically,borden2012phase} is used to ensure the irreversibility condition during compression or unloading. In this paper, the following strain-history field  $H(\bm x,t)$ is introduced:

	\begin{equation}
	H(\bm x,t) = \max \limits_{s\in[0,t]}\psi_\varepsilon^+\left(\bm\varepsilon(\bm x,s)\right), \hspace{0.5cm} \mathrm{in} \hspace{0.1cm} \Omega\times(0,T]
	\end{equation}

Note that the history field $H(\bm x,t)$ satisfies the Kuhn-Tucker condition \citep{borden2012phase} during loading and unloading. Therefore, replacing $\psi_\varepsilon^+$  by  $H(\bm x,t)$  in Eq. \eqref{governing equations 0}, the strong forms of the displacement and phase field are rewritten as
	\begin{equation}
	  \left\{
	   \begin{aligned}
	\frac {\partial {\sigma_{ij}^{por}}}{\partial x_i}+b_i=\rho \ddot{u_i}, &\hspace{0.5cm} \mathrm{in} \hspace{0.1cm} \Omega\times(0,T]
		\\ \left[\frac{2l_0(1-k)H}{G_c}+1\right]\phi-l_0^2\frac{\partial^2 \phi}{\partial {x^2}}=\frac{2l_0(1-k)H}{G_c}, &\hspace{0.5cm} \mathrm{in} \hspace{0.1cm} \Omega\times(0,T]
	\label{governing equation1}
	   \end{aligned}\right.
	\end{equation}

\noindent with
	\begin{equation}
	  \left\{
	   \begin{aligned}
	&\sigma_{ij}^{por}m_j=f_i, \hspace{1cm} &\mathrm{on}\hspace{0.5cm} \partial\Omega_{h_i}\times(0,T]
	\\ &\frac{\partial \phi}{\partial x_i} m_i = 0, &\mathrm{on}\hspace{0.5cm} \partial\Omega\times(0,T]
	\label{boundary condition of the phase field}
	\end{aligned}\right.
	\end{equation}
\noindent where $m_j$ are the components of the outward-pointing normal vector of the boundary.

\subsection{Governing equations for fluid pressure}\label{Governing equations for fluid pressure}

The key novelty of this paper is to couple the dynamic phase field formulation to a flow field in order to study  fluid-driven fracture problems. Therefore, we assume Darcy flow in the porous domain $\Omega$. Other more complicated flow fields will be considered in future work. The domain is subdivided into three parts: $\Omega_R(t)$, $\Omega_F(t)$ and $\Omega_T(t)$. $\Omega_R(t)$ represents the reservoir domain (unbroken domain) and $\Omega_F(t)$ is the fractured domain. $\Omega_T(t)$ is the transition domain betwen $\Omega_R(t)$ and $\Omega_F(t)$. In this paper, we follow \citet{lee2016pressure} and use the phase field as an indicator function to separate the three flow domains. 

First of all, two thresholds $c_1$ and $c_2$ are set. A subdomain is considered as the reservoir domain $\Omega_R(t)$ if $\phi\le c_1$ and as the fracture domain $\Omega_F(t)$ if $\phi\ge c_2$. In the transition domain, $c_1<\phi<c_2$. Note that in engineering reservoir scale, the singular limit that corresponds to history terms in dynamic Biot system \citep{mikelic2012theory} is extremely small and a quasi-static Biot system is obtained. Therefore, in the reservoir domain $\Omega_R(t)$, mass conservation is expressed as
	\begin{equation}
	\frac{\partial}{\partial t}(\varepsilon_{pR}\rho_R)+\nabla\cdot(\rho_R\bm v_R)=q_R-\rho_R\alpha_R\frac{\partial \varepsilon_{vol}}{\partial t}
	\label{mass conservation of the reservoir domain}
	\end{equation}  

\noindent where $\rho_R$, $q_R$, $\varepsilon_{pR}$, and $\alpha_R$ are the density of fluid, source term, porosity, and Biot coefficient in the reservoir domain, respectively; $\varepsilon_{vol}=\nabla\cdot\bm u$ is the volumetric strain of $\Omega_R(t)$. 

Darcy's law related the fluid velocity in $\Omega_R(t)$ to the pressure gradient:
	\begin{equation}
	\bm v_R=-\frac{K_R}{\mu_R}(\nabla p+\rho_R\bm g)
	\label{velocity of the reservoir domain}
	\end{equation} 

\noindent where $K_R$ and $\mu_R$ are the permeability and fluid viscosity of $\Omega_R(t)$, respectively; $\bm g$ is the gravity vector. Taking advantage of the storage model in \citep{biot1962mechanics}, we have 
	\begin{equation}
	\frac{\partial}{\partial t}(\varepsilon_{pR}\rho_R) = \rho_R S_R \frac{\partial p}{\partial t}
	\end{equation}

\noindent where $S_R$, the storage coefficient of $\Omega_R$, is given by
 	\begin{equation}
	S_R=\varepsilon_{pR}c_R+\frac{(\alpha_R-\varepsilon_{pR})(1-\alpha_R)}{K_{VR}}
	\end{equation}

\noindent with $c_R$ the fluid compressibility and $K_{VR}$ the bulk modulus of the reservoir domain. Thus, the equation of mass conservation \eqref{mass conservation of the reservoir domain} reads
	\begin{equation}
	\rho_R S_R \frac{\partial p}{\partial t}+\nabla\cdot(\rho_R\bm v_R)=q_R-\rho_R\alpha_R\frac{\partial \varepsilon_{vol}}{\partial t}
	\label{governing equation of the reservoir domain}
	\end{equation}  

In the fracture domain $\Omega_F(t)$, the volumetric strain $\varepsilon$ vanishes from the equation of mass conservation:
	\begin{equation}
	\rho_F S_F \frac{\partial p}{\partial t}+\nabla\cdot(\rho_F\bm v_F)=q_F
	\label{governing equation of the fracture domain}
	\end{equation} 

\noindent where $\rho_F$, $S_F$, and $q_F$ are the fluid density, storage coefficient, and source term in the fracture domain $\Omega_F(t)$. 

The storage coefficient $S_F$ is equal to the fluid compressibility $c_F$ and the Darcy's velocity in $\Omega_F(t)$ is given by 
	\begin{equation}
	\bm v_F=-\frac{K_F}{\mu_F}(\nabla p+\rho_F\bm g)
	\label{velocity of the fracture domain}
	\end{equation} 

\noindent where $K_F$ and $\mu_F$ denote the permeability and fluid viscosity of $\Omega_F(t)$, respectively. The transition functions link the governing equations of the reservoir and fracture domains. For simplicity, we follow \citet{lee2016pressure} and define two linear indicator functions: $\chi_R$ and $\chi_F$:
	\begin{equation}
	\chi_R(\cdot,\phi):= \chi_R(\bm x, t,\phi)=1 \quad in \quad \Omega_R(t),\quad \mathrm{and} \quad \chi_R(\cdot,\phi)= 0 \quad in \quad \Omega_F(t)
	\end{equation} 
	\begin{equation}
		\chi_F(\cdot,\phi):= \chi_F(\bm x, t,\phi)=1 \quad in \quad \Omega_F(t),\quad \mathrm{and} \quad \chi_F(\cdot,\phi)= 0 \quad in \quad \Omega_R(t)
		\end{equation} 

In the transition domain, the indicator functions depend on the phase field as follows:
	\begin{equation}
	\chi_R(\cdot,\phi)=\frac{c_2-\phi}{c_2-c_1},\quad \mathrm{and} \quad \chi_F(\cdot,\phi)=\frac{\phi-c_1}{c_2-c_1}
	\label{functions}
	\end{equation} 

\begin{figure}[htbp]
	\centering
	\subfigure[]{\includegraphics[height = 5cm]{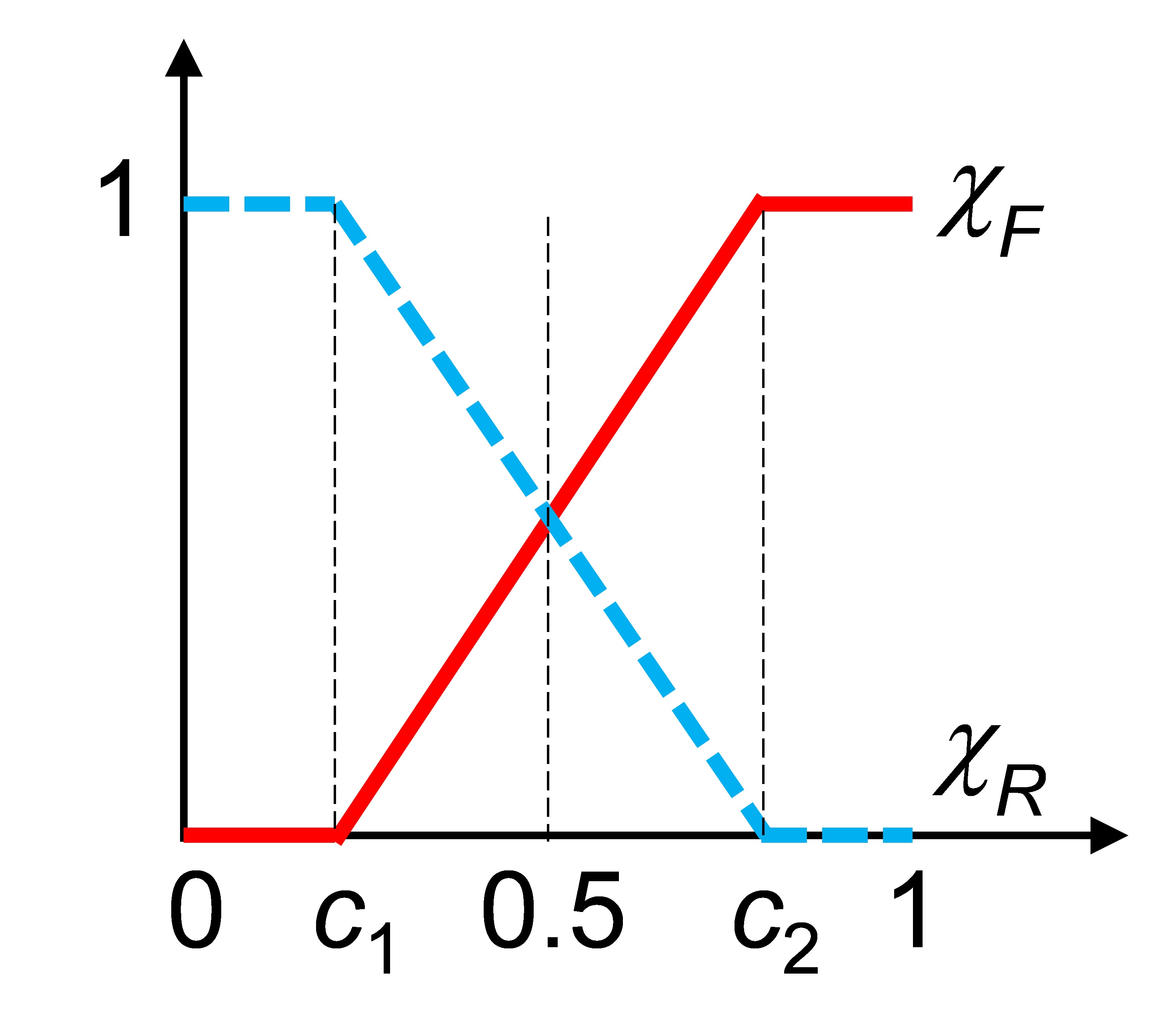}}\hspace{1cm}
	\subfigure[]{\includegraphics[height = 5cm]{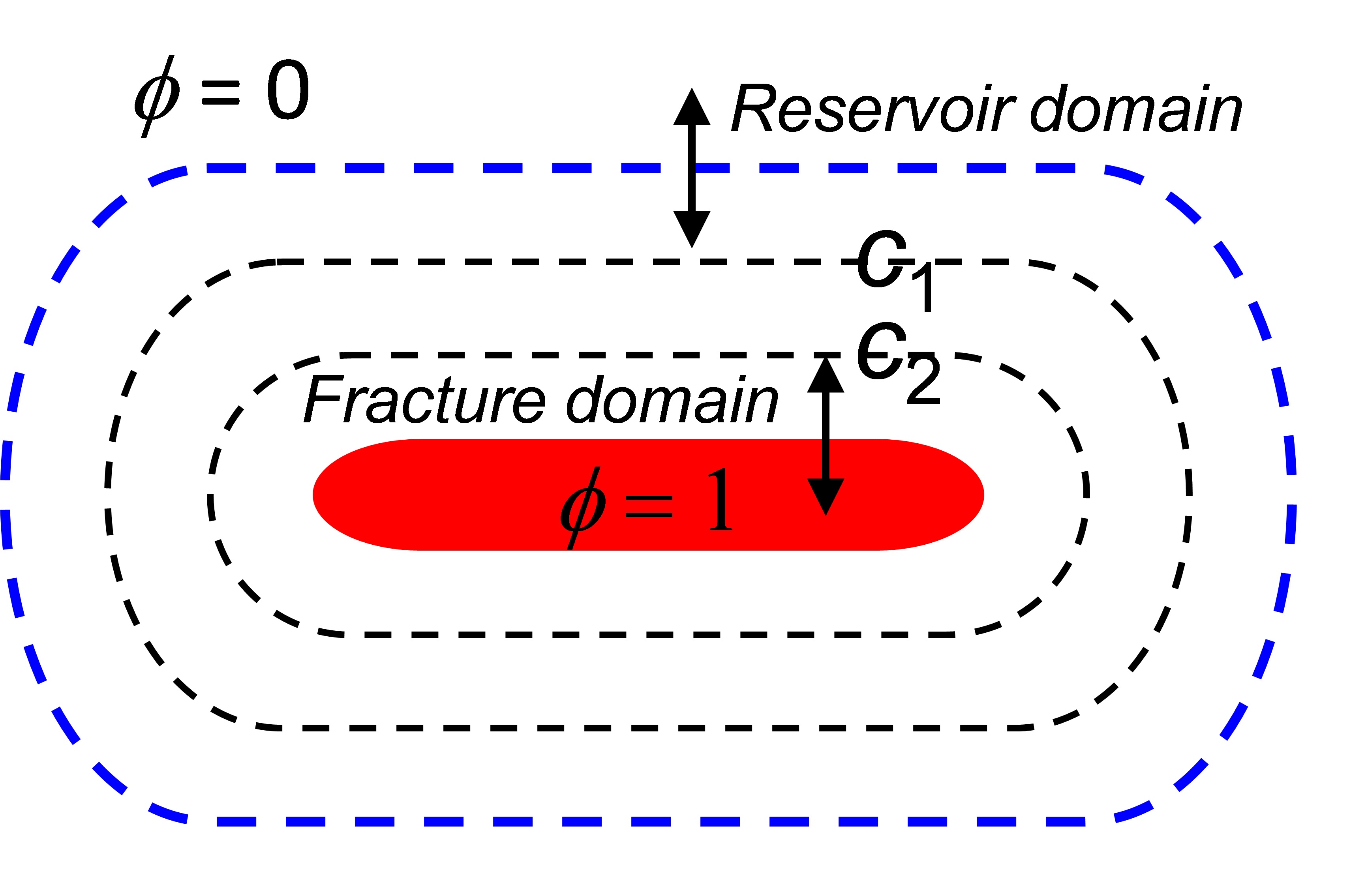}}
	\caption{(a) Linear indicator functions $\chi_R$ and $\chi_F$ and (b) the reservoir and fracture domains}
	\label{(a) Linear indicator functions and (b) the reservoir and fracture domains}
	\end{figure}

Figure \ref{(a) Linear indicator functions and (b) the reservoir and fracture domains} shows the linear indicator functions $\chi_R$ and $\chi_F$. The fracture, reservoir and transition domains are also illustrated based on the thresholds of the phase field in Fig. \ref{(a) Linear indicator functions and (b) the reservoir and fracture domains}. Thus, we obtain the fluid and solid properties of the transition domain $\Omega_T$ from interpolation of the reservoir and fracture domains with the indicator functions $\chi_R$ and $\chi_F$. Then, the mass conservation in the transition domain is given by
	\begin{equation}
	\rho S \frac{\partial p}{\partial t}+\nabla\cdot(\rho\bm v)=q_m-\rho\alpha\chi_R\frac{\partial \varepsilon_{vol}}{\partial t}
	\label{mass conservation of the whole domain}
	\end{equation}  

\noindent  with $\rho=\rho_R\chi_R+\rho_F\chi_F$, $\alpha=\alpha_R\chi_R+\alpha_F\chi_F$ and $q_m$ designates the source term. The storage coefficient $S$ is then replaced by
 	\begin{equation}
	S=\varepsilon_pc+\frac{(\alpha-\varepsilon_p)(1-\alpha)}{K_{VR}}
	\end{equation}

\noindent with $c=c_R\chi_R+c_F\chi_F$. Note that $\varepsilon_p=0$ and $\alpha=1$ for the fracture domain and thereby $\varepsilon_p=\varepsilon_{pR}\chi_R$ and $\alpha=\alpha_{R}\chi_R+\chi_F$. 

Now, we use Eq. \eqref{mass conservation of the whole domain} as one of the main equations for the whole domain $\Omega$ because it can be degenerated into Eq. \eqref{governing equation of the reservoir domain} for $\Omega_R$ and Eq. \eqref{governing equation of the fracture domain} for $\Omega_F$.  The Darcy's velocity $\bm v$ is then calculated by
	\begin{equation}
	\bm v=-\frac{K}{\mu}(\nabla p+\rho\bm g)
	\label{velocity of the whole domain}
	\end{equation} 

\noindent where $K=K_{R}\chi_R+K_F\chi_F$ is the effective permeability and $\mu=\mu_{R}\chi_R+\mu_F \chi_F$ is the effective fluid viscosity. Finally, we rewrite the governing equation for the flow field in the porous domain $\Omega$ in terms of the fluid pressure $p$:
	\begin{equation}
	\rho S \frac{\partial p}{\partial t}-\nabla\cdot \frac{\rho K}{\mu}(\nabla p+\rho\bm g)=q_m-\rho\alpha\chi_R\frac{\partial \varepsilon_{vol}}{\partial t}
	\label{governing equation of the whole domain}
	\end{equation}
	
Note that Eqs. \ref{mass conservation of the reservoir domain} to \ref{governing equation of the whole domain} give smooth transition for the governing equation and fluid property form the fully fractured to the intact reservoir domains by using the indicator functions and in the simulation all the fields are continuous, thereby avoiding special treatment for the singularity at the fracture tip \citep{mikelic2015phase3}. In addition, the focus of this work is to present the framework of a phase-field modeling approach for dynamic fluid-driven cracks that can be easily implemented and applied. Thus, the determination of the coefficients including the permeability model that describe hydro-mechanical responses and fracture behavior will be tackled in future research. For example, the permeability of the fractured domain on crack opening \citep{miehe2016phase, mikelic2015phase3} or on the volumetric strain \citep{Zhuang2017} can be applied. However, as will be described in Section \ref{Numerical implementation}, the PFM is implemented within the COMSOL environment where the crack opening cannot be extracted in a straight forward manner due to smeared representation of the sharp crack shape. Therefore, for simplicity, we apply an unchanged fluid property for the fracture domain in our presented examples, which also show favorable results.

\subsection{Initial and boundary conditions}\label{Initial and boundary conditions}

The following initial conditions are imposed:
	\begin{equation}
	  \left\{
	   \begin{aligned}
	&\bm u(\bm x,0)=\bm u_0(\bm x)\hspace{2cm} &\bm x\in\Omega
     \\&\bm v(\bm x,0)=\bm v_0(\bm x)\hspace{2cm} &\bm x\in\Omega
	\\ &p(\bm x,0)=p_0(\bm x)\hspace{2cm} &\bm x\in\Omega
	\\&\phi(\bm x,0)=\phi_0(\bm x)\hspace{2cm} &\bm x\in\Omega
	\end{aligned}\right.
	\label{Initial condition}
	\end{equation}

The initial phase field  $\phi_0=1$ in a local domain models a pre-existing crack \citep{borden2012phase}. For the displacement field and phase field, the boundary conditions are given in Subsections \ref{Theory of brittle fracture} and \ref{Governing equations for evolution of the phase field}. Likewise, the Dirichlet  boundary condition on $\partial\Omega_D$ and von Neumann boundary condition on $\partial\Omega_N$ with $\partial\Omega_D\cap \partial\Omega_N =\emptyset$ are prescribed for the fluid pressure field:
	\begin{equation}
	p=p_D \hspace{2cm}\mathrm{on}\quad\partial\Omega_D\times(0,T]
	\end{equation}
	\begin{equation}
	-\bm m\cdot \rho\bm v=M_N  \hspace{2cm}\mathrm{on}\quad\partial\Omega_N\times(0,T]
	\end{equation}

\noindent with the prescribed pressure $p_D$  on the Dirichlet boundary and $M_N$ is the mass flux on the Neumann boundary.

\section {Numerical implementation}\label{Numerical implementation}
\subsection{Finite element discretization}\label{Finite element discretization}
Find ${\bf u} \in \mathcal{U}$ $\forall \delta {\bf u} \in \mathcal{U}_0$,  ${\phi} \in \mathcal{V}$ $\forall \delta {\phi} \in \mathcal{V}_0$ and $p \in \mathcal{P}$ $\forall \delta p \in \mathcal{P}_0$ such that the weak forms of the governing equations are given by: 
	\begin{equation}
	\int_{\Omega}\left[-\rho\ddot{\bm u}\cdot\delta\bm u-(\bm\sigma-\alpha p\bm I):\delta \bm {\varepsilon}\right] \mathrm{d}\Omega +\int_{\Omega}\bm b \cdot \delta \bm u \mathrm{d}\Omega +\int_{\Omega_{h_i}}\bm f \cdot \delta \bm u \mathrm{d}S=0
	\label{weak form 1}
	\end{equation}
\noindent ,
	\begin{equation}
	\int_{\Omega}-2(1-k)H(1-\phi)\delta\phi\mathrm{d}\Omega+\int_{\Omega}G_c\left(l_0\nabla\phi\cdot\nabla\delta\phi+\frac{1}{l_0}\phi\delta\phi\right)\mathrm{d}\Omega=0
	\label{weak form 2}
	\end{equation}
\noindent and
	\begin{equation}
	\int_{\Omega} \rho S \frac{\partial p}{\partial t}\delta p\mathrm{d}\Omega-\int_{\Omega} \rho \bm v \cdot\nabla\delta p \mathrm{d}\Omega=\int_{\partial\Omega}M_n\mathrm{d}S+\int_{\Omega}\left(q_m-\rho\alpha\chi_R\frac{\partial \varepsilon_{vol}}{\partial t}\right)\mathrm{d}\Omega
	\label{weak form 3}
	\end{equation}  
	
\noindent with approximation spaces
\begin{eqnarray}	
  \mathcal{U}&=& \left\{ {\bf u} \in \mathcal{C}^0 | {\bf u} = \bar{\bf u} \,\, on \,\ \partial \Omega_D    \right\} \nonumber \\
    \mathcal{U}_0&=&  \left\{ \delta {\bf u} \in \mathcal{C}^0 | \delta {\bf u} = 0 \,\, on \,\ \partial \Omega_D   \right\} \nonumber \\
      \mathcal{V}&=& \left\{ {\phi} \in \mathcal{C}^0 | {\phi} = \bar{\phi} \,\, on \,\ \partial \Omega_c    \right\} \nonumber \\
    \mathcal{V}_0&=&  \left\{ \delta {\phi} \in \mathcal{C}^0 | \delta {\phi} = 0 \,\, on \,\ \partial \Omega_c   \right\} \nonumber \\    
  \mathcal{P}&=&   \left\{ p \in \mathcal{C}^{-1} | p = \bar{p} \,\, on \,\ \partial \Omega_t   \right\} \nonumber \\
  \mathcal{P}_0&=&   \left\{  \delta p \in \mathcal{C}^{-1} | \delta p = 0 \,\, on \,\ \partial \Omega_t   \right\}
\end{eqnarray}

Defining the nodal values for the three fields ($\bm u$, $\phi$, and $p$) with $\bm u_i$, $\phi_i$, and $p_i$, their approximation can be written as
	\begin{equation}
	\bm u = \sum_i^{n}N_i \bm u_i,\hspace{0.5cm} \phi = \sum_i^{n}N_i \phi_i,\hspace{0.5cm} p = \sum_i^{n}N_i p_i
	\end{equation}
\noindent where $n$ is the number of nodes in each element and $N_i$ the shape function of node $i$. The gradients of the three fields are given by
	\begin{equation}
	\bm \varepsilon = \sum_i^{n}\bm B_i^u \bm u_i,\hspace{0.5cm} \nabla\phi = \sum_i^{n}\bm B_i^{\phi} \phi_i,\hspace{0.5cm} \nabla p = \sum_i^{n}\bm B_i^{p} p_i
\label{trial functions}
	\end{equation}
\noindent where $\bm B_i^u$, $\bm B_i^\phi$, and $\bm B_i^p$ are matrices containing the derivatives of the shape functions:
	\begin{equation}
	\bm B_i^u=\left[
	\begin{array}{cc}
	N_{i,x}&0\\
	0&N_{i,y}\\
	N_{i,y}&N_{i,x}
	\end{array}\right],
	\hspace{0.5cm}\bm B_i^{\phi}=\bm B_i^p=\left[
	\begin{array}{ccc}
	N_{i,x}\\
	N_{i,y}\\
	\end{array}\right]\label{gradient of shape function 2D}
	\end{equation}

Substituting the trial functions, Eqs. \eqref{trial functions} and associated test functions which have a similar structure into the weak form, Eqs. \eqref{weak form 1} to \eqref{weak form 3} leads the following system of equation:
	\begin{equation} 
	\left\{
	\begin{aligned}
	\bm R_i^{u}&=\bm F_i^{u,ext}-\bm F_i^{u,int}-\bm F_i^{u,ine}\\
	R_i^{\phi}&=- F_i^{\phi,int}\\
	R_i^{p}&=F_i^{p,ext}- F_i^{p,int}- F_i^{p,vis}
	\end{aligned}
	\right.
	\end{equation}

\noindent where $\bm R_i^{u}$, $R_i^{\phi}$, and $R_i^{p}$ are residuals of the three fields with external force vector $\bm F_i^{u,ext}$, inner force vector $\bm F_i^{u,int}$ and inertia force vector $\bm F_i^{u,ine}$ described by
	\begin{equation}
	\left\{
	\begin{aligned}
	\bm F_i^{u,ext} &= \int_{\Omega}N_i\bm b \mathrm{d}\Omega+ \int_{\Omega_{h_i}}N_i\bm f \mathrm{d}S+ \int_{\Omega}[\bm B_{i}^{u}]^{\mathrm T}\alpha p \bm I \mathrm{d}\Omega\\
	\bm F_i^{u,int} &= \int_{\Omega}[\bm B_i^u]^{\mathrm T}\bm\sigma \mathrm{d}\Omega\\
	\bm F_i^{u,ine} &= \int_{\Omega}\rho N_i \ddot{\bm u} \mathrm{d}\Omega
	\end{aligned}
	\right.
	\end{equation}

The inner force term of the phase field is given by
	\begin{equation} 
	F_i^{\phi,int} = \int_{\Omega} \left\{-2(1-k)(1-\phi)H N_i+G_c\left(l_0[\bm B_i^{\phi}]^{\mathrm T}\nabla \phi+\frac 1 {l_0}\phi N_i \right )\right\} \mathrm{d}\Omega
	\end{equation}

Neglecting gravity, we  derive the inner force $ F_i^{p,int}$, viscous force $F_i^{p,vis}$, and external force $F_i^{p,ext}$ of the pressure field as follows
	\begin{equation}
	\left\{
	\begin{aligned}
	F_i^{p,int} &= \int_{\Omega}[\bm B_i^p]^{\mathrm T}\frac{\rho K}{\mu}\nabla p \mathrm{d}\Omega\\
	F_i^{p,vis} &= \int_{\Omega}N_i \rho S \frac{\partial p}{\partial t} \mathrm{d}\Omega\\
	F_i^{p,ext} &=  \int_{\Omega}N_i \left( q_m-\rho\alpha\chi_R\frac{\partial \varepsilon_{vol}}{\partial t} \right) \mathrm{d}\Omega+ \int_{\partial\Omega_{N}}N_i M_N \mathrm{d}S
	\end{aligned}
	\right.
	\end{equation}


In this paper, we use the staggered scheme to solve for the displacement, phase field and fluid pressure. The Newton-Raphson approach is adopted and the tangents on the element level are calculated by
		\begin{equation}
		\left\{
		\begin{aligned} 
		\bm K_{ij}^{uu}&=\frac{\partial \bm F_i^{u,int}}{\partial \bm u_j}=\int_{\Omega}[\bm B_i^u]^{\mathrm T}\bm D [\bm B_j^u] \mathrm{d} \Omega\\
		\bm K_{ij}^{\phi\phi}&=\frac{\partial F_i^{\phi,int}}{\partial \phi_j}=\int_{\Omega}\left\{[\bm B_i^\phi]^{\mathrm T} G_cl_0 [\bm B_j^\phi]+N_i\left(2(1-k)H+\frac{G_c}{l_0}\right)N_j\right\} \mathrm{d} \Omega\\
		\bm K_{ij}^{pp}&=\frac{\partial F_i^{p,int}}{\partial p_j}=\int_{\Omega}[\bm B_i^p]^{\mathrm T} \frac{\rho K}{\mu} [\bm B_j^p] \mathrm{d} \Omega
		\end{aligned}
		\right.
		\end{equation}

\noindent where $\bm D$ is the fourth order elasticity tensor given by
	\begin{equation}\begin{aligned}
	\bm D = \frac {\partial \bm\sigma}{\partial \bm\varepsilon}=\lambda\left\{ \left[(1-k)(1-\phi)^2+k \right]H_\varepsilon(tr(\bm\varepsilon))+H_\varepsilon(-tr(\bm\varepsilon))\right\}\bm J +\\
2\mu\left\{\left[(1-k)(1-\phi)^2+k \right]\frac{\partial \bm\varepsilon_+}{\partial \bm\varepsilon}+\frac{\partial \bm\varepsilon_-}{\partial \bm\varepsilon}\right\}
	\end{aligned}
	\end{equation}

\noindent where $H_\varepsilon \langle x \rangle$  is the Heaviside function: $H_\varepsilon \langle x \rangle=1$  if  $x>0$ and $H_\varepsilon \langle x \rangle=0$  if $x\leq 0$, and $J_{ijkl}=\delta_{ij}\delta_{kl}$  where $\delta_{ij}$ is the Kronecker delta. We decompose $D_{ijkl}$ as $D_{ijkl}=\bar D_{ijkl}+\tilde D_{ijkl}$; $\bar D_{ijkl}$ is related to the trace of the strain tensor $tr(\bm \varepsilon)$:
	\begin{equation}
	\bar D_{ijkl}=\lambda\left\{ \left[(1-k)(1-\phi)^2+k \right]H_\varepsilon(tr(\bm\varepsilon))+H_\varepsilon(-tr(\bm\varepsilon))\right\} \delta_{ij}\delta_{kl}
	\end{equation}

and
	\begin{equation}
	\tilde D_{ijkl}=2\mu\left\{ \left[(1-k)(1-\phi)^2+k \right]P_{ijkl}^++P_{ijkl}^-\right\}
	\label{Pijkl} 
	\end{equation}

\noindent where	
	\begin{equation}
	P_{ijkl}^\pm = P_{1ijkl}^\pm+P_{2ijkl}^\pm
	\end{equation}
	\noindent and
	\begin{equation}
	\left\{
	\begin{aligned}
	&P_{1ijkl}^\pm = \sum_{a=1}^3\sum_{b=1}^3 H_\varepsilon(\pm\varepsilon_a)\delta_{ab}n_{ai}n_{aj}n_{bk}n_{bl}\\
	&P_{2ijkl}^\pm = \sum_{a=1}^3\sum_{b\neq a}^3 \frac 1 2 \frac {\langle \varepsilon_a\rangle_\pm - \langle \varepsilon_b\rangle_\pm}{\varepsilon_a-\varepsilon_b}n_{ai}n_{bj}(n_{ak}n_{bl}+n_{bk}n_{al})
	\end{aligned}\right.
	\end{equation}

\noindent with $n_{ai}$  being the $i$-th component of vector $\bm n_a$. To avoid singularity in the calculation when $\varepsilon_a=\varepsilon_b$, we refer to \citep{miehe1993computation} and use a ``perturbation''  for the principal strains:
	\begin{equation}
	  \left\{
	   \begin{aligned}
	&\varepsilon_1 = \varepsilon_1(1+\delta)\hspace{0.5cm} &if\hspace{0.1cm}\varepsilon_1 = \varepsilon_2
	\\ &\varepsilon_3 = \varepsilon_3(1-\delta)\hspace{0.5cm} &if\hspace{0.1cm}\varepsilon_2 = \varepsilon_3
	\end{aligned}\right.
	\end{equation}

We set the perturbation $\delta=1\times 10^{-9}$ for this paper and the second principal strain and volumetric strain are unchanged. 

\subsection{COMSOL implementation}\label{COMSOL implementation}

We implemented our approach into the commercial software COMSOL Multiphysics. Therefore, we establish five modules: Solid Mechanics Module, Darcy Flow Module, Phase Field Module, History-strain Module, and Storage Module. All the established modules are written in strong forms and solved based on the standard finite element discretization in space domain and finite difference discretization in time domain.

The Solid Mechanics and Darcy Flow Modules solve for the displacement and fluid pressure, respectively. The Solid Mechanics Module contains the linear elastic material model and the transient formulation of Darcy's law is used in the Darcy Flow Module. The boundary and initial conditions in Section \ref{Mathematical models of fracture in porous media} are implemented in the Solid Mechanics and Darcy Flow Modules.  The  Phase Field Module and History-strain Module are constructed to solve the other two fields $\phi$ and $H$. We establish the Phase Field Module by revising the Helmholtz equation in COMSOL. The coefficients of the Helmholtz equation have the same form as the governing equation \eqref{governing equation1}. The boundary condition in Eq. \eqref{boundary condition of the phase field} and initial condition \eqref{Initial condition} are also implemented in this module. 

In COMSOL, the ODEs and DAEs Interface can be used to solve distributed ordinary differential equations (ODE) and differential-algebraic equations (DAE). Thus, the History-strain Module is established based on the Distributed ODEs and DAEs Interface. The history-strain field is implemented by establishing a ``previous solution" function in the COMSOL solvers to record and update $H(\bm x, t)$. The format of the equations written into the History-strain Module is shown in the code in ``https://sourceforge.net/projects/phasefieldmodelingcomsol/" where $H_0(\bm x)=0$ is the initial condition of the History-strain Module while pre-existing cracks can be generated by introducing the following initial conditions \citep{borden2012phase}:
	\begin{equation}
	H_0(\bm x) =  \left\{
	   \begin{aligned}
	&\frac {BG_c}{2l_0}\left [1-\frac {2d(\bm x,l)}{l_0}\right ],\hspace{0.1cm}&d(\bm x, l)\leq \frac {l_0} 2
	\\&0,\hspace{0.5cm} & d(\bm x, l)> \frac {l_0} 2
	\end{aligned}\right.
	\end{equation}

\noindent with $B=1\times 10^3$ satisfying Eq. \eqref{governing equation1}. Thus, we succeed in creating an initial $\phi_0=1$ for initial cracks.

\subsection{Staggered scheme}\label{Segregated scheme}

\begin{figure}[htbp]
	\centering
	\includegraphics[width = 12cm]{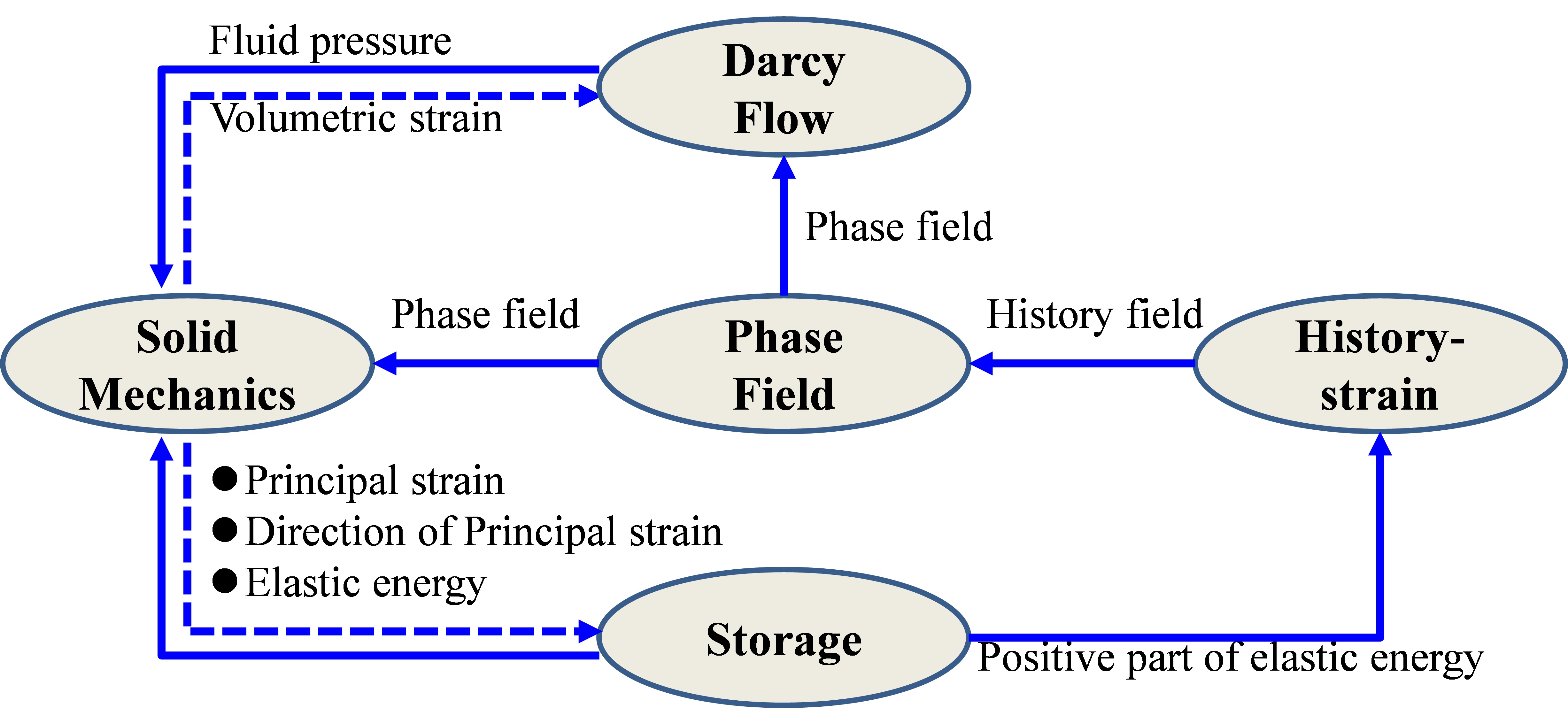}
	\caption{Relationship between the established modules}
	\label{Relationship between the established modules}
	\end{figure}

To facilitate the implementation in COMSOL, a pre-set Storage Module is also used. The relationship between all the established modules is shown in Fig. \ref{Relationship between the established modules}. The Storage Module stores the principal strains as well as the direction of the principal strain from the Solid Mechanics Module during each time step. Some temporary variables such as the elastic energy and component $D_{ijkl}$ are also calculated in the Storage Module. The positive part of the elastic energy $\psi_{\varepsilon}^+$ from the Storage Module is then used to solve and update the local history-strain field in the History-strain Module. Afterwards, the Phase Field Module utilizes the updated history strain $H$ to solve for the phase field. The resulting volumetric strain from the Solid Mechanics Module and the phase field are then employed to compute the pressure field in the Darcy Flow Module.  After assembling the stiffness matrix according to the previously stored intermediate variables and the updated phase field, the Solid Mechanics computes the displacement field by using the updated pressure field.

Figure \ref{Relationship between the established modules} shows the coupling between the established modules within the staggered scheme. It has been shown in \citep{miehe2010phase} that a staggered scheme has advantages over the previously developed monolithic approach \citep{miehe2010thermodynamically} because the latter one does not guarantee convexity of the potential which leads to loss of robustness. However, in recent years, a lot of work has been done on quasi-monolithic (explicit convexification) \citep{heister2015primal} and fully monolithic methods \citep{gerasimov2016line, wick2017modified}. It should be noted that the recently developed monolithic methods are more robust and can be more efficient \citep{gerasimov2016line}. In addition, the monolithic schemes have more accuracy and stronger coupling conditions; therefore, they can achieve a better crack tip velocity \citep{wick2017error}. Although reasonable results are shown for both the monolithic and staggered schemes, only for a monolithic system of phase-field fracture in porous media a mathematical existence result could be obtained \citep{mikelic2015phase}; in staggered schemes there is no final evidence that the numerical solution is valid from a mathematical point of view.

\begin{figure}[htbp]
	\centering
	\includegraphics[width = 15cm]{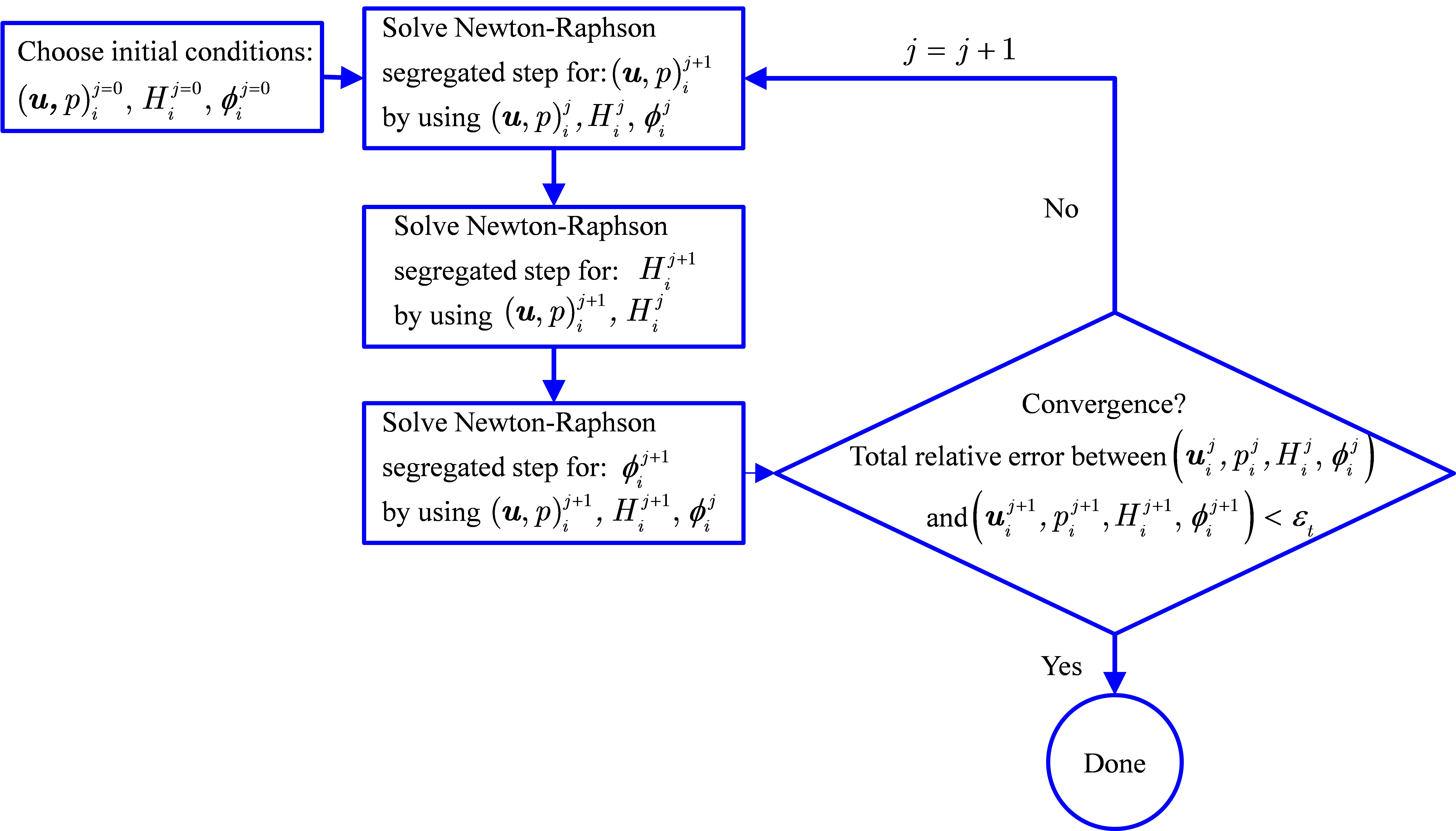}
	\caption{Segregated scheme for the coupled calculation in phase field modeling}
	\label{Segregated scheme for the coupled calculation in phase field modeling}
	\end{figure}

In this paper, we focus on the staggered scheme due to its ease of implementation in COMSOL. More details of the staggered scheme are highlighted in Figure \ref{Segregated scheme for the coupled calculation in phase field modeling}. The displacement and pressure are placed in one staggered step and solved together while the history strain and phase field are in other two staggered steps. The implicit Generalized-$\alpha$ method \citep{borden2012phase} -- ensuring unconditional stability -- is used for time integration. When the time reaches $t_i$, linear extrapolation of the solution in the former time step provides the initial guess for the three staggered steps. Then, the steps are solved sequentially based on the updated results from the previous step. After the solution of all  three fields, the relative error $\varepsilon_r$ is estimated. If $\varepsilon_r$ is less than a prescribed tolerance $\varepsilon_t$, the calculation proceeds from time step $i$ to $t_{i+1}$. Otherwise, a new iteration step will be started until $\varepsilon_r<\varepsilon_t$. In this paper, we set $\varepsilon_t=1\times 10^{-3}$.

\begin{figure}[htbp]
	\centering
	\includegraphics[width = 12cm]{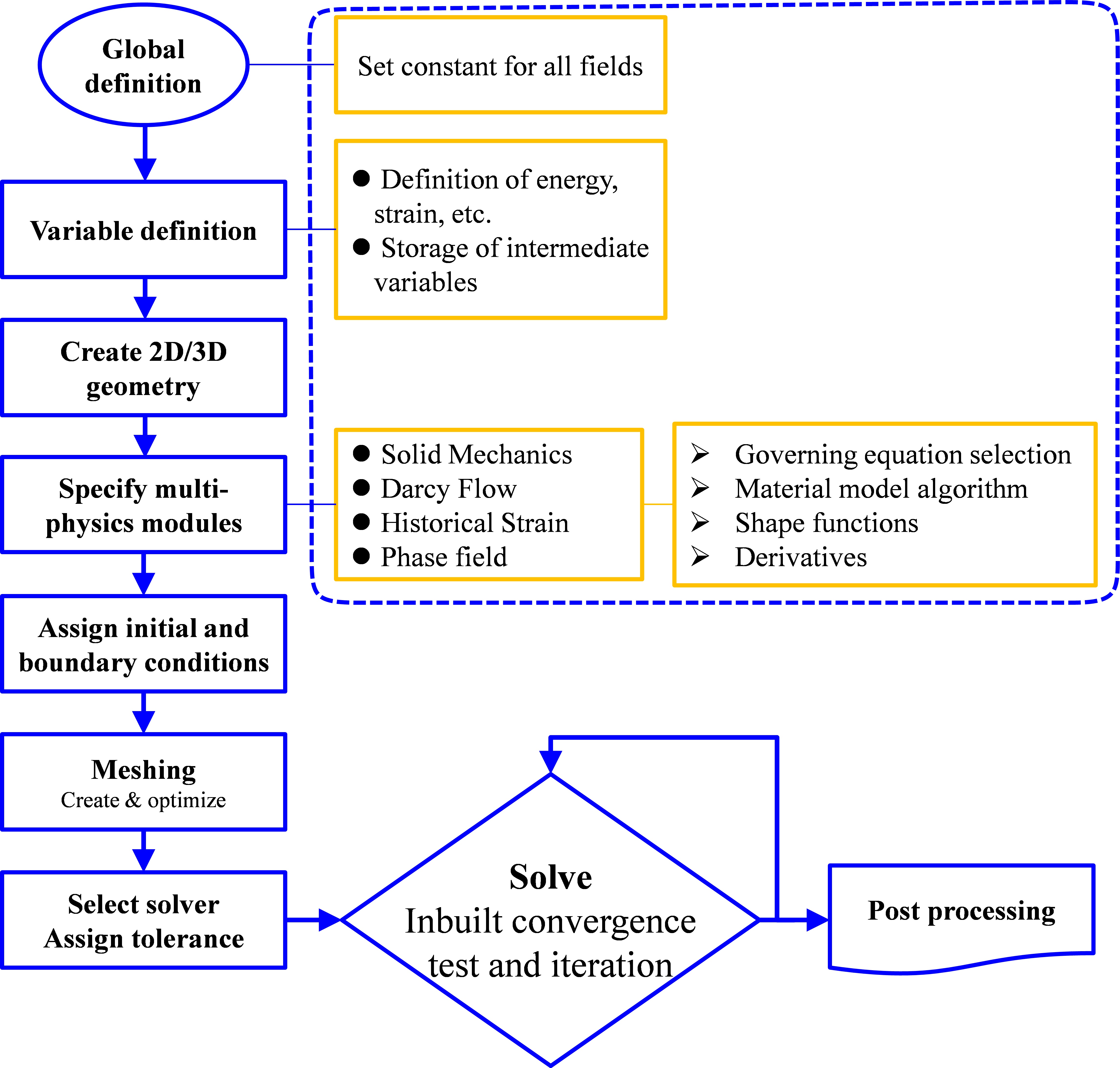}
	\caption{COMSOL implementation of phase field modeling for dynamic fluid-driven crack problems}
	\label{COMSOL implementation of phase field modeling for dynamic fluid-driven crack problems}
	\end{figure}

We take advantage of the Anderson acceleration method \citep{COMSOL2005COMSOL} available in COMSOL to accelerate convergence. The convergence acceleration technology uses the iteration information from the previous Newton iterations and the dimension of the iteration space field is set larger than 300 to control the number of iteration increments. Finally, Fig. \ref{COMSOL implementation of phase field modeling for dynamic fluid-driven crack problems} depicts the implementation of the phase field model for dynamic fluid-driven fracture problems.

\section {Verification of the proposed approach}\label{Verification of the proposed approach}

In this section, three numerical examples are presented to demonstrate the correctness and accuracy of the implemented phase field model. Therefore, the results obtained from the phase field simulations are compared with some existing analytical solutions.

\subsection{Dynamic consolidation}

\begin{figure}[htbp]
	\centering
	\includegraphics[height = 6cm]{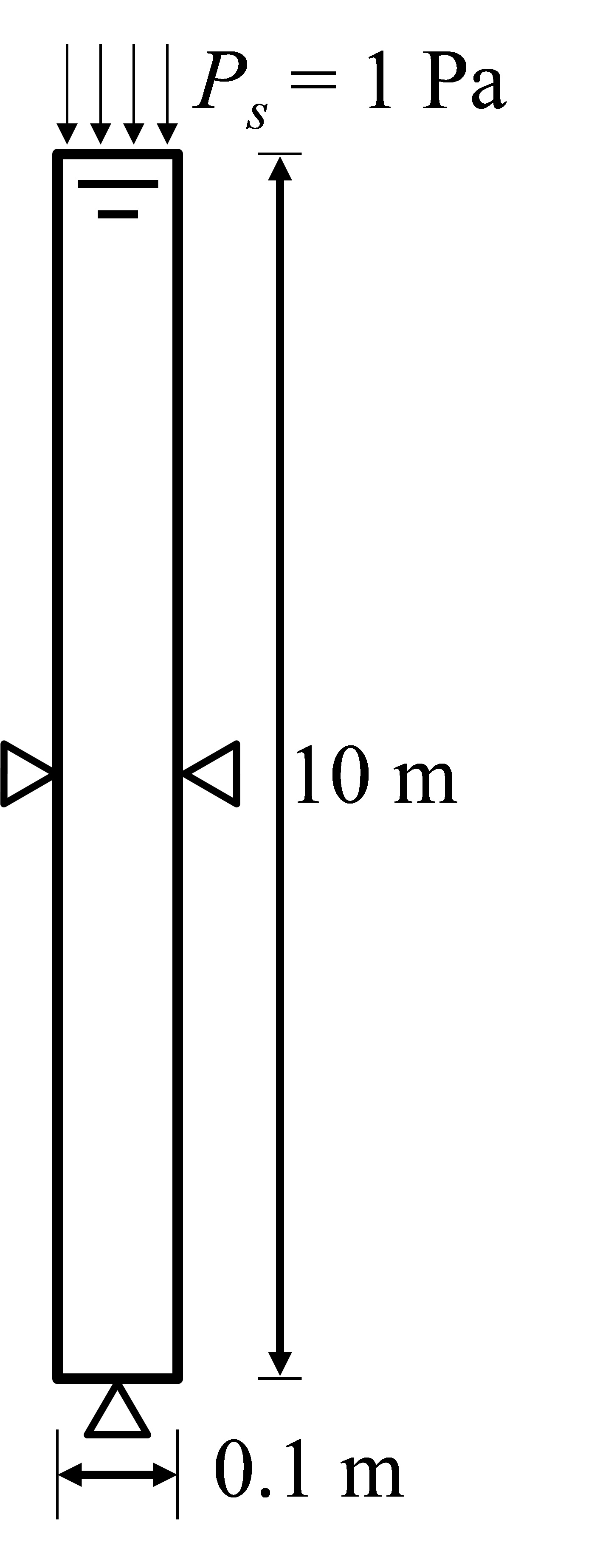}
	\caption{Geometry and boundary conditions of the dynamic consolidation}
	\label{Geometry and boundary conditions of the dynamic consolidation}
	\end{figure}

The first example is presented for dynamic consolidation problem. For 1D dynamic consolidation, \citet{schanz2000transient} showed the transient wave propagation and obtained the analytical solution for the displacement and pressure distribution. To compare with the 1D analytical solution proposed by \citep{schanz2000transient}, a 2D domain is constructed in Fig. \ref{Geometry and boundary conditions of the dynamic consolidation} along with the corresponding boundary conditions. The left and right boundaries as well as the bottom of the domain are constrained in the normal direction. The top boundary of the domain is permeable with $p=0$ and subjected to a sudden pressure $P_s=$ 1 kPa when $t = 0$. Q4 elements with size $h=0.1$ m are used to discretize all the fields and the parameters used for calculation are listed in Table \ref{Calculation parameters for the example of dynamic consolidation}. In addition, a mandatory condition of $\phi = 0$ is set for the porous domain to compare well with the 1D analytical solution. 

	\begin{table}[htbp]
	\small
	\caption{Calculation parameters for the example of dynamic consolidation}
	\label{Calculation parameters for the example of dynamic consolidation}
	\centering
	\begin{tabular}{llllll}
	\hline
	$E$ & 0.254 GPa & $\nu$ &0.3 & $c_1$ & 0.5 \\
	 $c_2$ & 1.0 & $\varepsilon_{pR}$ & 0.48 & $\rho_{R}$,  $\rho_{F}$& $1.0\times10^{3}$ kg/m$^3$  \\
	 $\alpha_R$ &0.981 & $q_R$ & 0 & $q_F$ & 0 \\ 
	$K_R$ & $3.55\times10^{-12}$ m$^2$ & $K_F$ & $3.55\times10^{-12}$ m$^2$ & $c_R$ & $3.33\times10^{-10}$ 1/Pa \\
	 $c_F$ & $3.33\times10^{-10}$ 1/Pa & $\mu_R$ & $1\times10^{-3}$ Pa$\cdot$s & $\mu_F$ & $1\times10^{-3}$ Pa$\cdot$s\\
	\hline
	\normalsize
	\end{tabular}
	\end{table}

\begin{figure}[htbp]
	\centering
	\includegraphics[width= 10cm]{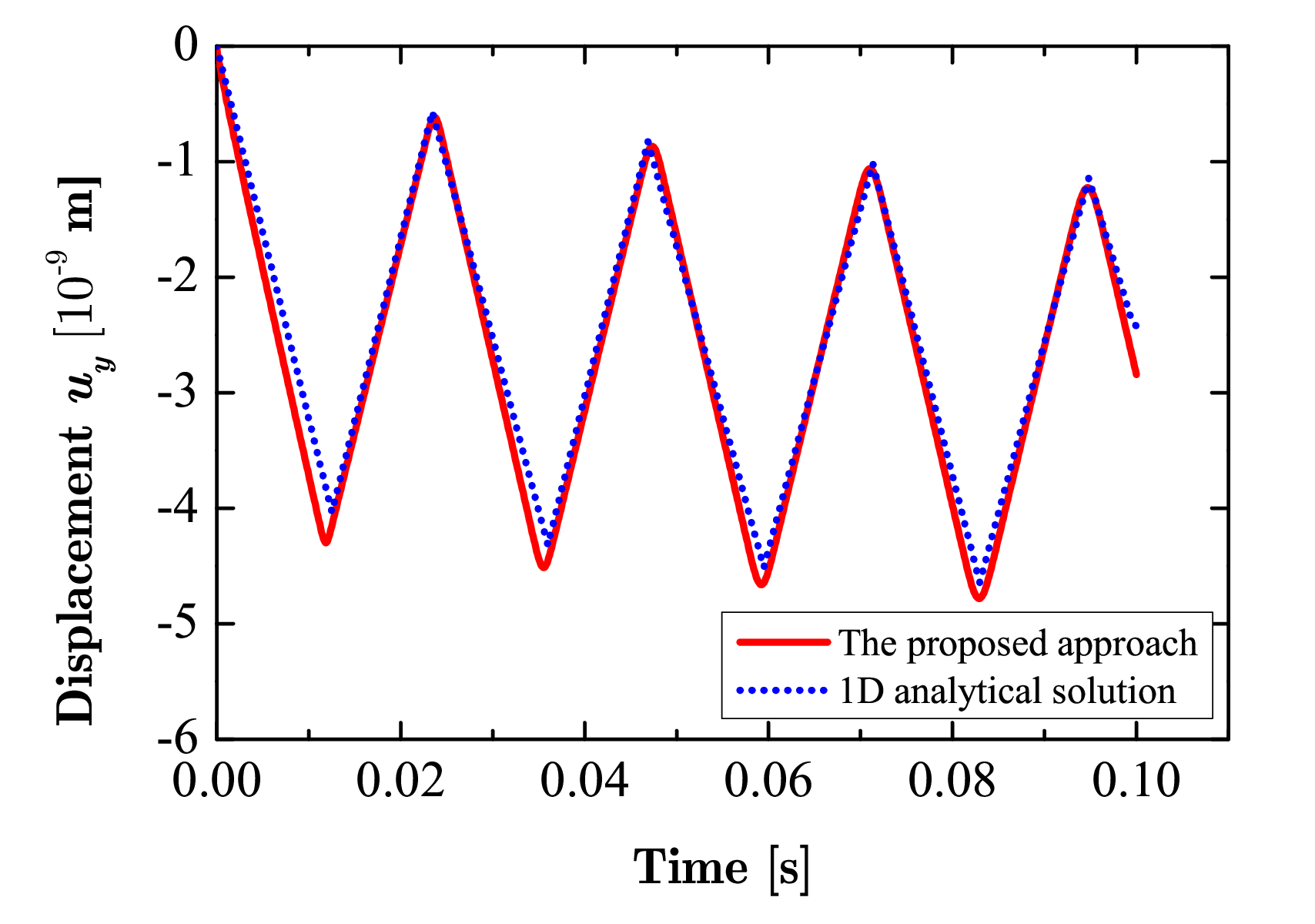}
	\caption{Vertical displacement on the top edge for the example of the dynamic consolidation}
	\label{Vertical displacement on the top edge for the example of the dynamic consolidation}
	\end{figure}

\begin{figure}[htbp]
	\centering
	\includegraphics[width= 10cm]{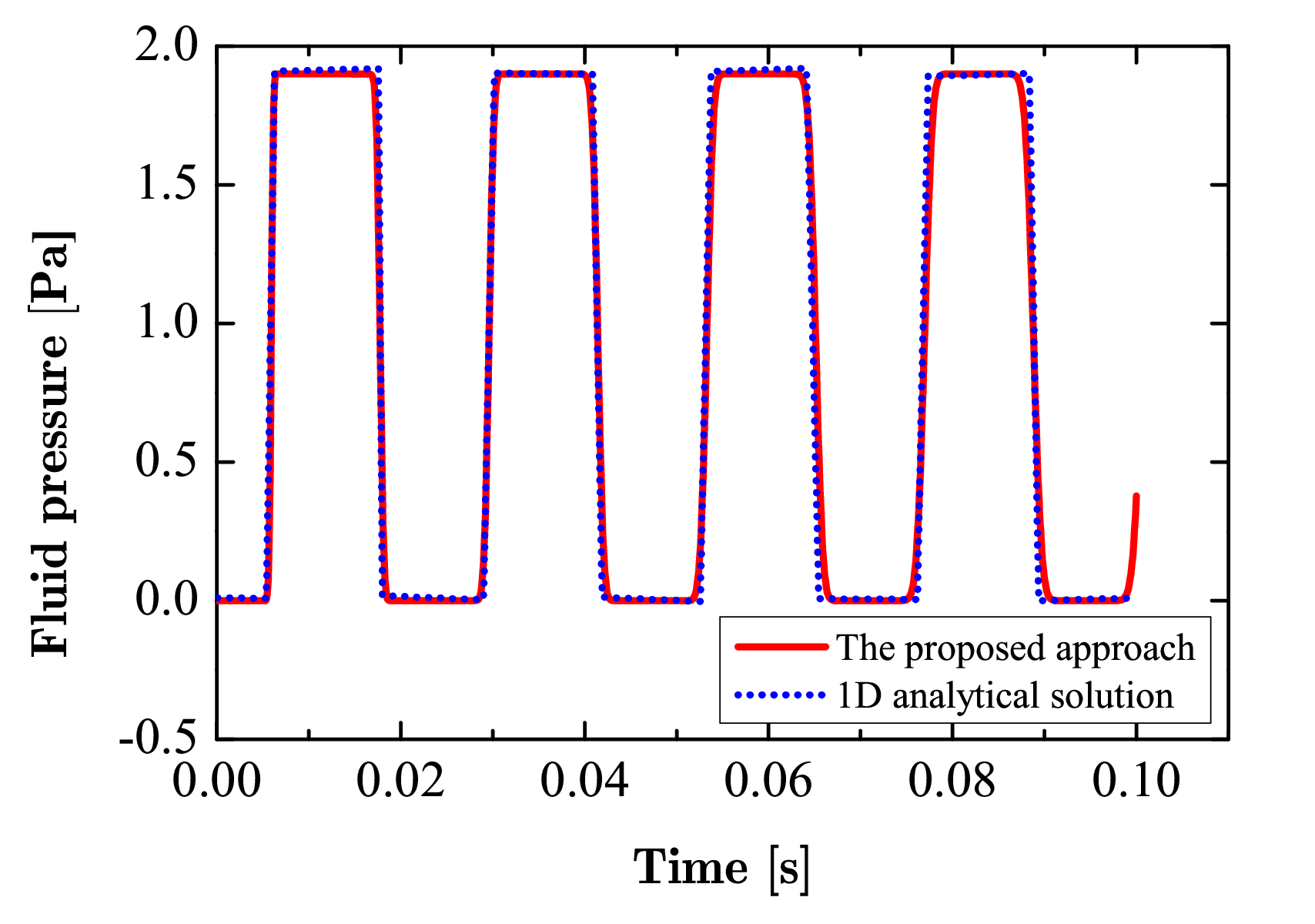}
	\caption{Fluid pressure on the bottom for the example of the dynamic consolidation}
	\label{Fluid pressure on the bottom for the example of the dynamic consolidation}
	\end{figure}

We set the time step as $\Delta t=1\times 10^{-5}$ s and comparison of the vertical displacement on the top edge of the domain by using the proposed approach and the 1D analytical solution is shown in Fig. \ref{Vertical displacement on the top edge for the example of the dynamic consolidation}. For dynamic consolidation, the transient wave from sudden load causes oscillation in the displacement as shown in Fig. \ref{Vertical displacement on the top edge for the example of the dynamic consolidation}. The fluid pressure on the bottom of the domain is shown in Fig. \ref{Fluid pressure on the bottom for the example of the dynamic consolidation}. The fluid pressure varies between 0 and 1.8 Pa (close to twice the external load). The reason is the wave reflection on the bottom of the domain. The displacement and fluid pressure obtained by the proposed approach are in good agreement with the results by using the 1D analytical solution, thereby initially indicating the reliability and feasibility of the proposed approach.

\subsection{Pressure distribution in a single crack}\label{Pressure distribution in a single fracture}

\begin{figure}[htbp]
	\centering
	\includegraphics[width = 12cm]{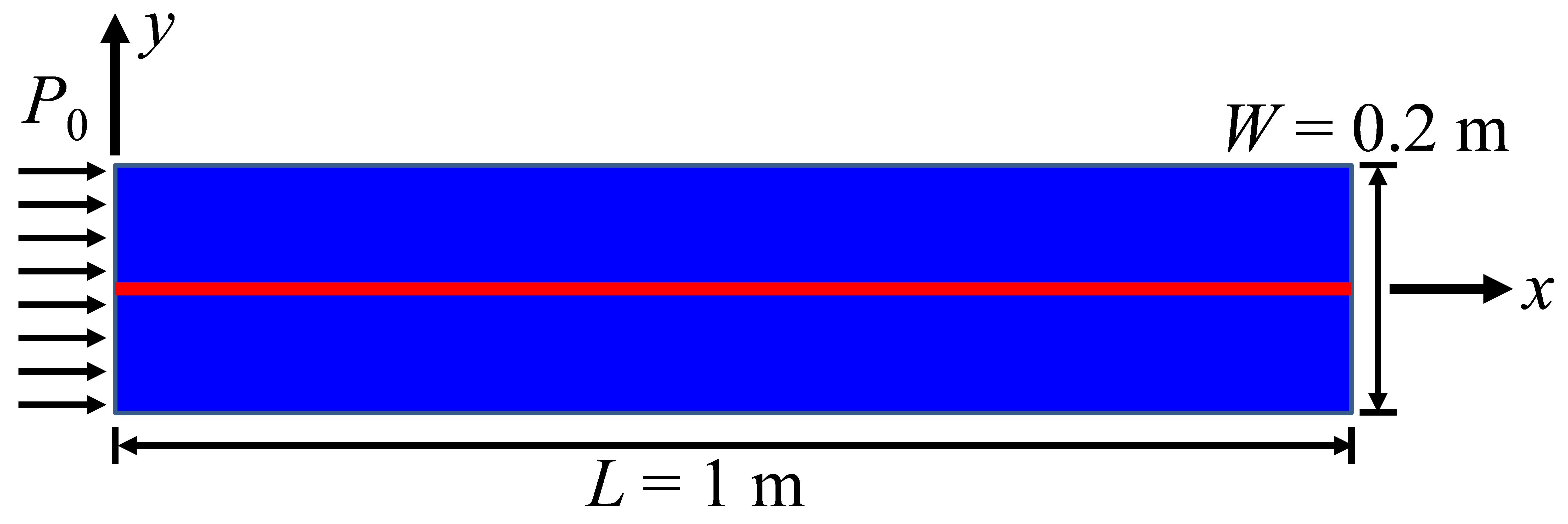}
	\caption{Geometry and boundary conditions of the rectangular plate with a pre-existing crack subjected to sudden fluid pressure}
	\label{Geometry and boundary conditions of the rectangular plate with a pre-existing crack subjected to sudden fluid pressure}
	\end{figure}

The second example is a rectangular plate with a pre-existing crack subjected to sudden fluid pressure $P_0$. We test this example to obtain the pressure distribution along the crack. This example has been tested by \citet{yang2017hydraulic} by using an enriched numerical manifold method. The geometry and boundary conditions are shown in Fig. \ref{Geometry and boundary conditions of the rectangular plate with a pre-existing crack subjected to sudden fluid pressure}. For the displacement field, the right end of the plate is fixed while the left end is subjected to pressure $P_0$. For the flow field, all the rest boundaries of the plate are impermeable except the pressure boundary on the left end. 

This problem is a well-known unsteady-state problem and has an analytical solution \citep{yang2017hydraulic} for the pressure distribution along the existing crack:
	\begin{equation}
	\frac{P}{P_0}=1+\frac 4 {\pi}\sum_{n=0}^{\infty}\left\{\mathrm{exp}\left[-(2n+1)^2(T_d/4)\pi^2\right]\mathrm{cos}\left[\frac{(2n+1)\pi}{2}\xi\right]\left[\frac{(-1)^{n+1}}{2n+1}\right]\right\}
	\end{equation}

\noindent where $\xi = (L-x)/x$, $P$ is the fluid pressure at the point with a distance of $x$ from the left end of the plate (i.e., coordinate $x$ in Fig. \ref{Geometry and boundary conditions of the rectangular plate with a pre-existing crack subjected to sudden fluid pressure}), $L$ is the length of the plate, and $T_d$ is a dimensionless time given by
	\begin{equation}
	T_d=K_w\frac{a t}{12\mu_F L^2}
	\end{equation}

\noindent with $K_w=1/c_F$ bulk modulus of water and $a$ aperture of the existing crack.

To compare well the analytical solution, we adopt the calculation parameters as listed in Table \ref{Parameters for the rectangular plate with a pre-existing crack subjected to sudden fluid pressure} and we neglect the mutual coupling between the Solid Mechanics Module and Darcy Flow Module. The aperture $a=3\times10^{-5}$ m is used to match the parameters in Table \ref{Parameters for the rectangular plate with a pre-existing crack subjected to sudden fluid pressure}. Density of the material is $\rho=2700$ kg/m$^3$. By exploiting symmetry, we take half of the plate to calculate the numerical results. Q4 elements with size $h=5\times10^{-3}$ m is used to discretize all the fields and the time step is chosen as $\Delta t = 1\times10^{-6}$ s for calculation.

\begin{figure}[htbp]
	\centering
	\includegraphics[width = 12cm]{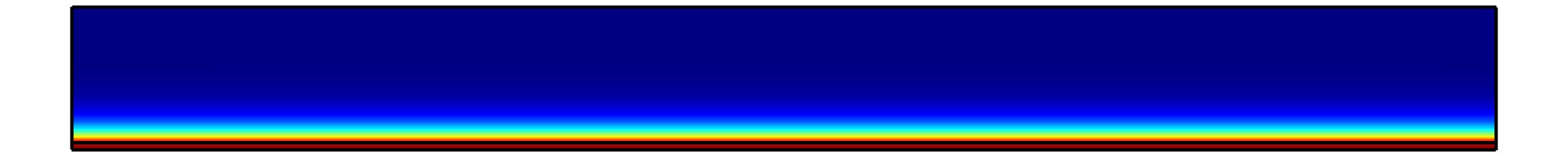}
	\caption{Induced phase field for the pre-existing crack}
	\label{Induced phase field for the pre-existing crack}
	\end{figure}

\begin{figure}[htbp]
	\centering
	\subfigure[$T_d=0.1$]{\includegraphics[width = 12cm]{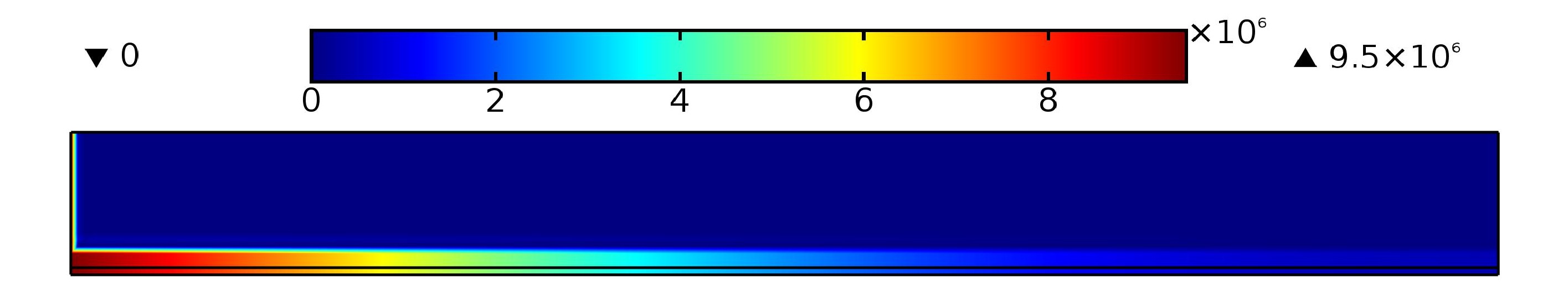}}\\
	\subfigure[$T_d=0.2$]{\includegraphics[width = 12cm]{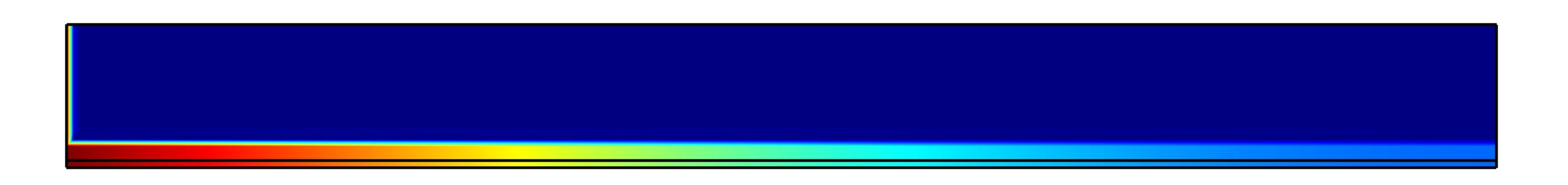}}\\
	\subfigure[$T_d=0.3$]{\includegraphics[width = 12cm]{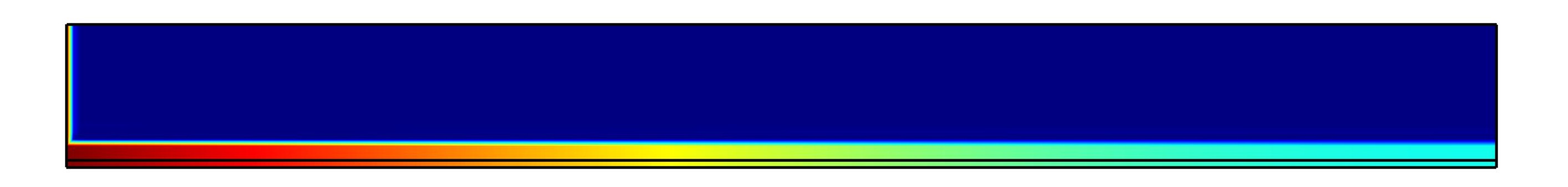}}\\
	\subfigure[$T_d=0.5$]{\includegraphics[width = 12cm]{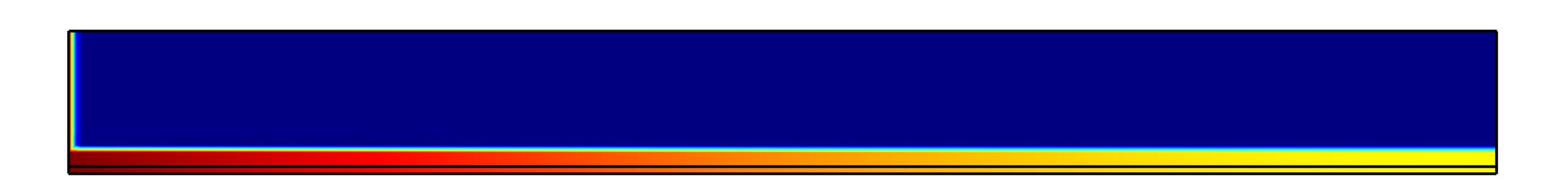}}
	\caption{Pressure distribution of the plate at different time}
	\label{Pressure distribution of the plate at different time}
	\end{figure}

We use the initial history strain field to create the pre-existing crack. The resulting pre-existing crack is shown in Fig. \ref{Induced phase field for the pre-existing crack}, while Fig. \ref{Pressure distribution of the plate at different time} illustrates the pressure distribution along the pre-existing crack at different time $T_d$. As observed, the pressure extends along the crack and towards the right end of the plate as the time increases. Meanwhile, the fluid flow is much easier to penetrate along the crack than perpendicular to the crack.

	\begin{table}[htbp]
	\small
	\caption{Parameters for the rectangular plate with a pre-existing crack subjected to sudden fluid pressure}
	\label{Parameters for the rectangular plate with a pre-existing crack subjected to sudden fluid pressure}
	\centering
	\begin{tabular}{llllll}
	\hline
	$E$ & 55.8 GPa & $\nu$ &0.25 & $G_c$ &$1\times10^{4}$ N/m\\
	$k$ &$1\times10^{-9}$ & $l_0$ & $1\times10^{-2}$ m & $c_1$ & 0.5 \\
	 $c_2$ & 1.0 & $\varepsilon_{pR}$ & $2\times10^{-5}$ & $\rho_{R}$,  $\rho_{F}$& $1.0\times10^{3}$ kg/m$^3$  \\
	 $\alpha_R$ &$2\times10^{-5}$ & $q_R$ & 0 & $q_F$ & 0 \\ 
	$K_R$ & $1\times10^{-20}$ m$^2$ & $K_F$ & $7.5\times10^{-11}$ m$^2$ & $c_R$ & $4.55\times10^{-10}$ 1/Pa \\
	 $c_F$ & $4.55\times10^{-10}$ 1/Pa & $\mu_R$ & $1\times10^{-3}$ Pa$\cdot$s & $\mu_F$ & $1\times10^{-3}$ Pa$\cdot$s\\
	\hline
	\normalsize
	\end{tabular}
	\end{table}

\begin{figure}[htbp]
	\centering
	\subfigure[$T_d=0.1$]{\includegraphics[width = 7cm]{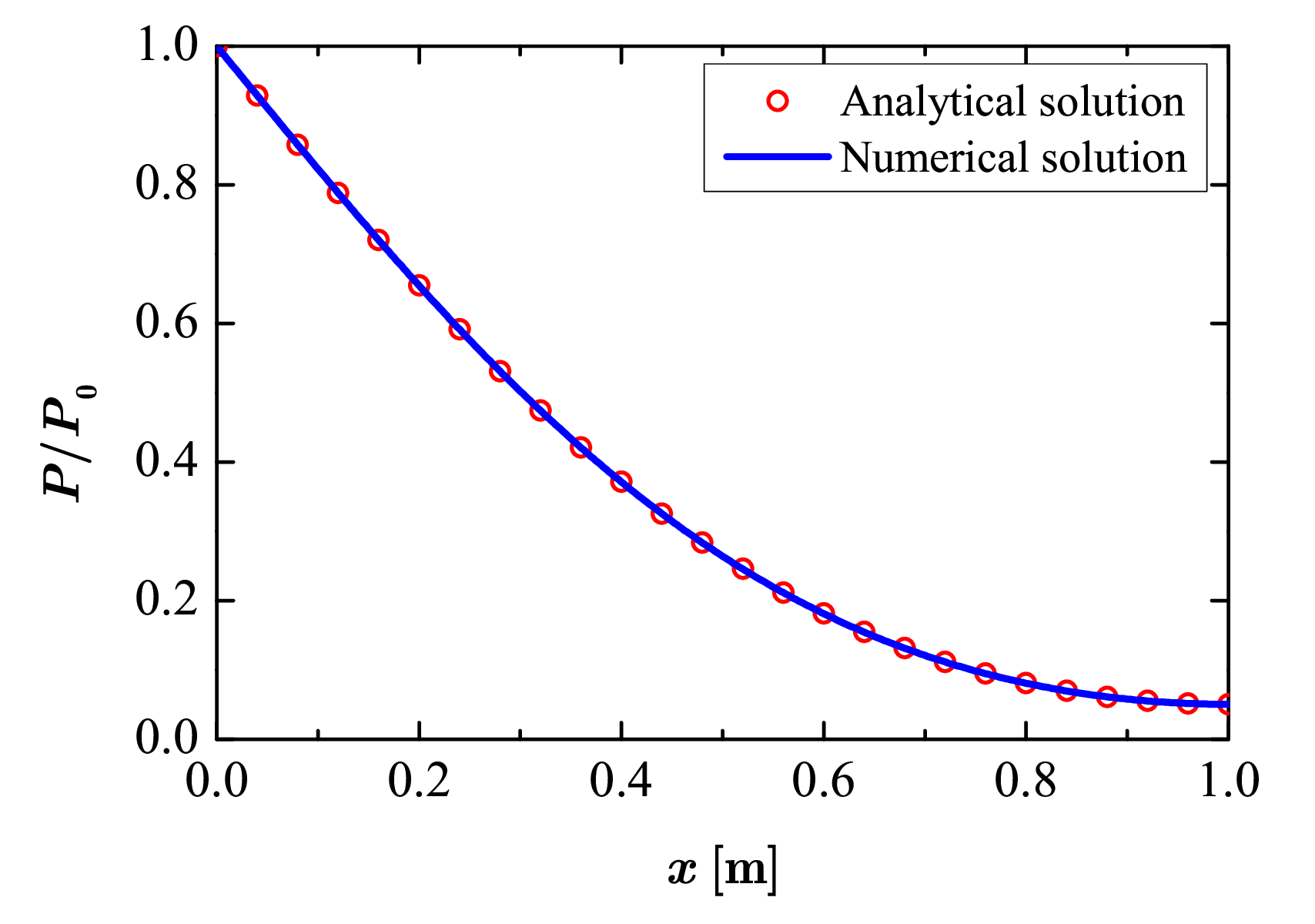}}
	\subfigure[$T_d=0.2$]{\includegraphics[width = 7cm]{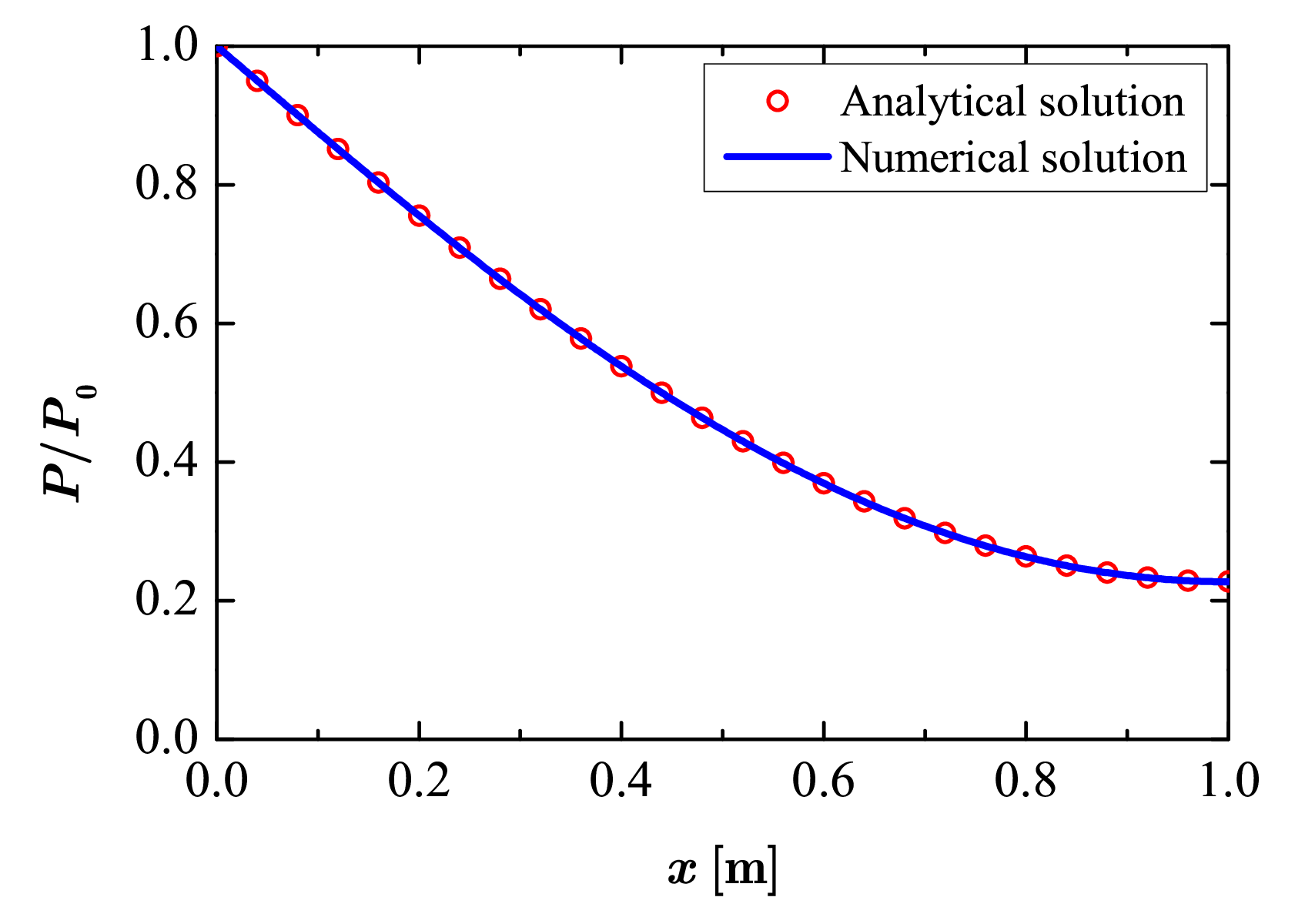}}\\
	\subfigure[$T_d=0.3$]{\includegraphics[width = 7cm]{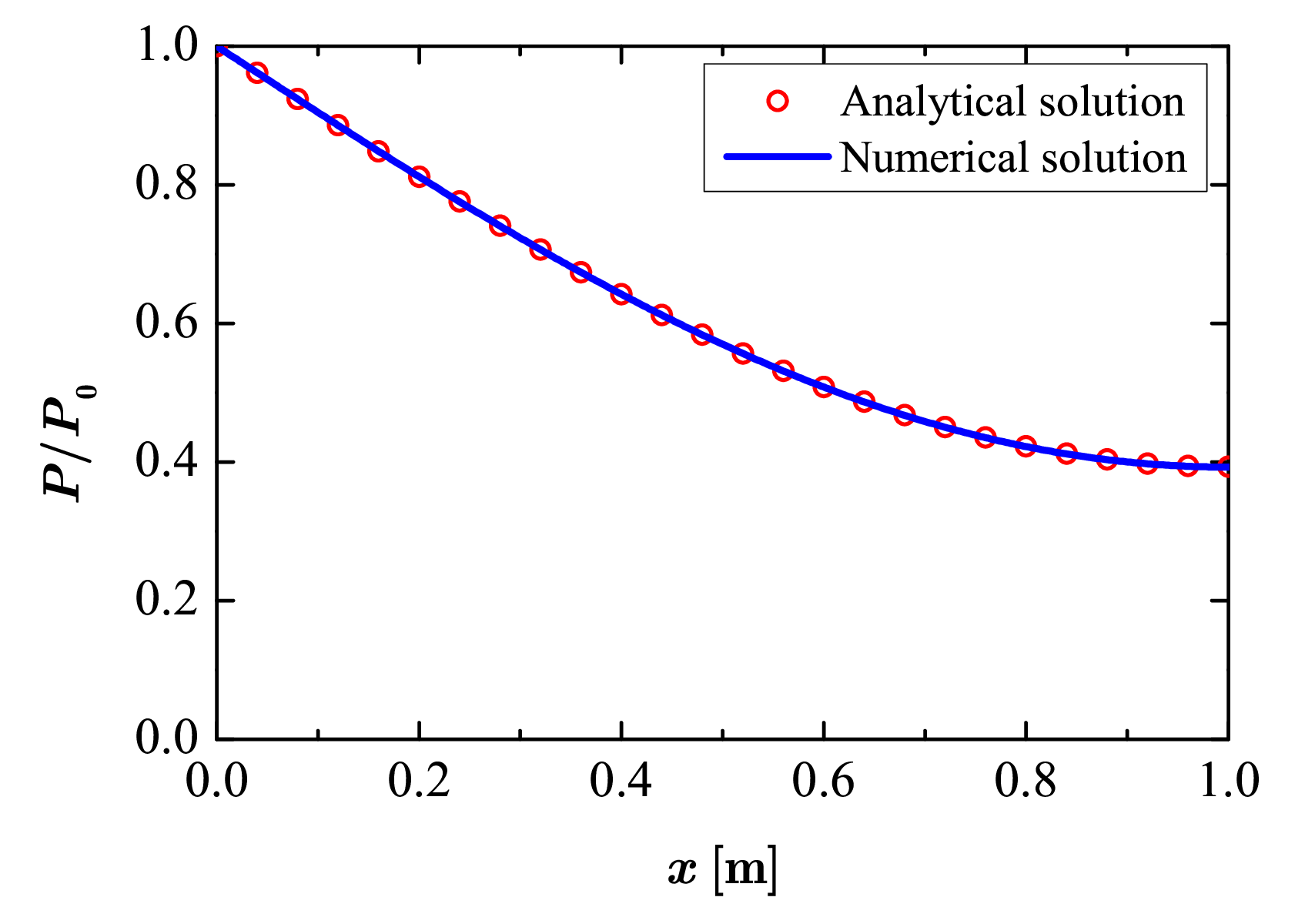}}
	\subfigure[$T_d=0.5$]{\includegraphics[width = 7cm]{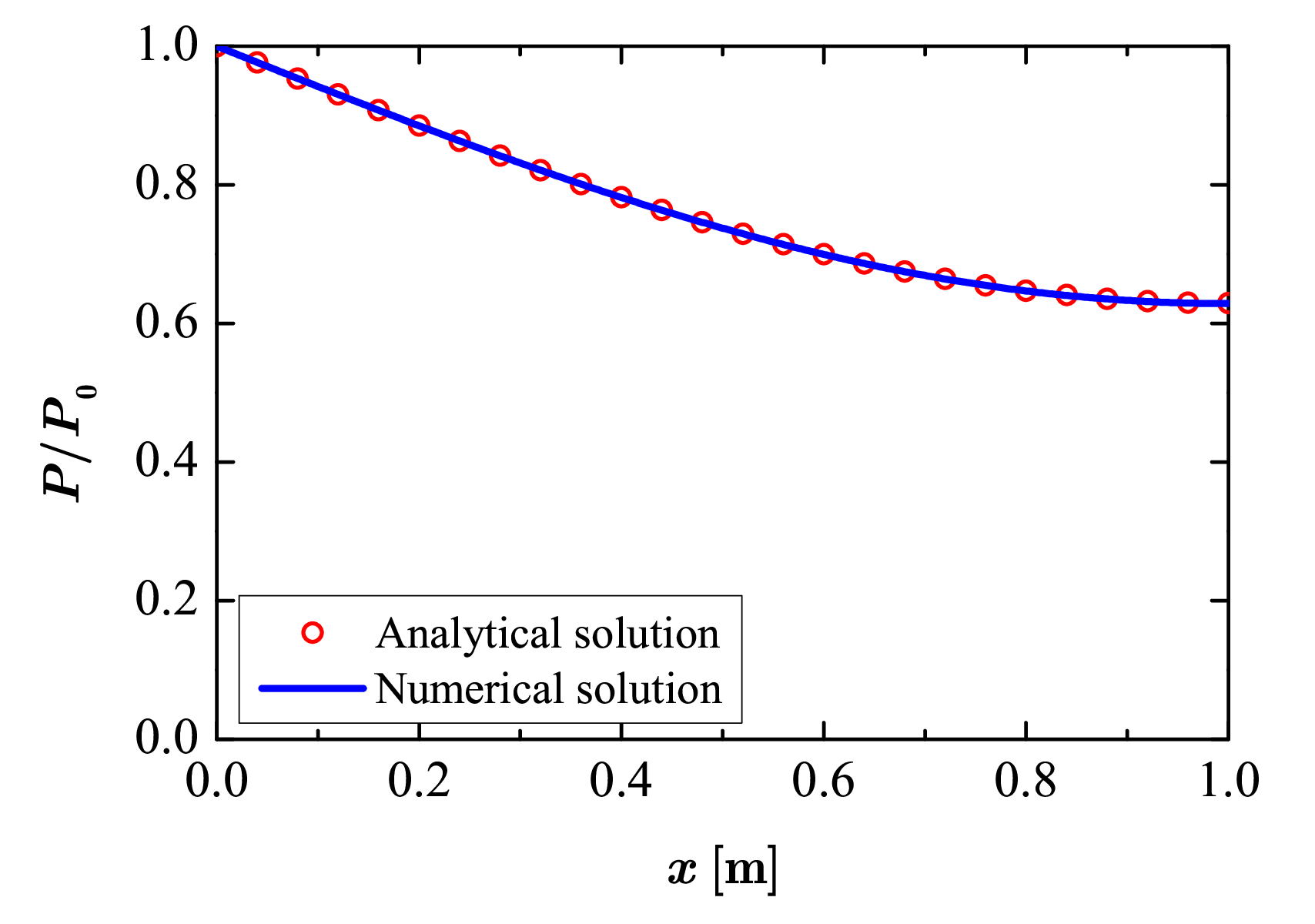}}
	\caption{Comparison of the numerical and analytical pressure distribution along the crack}
	\label{Comparison of the numerical and analytical pressure distribution along the crack}
	\end{figure}

\begin{figure}[htbp]
	\centering
	\includegraphics[width = 7cm]{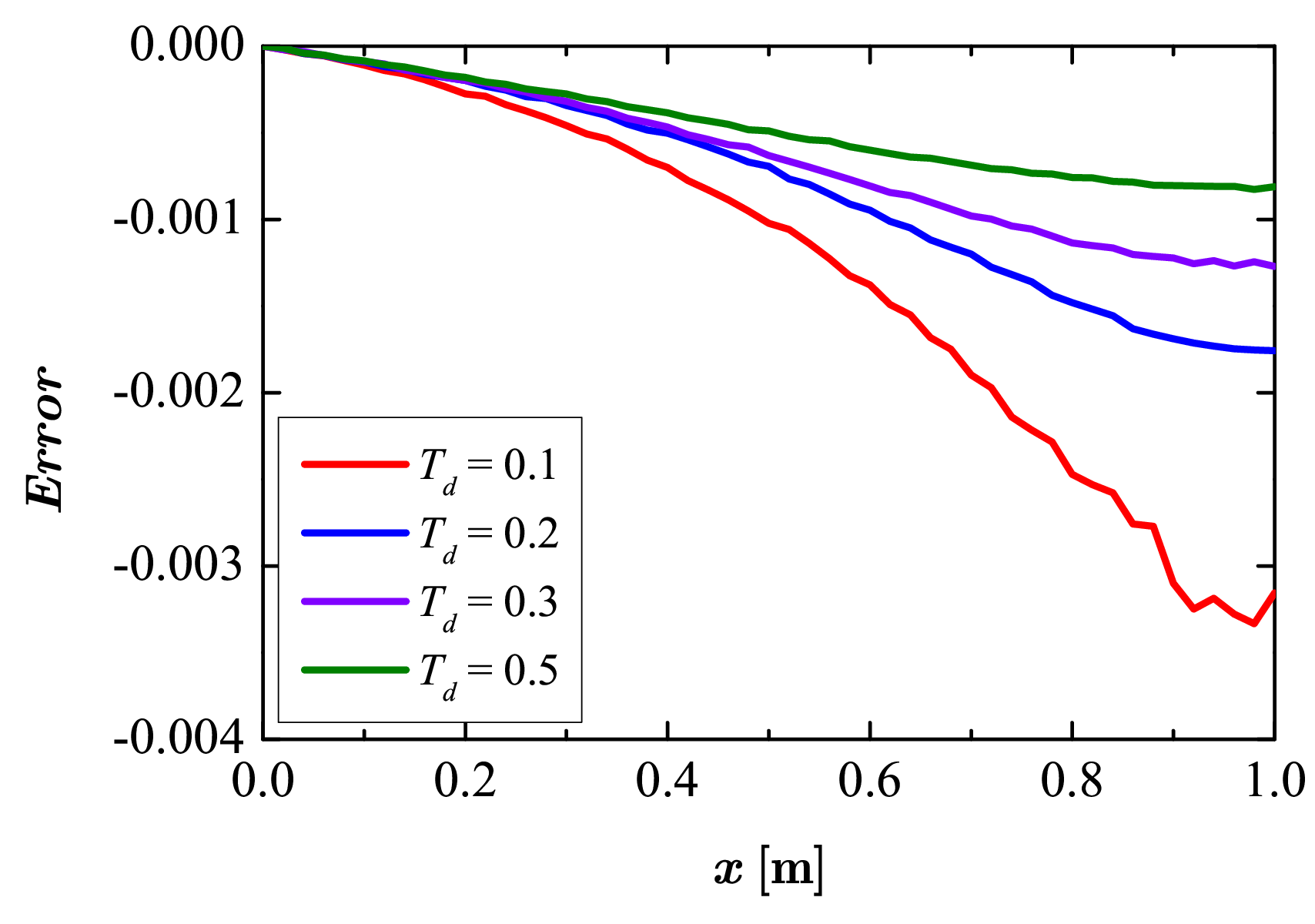}
	\caption{Error of the numerical pressure along the crack compared with the analytical solution}
	\label{Error of the numerical pressure along the crack compared with the analytical solution}
	\end{figure}

Furthermore, comparison of the hydraulic pressure distribution at different time obtained from the phase field modeling and from the analytical solution is shown in Fig. \ref{Comparison of the numerical and analytical pressure distribution along the crack}. The numerical results by the phase field modeling are in good agreement with the analytical solution, showing the accuracy of the phase field method. In addition, as observed, the fluid pressure gradually decreases along the pre-existing crack, while the pressure at the same point gradually increases with time. Figure \ref{Comparison of the numerical and analytical pressure distribution along the crack} also shows that the pressure in the cracked region tends to $P_0$ with time. Figure \ref{Error of the numerical pressure along the crack compared with the analytical solution} shows the error of the numerical pressure along the crack compared with the analytical solution. The error is quite small and less than 0.004. The right end of the crack has the maximum error while the error increases as time increases. The error analysis also reflects the consistency of the phase field modeling and the analytical solution.

\subsection{Pressure distribution in a specimen with two sets of joints}

\begin{figure}[htbp]
	\centering
	\includegraphics[width = 8cm]{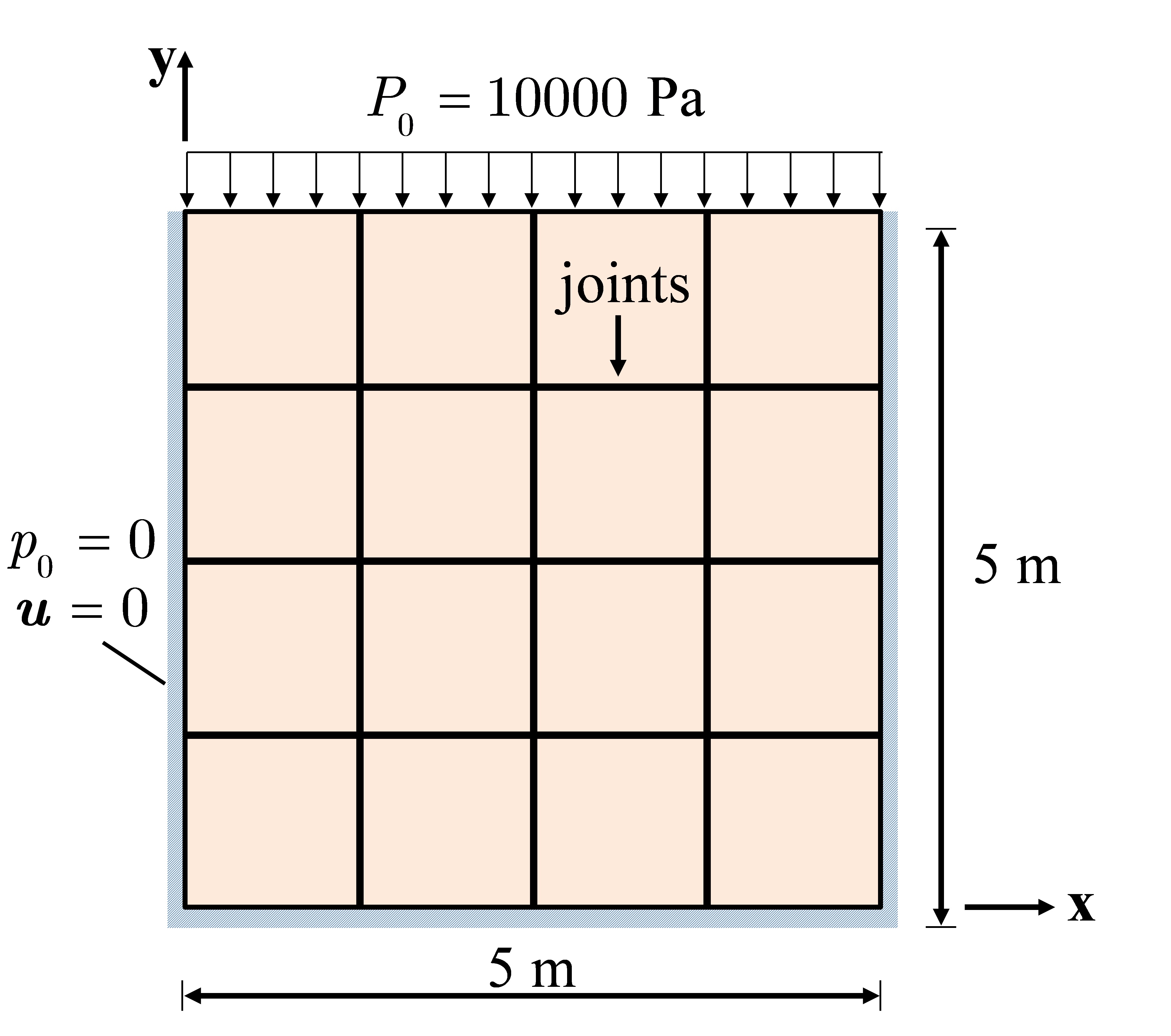}
	\caption{Geometry and boundary conditions of the specimen with two sets of joints}
	\label{Geometry and boundary conditions of the specimen with two sets of joints}
	\end{figure}

The third example tests the pressure distribution in a specimen subjected to hydraulic pressure $P_0$ = 10000 Pa. The specimen has two sets of perpendicular joints. The geometry and boundary conditions are depicted in Fig. \ref{Geometry and boundary conditions of the specimen with two sets of joints}. The density of the specimen is 2700 kg/m$^3$. For the pressure field, all boundaries are impermeable except the pressure boundary $P_0$ on the top. Similar examples have been studied by \citet{jiao2015two} with discontinuous deformation analysis (DDA) and \citet{yang2017hydraulic} with numerical manifold method (NMM). In this paper, all the parameters for calculation are listed in Table \ref{Parameters for the specimen with two sets of joints}. Besides, gravity acceleration $g=10$ m/s$^2$ is adopted.
 
	\begin{table}[htbp]
	\small
	\caption{Parameters for the specimen with two sets of joints}
	\label{Parameters for the specimen with two sets of joints}
	\centering
	\begin{tabular}{llllll}
	\hline
	$E$ & 55.8 GPa & $\nu$ &0.25 & $G_c$ &$1\times10^{-2}$ N/m\\
	$k$ &$1\times10^{-9}$ & $l_0$ & $2.5\times10^{-2}$ m & $c_1$ & 0.5 \\
	 $c_2$ & 1.0 & $\varepsilon_{pR}$ & $1\times10^{-5}$ & $\rho_{R}$,  $\rho_{F}$& $1.0\times10^{3}$ kg/m$^3$  \\
	 $\alpha_R$ &$1\times10^{-5}$ & $q_R$ & 0 & $q_F$ & 0 \\ 
	$K_R$ & $1\times10^{-20}$ m$^2$ & $K_F$ & $2.08\times10^{-8}$ m$^2$ & $c_R$ & $4.55\times10^{-10}$ 1/Pa \\
	 $c_F$ & $4.55\times10^{-10}$ 1/Pa & $\mu_R$ & $1\times10^{-3}$ Pa$\cdot$s & $\mu_F$ & $1\times10^{-3}$ Pa$\cdot$s\\
	\hline
	\normalsize
	\end{tabular}
	\end{table}

\begin{figure}[htbp]
	\centering
	\includegraphics[width = 6cm]{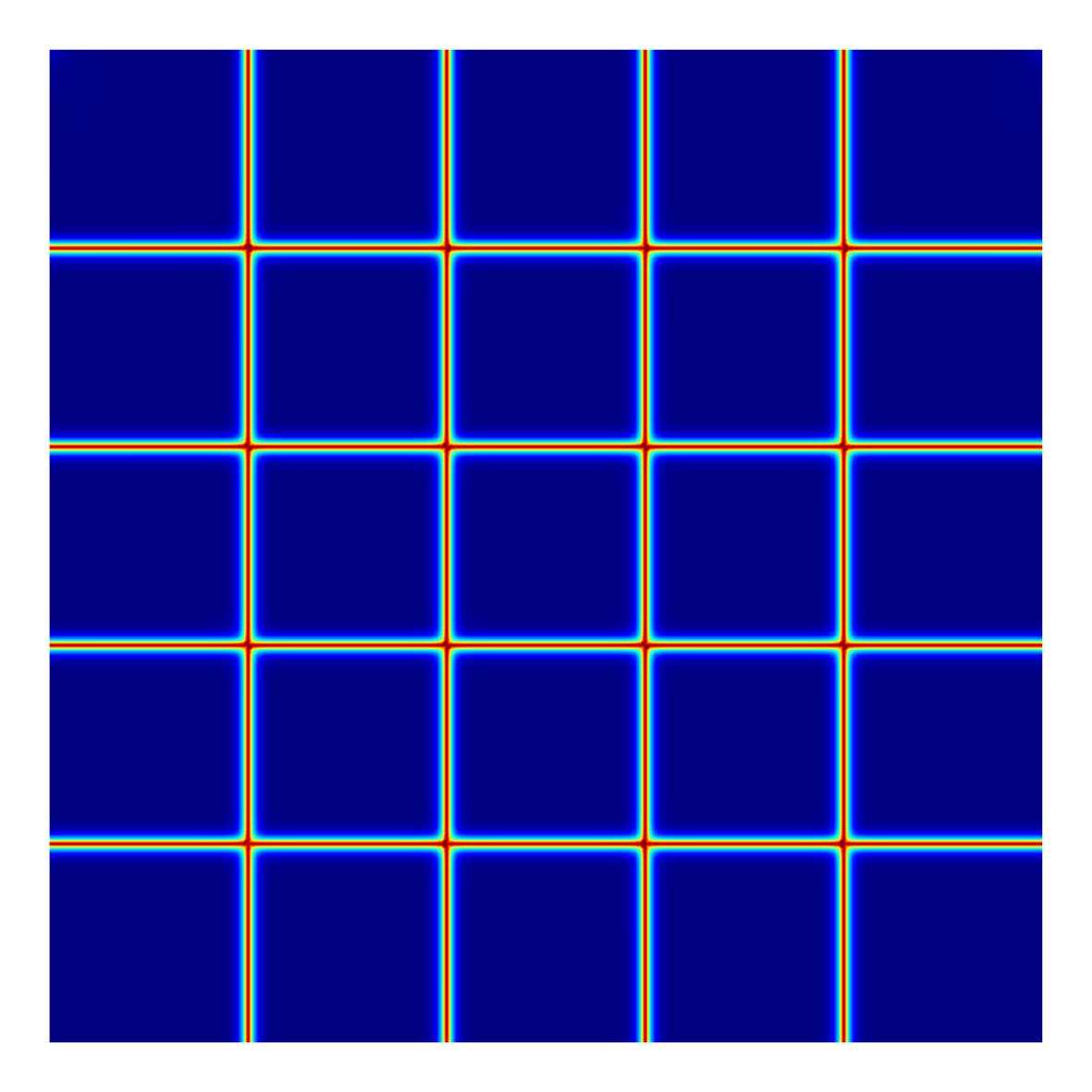}
	\caption{Induced phase field for the two sets of joints}
	\label{Induced phase field for the two sets of joints}
	\end{figure}

\begin{figure}[htbp]
	\centering
	\includegraphics[width = 7cm]{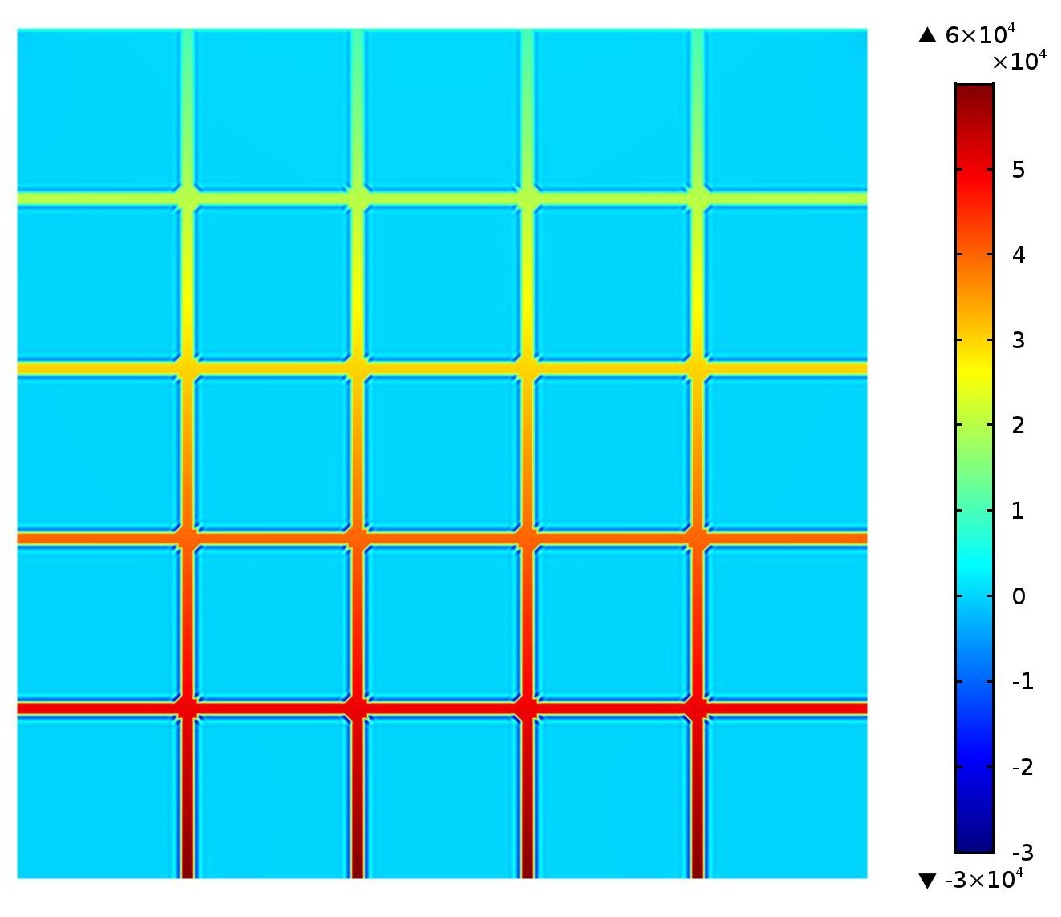}
	\caption{Pressure distribution for the specimen with the two sets of joints at $5\times10^{-3}$ s (unit: Pa)}
	\label{Pressure distribution for the specimen with the two sets of joints}
	\end{figure}

We discretize the domain with uniform Q4 elements of $h=2.5\times10^{-2}$ m and the joints are created by the mandatory condition $\phi=1$. The induced phase field for the two sets of joints is shown in Fig. \ref{Induced phase field for the two sets of joints} and the pressure distribution by using the proposed approach is shown in Fig. \ref{Pressure distribution for the specimen with the two sets of joints}. As expected, the pressure along the joints has a linear relationship with the $y$ coordinate. Because of zero fluid flow and fluid pressure in Fig. \ref{Geometry and boundary conditions of the specimen with two sets of joints}, the analytical solution for the water head along the two sets of joints can be obtained easily. The analytical solution of water header at all the nodes of the joint network is $H_w=p/\rho g + y=6$ m. For the numerical simulation, the water head at all the nodes is still 6 m and in good agreement with analytical solution, indicating the fluid flow algorithm coupled with the phase field method is feasible and correct. In summary,  the numerical results by the proposed approach all match well the analytical solution for the presented three examples in this section, thereby showing feasibility and rationality of the proposed approach.

\section{Examples of dynamic fluid-driven crack branching}\label{Examples of dynamic fluid-driven crack branching}

In this section,  examples of porous medium subjected to internal fluid injection are presented to demonstrate capability of the proposed modeling approach. The crack branching caused by hydraulic fracturing is clearly observed in these examples. 

\subsection{2D examples}\label{2D examples}

\begin{figure}[htbp]
	\centering
	\includegraphics[height = 6cm]{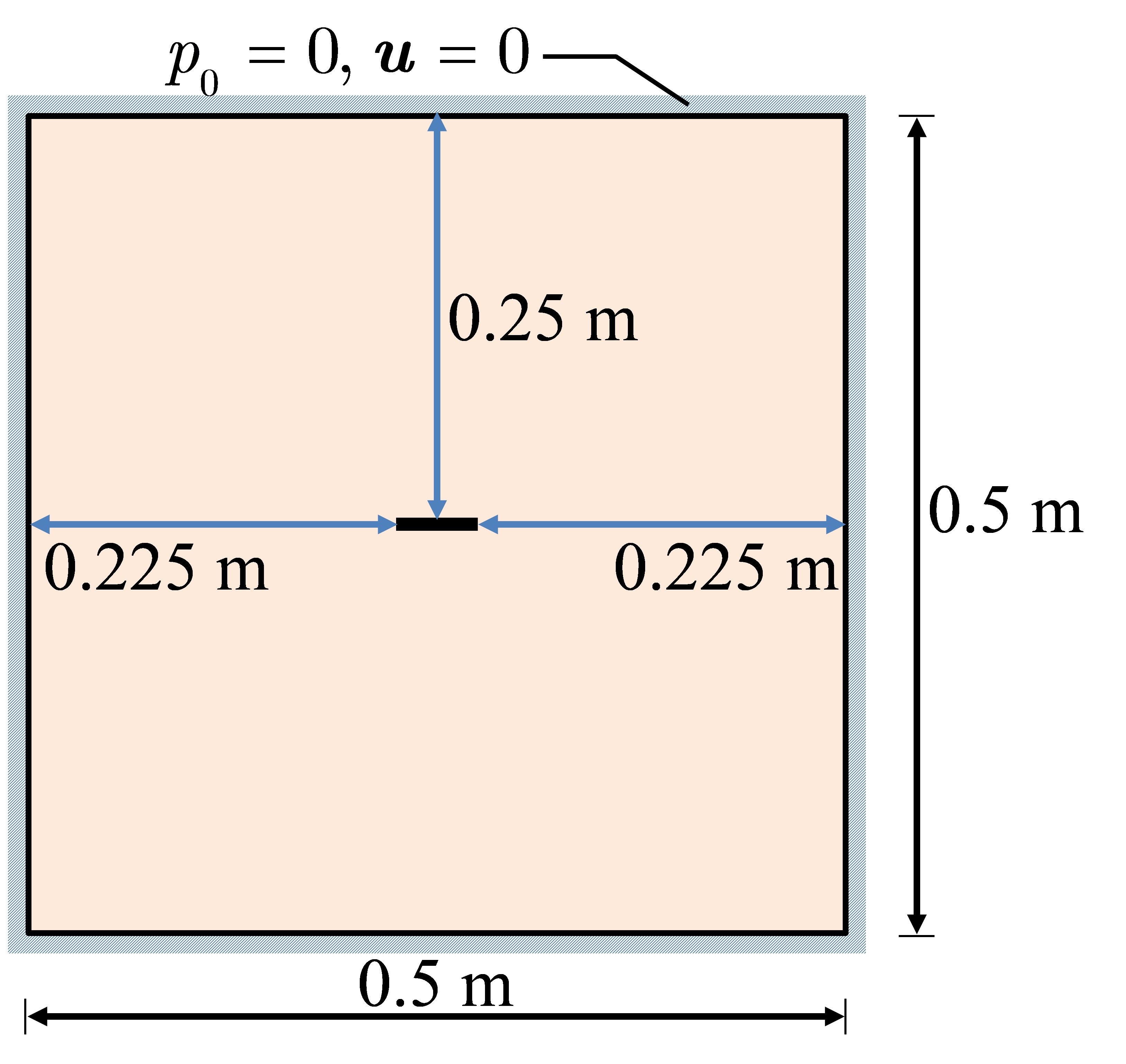}
	\caption{Geometry and boundary conditions of the porous specimen subjected to internal fluid pressure}
	\label{Geometry and boundary conditions of the porous specimen subjected to internal fluid pressure}
	\end{figure}

Geometry and boundary conditions of the 2D examples are shown in Fig. \ref{Geometry and boundary conditions of the porous specimen subjected to internal fluid pressure}. The pre-existing crack is placed horizontally at the center of a square specimen of 0.5 m $\times$ 0.5 m. The initial length of the crack is 0.05 m. We use the initial history field to induce the pre-existing crack and fluid source term $q_{F}=$ 10000 kg/(m$^3\cdot \textrm s$) is set in the pre-existing crack. The crack propagation is therefore driven by the fluid injection in the pre-existing crack. 

\begin{figure}[htbp]
	\centering
	\subfigure[$t=0$]{\includegraphics[height = 5cm]{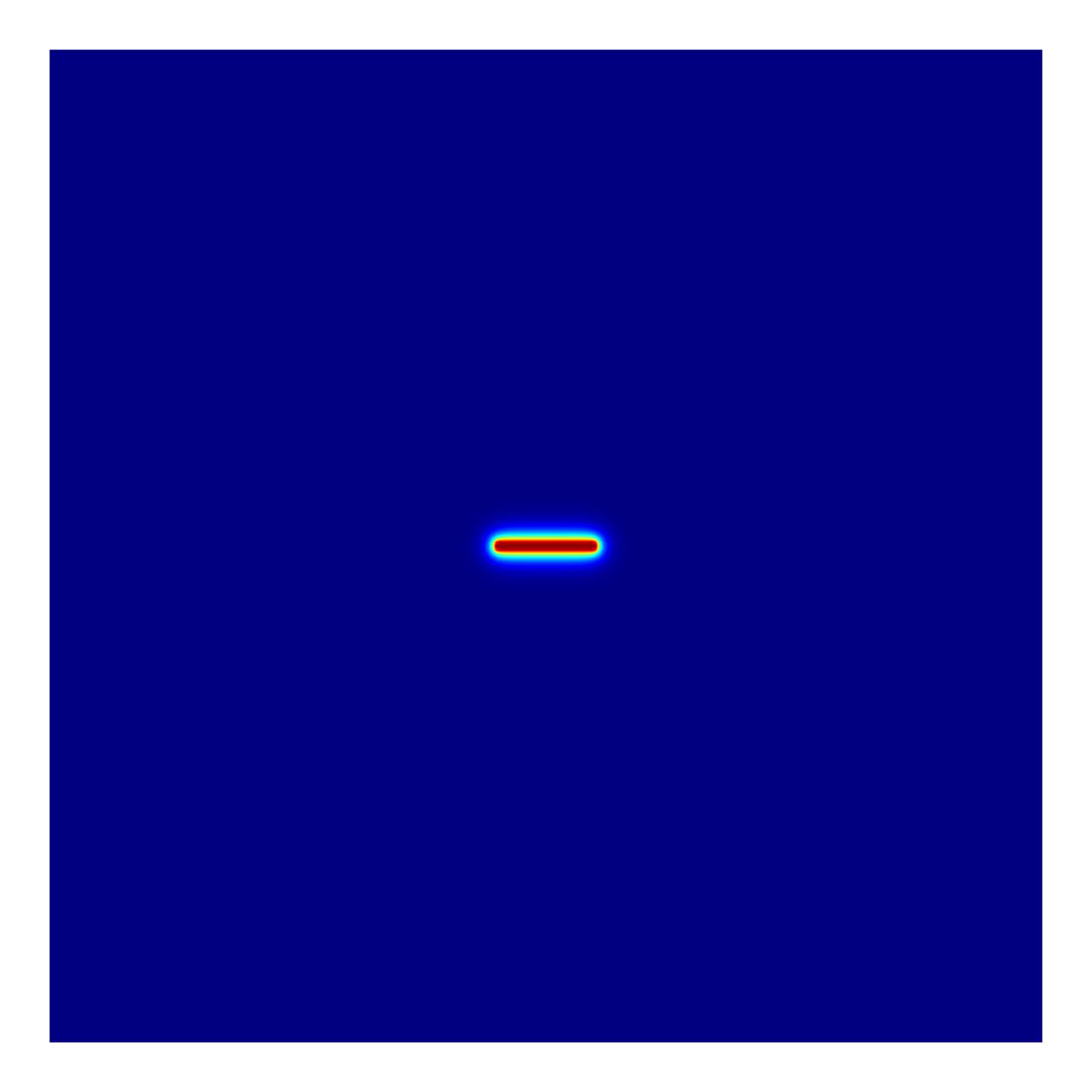}}
	\subfigure[$t=430$ $\mu$s]{\includegraphics[height = 5cm]{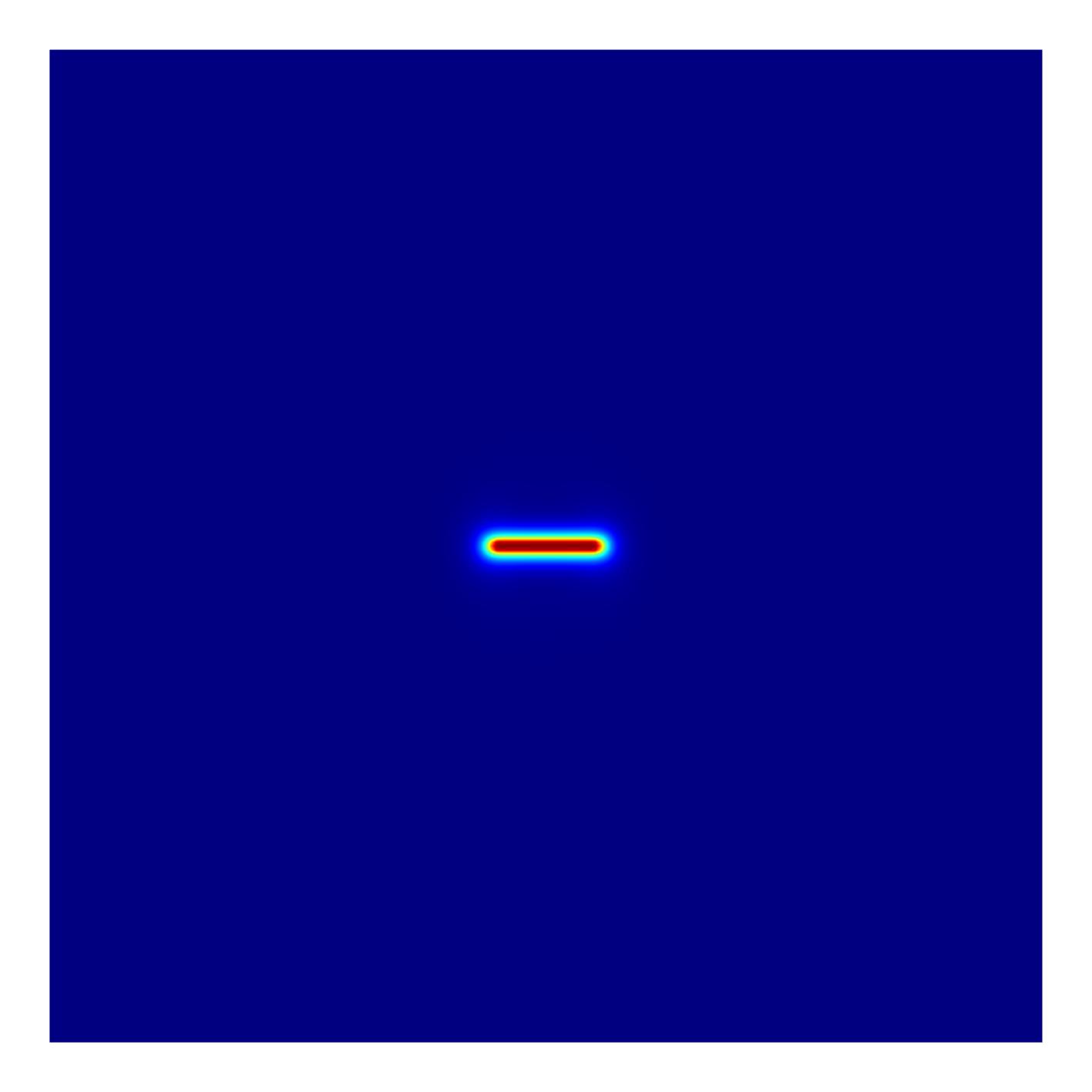}}
	\subfigure[$t=480$ $\mu$s]{\includegraphics[height = 5cm]{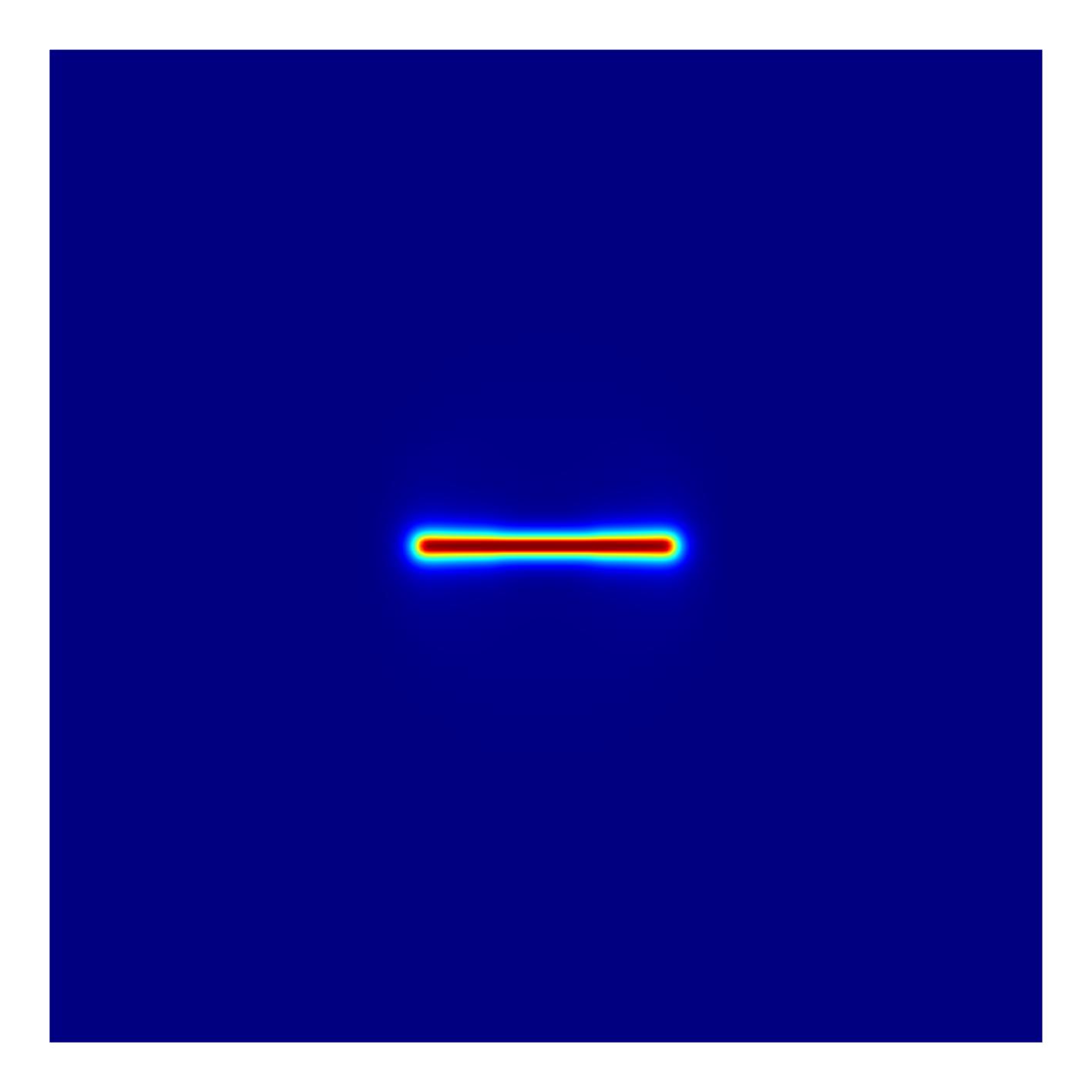}}\\
	\subfigure[$t=495$ $\mu$s]{\includegraphics[height = 5cm]{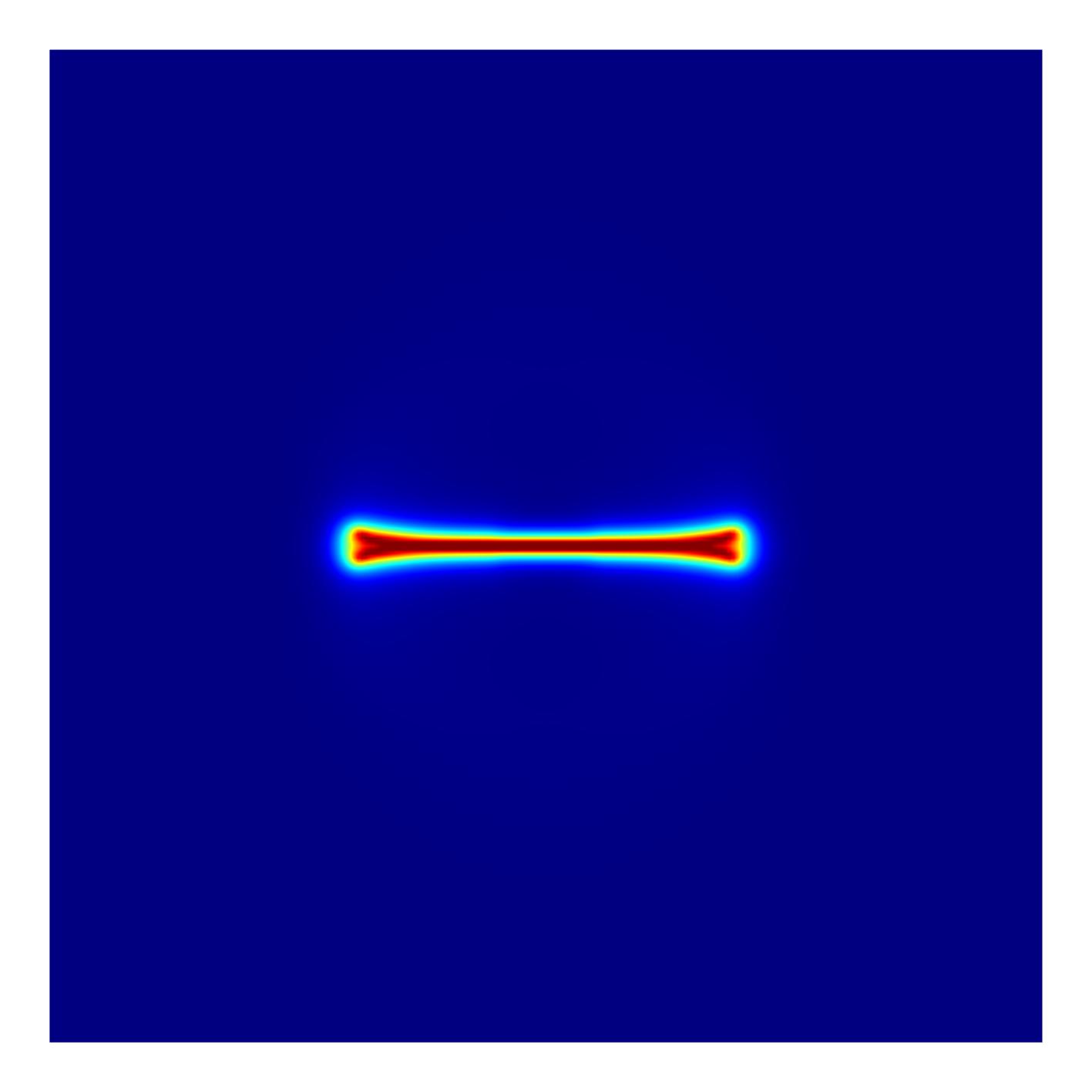}}
	\subfigure[$t=510$ $\mu$s]{\includegraphics[height = 5cm]{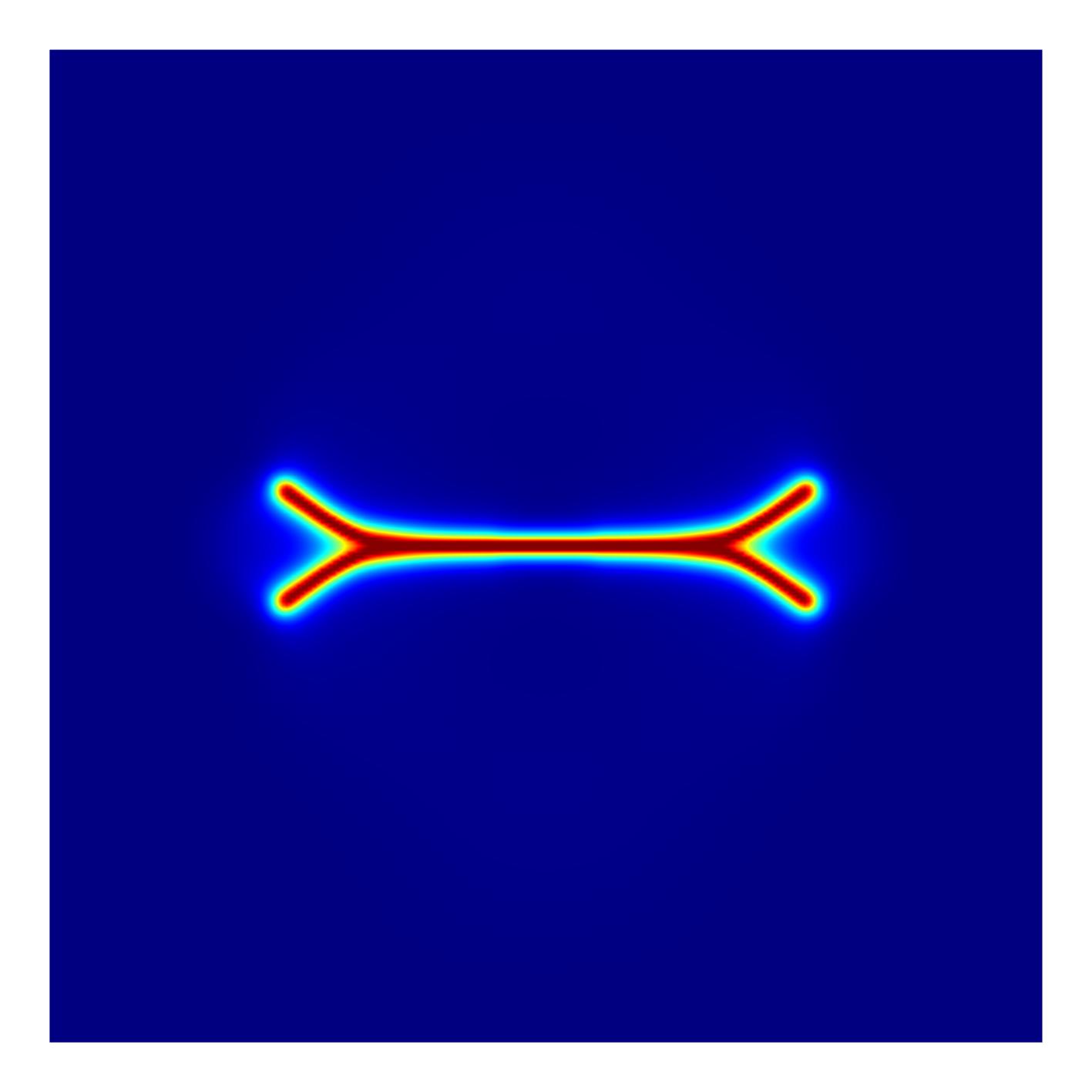}}
	\subfigure[$t=570$ $\mu$s]{\includegraphics[height = 5cm]{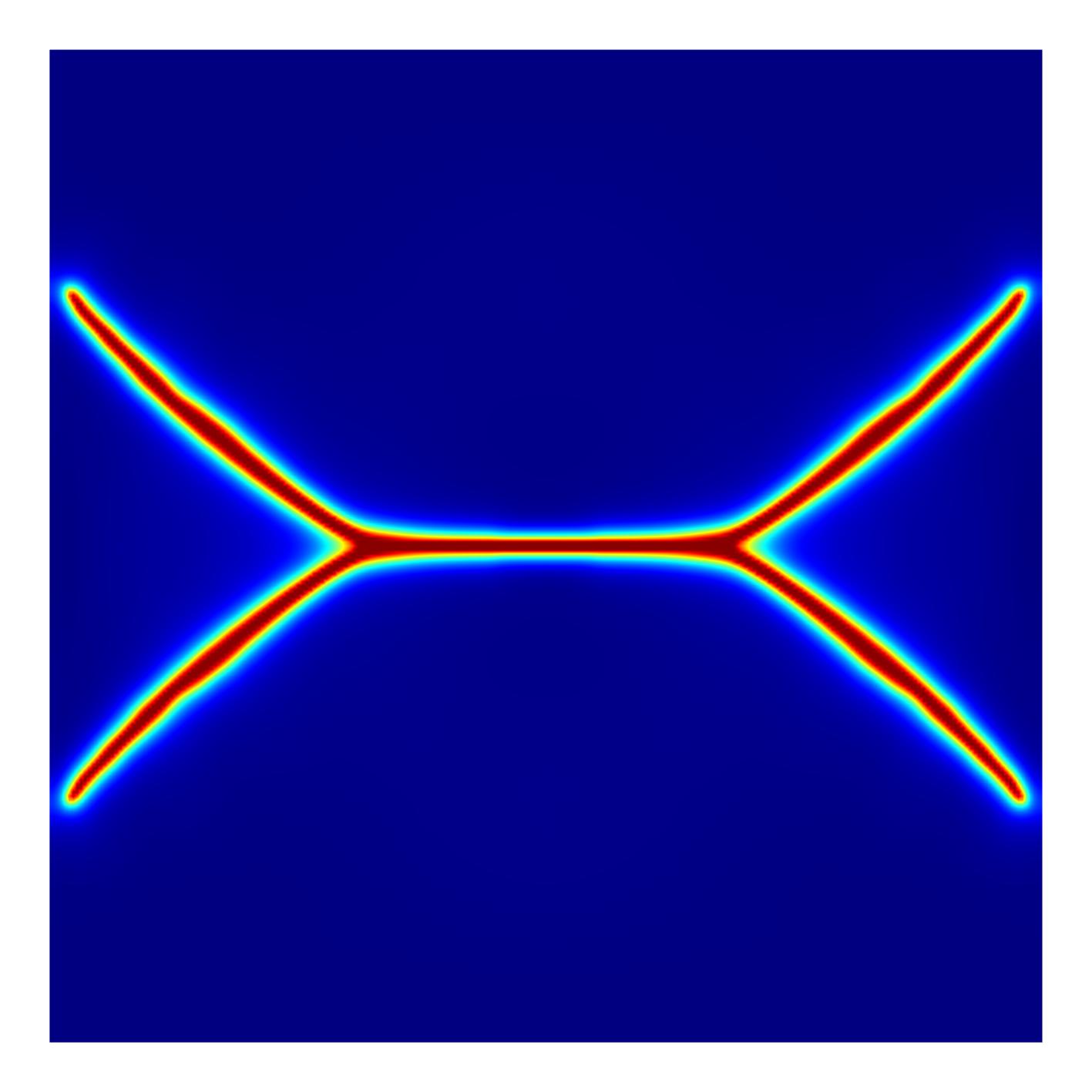}}
	\caption{Progressive crack propagation in the specimen with a pre-existing crack subjected to internal fluid injection}
	\label{Progressive crack propagation in the specimen with a pre-existing crack subjected to internal fluid injection}
	\end{figure}

The parameters used for calculation are listed in Table \ref{Base parameters for the specimen with a pre-existing crack subjected to internal fluid injection}. We first use uniform Q4 elements to discretize all the fields with the element size $h=$ $2\times 10^{-3}$ m. In addition, the time step $\Delta t$ is set as 0.1 $\mu$s. The progressive crack evolution with time $t$ is shown in Fig. \ref{Progressive crack propagation in the specimen with a pre-existing crack subjected to internal fluid injection}. It should be noted that for the dynamic fracture, the number of staggered iterations required in one time step is no more than 1 before fracture initiation and around 4 after fracture initiation and propagation, which is consistent with the observations in \citep{lee2017iterative}. When $t=0$, an induced crack is shown in Fig. \ref{Progressive crack propagation in the specimen with a pre-existing crack subjected to internal fluid injection}a while the crack starts to propagate when $t=430$ $\mu$s in Fig. \ref{Progressive crack propagation in the specimen with a pre-existing crack subjected to internal fluid injection}b. The initiated crack propagates along the horizontal direction when $t=480$ $\mu$s. When $t=495$ $\mu$s, the propagating crack starts to branch. Four similar cracks occur in Fig. \ref{Progressive crack propagation in the specimen with a pre-existing crack subjected to internal fluid injection}d when $t=510$ $\mu$s. The bifurcated cracks continue to propagate at a large angle with the horizontal direction and approach the left and right boundaries of the specimen when $t=595$ $\mu$s. The reason for crack branching is that the stress in a fracture zone is relatively high and the porous solid is unable to dissipate the energy that is driving the failure with a single crack, especially when the parameter $G_c$ is small. In addition, crack branches are observed to highly related to crack speed \citep{zhou2018phase}. Therefore, a faster loading rate, higher-energy release and higher stress in the failure part will cause a greater number of crack branching.

	\begin{table}[htbp]
	\small
	\caption{Base parameters for the specimen with a pre-existing crack subjected to internal fluid injection}
	\label{Base parameters for the specimen with a pre-existing crack subjected to internal fluid injection}
	\centering
	\begin{tabular}{llllll}
	\hline
	$E$ & 210 GPa & $\nu$ &0.3 & $G_c$ &$1\times10^{-2}$ N/m\\
	$k$ &$1\times10^{-9}$ & $l_0$ & $4\times10^{-3}$ m & $c_1$ & 0.4 \\
	 $c_2$ & 1.0 & $\varepsilon_{pR}$ & $2\times10^{-3}$ & $\rho_{R}$,  $\rho_{F}$& $1.0\times10^{3}$ kg/m$^3$  \\
	 $\alpha_R$ &$2\times10^{-3}$ & $q_R$ & 0 & $q_F$ & 10000 kg/(m$^3\cdot \textrm s$) \\ 
	$K_R$ & $1\times10^{-15}$ m$^2$ & $K_F$ & $1.333\times10^{-6}$ m$^2$ & $c_R$ & $1\times10^{-8}$ 1/Pa \\
	 $c_F$ & $1\times10^{-8}$ 1/Pa & $\mu_R$ & $1\times10^{-3}$ Pa$\cdot$s & $\mu_F$ & $1\times10^{-3}$ Pa$\cdot$s\\
	\hline
	\normalsize
	\end{tabular}
	\end{table}

\begin{figure}[htbp]
	\centering
	\subfigure[$\Delta t=0.1$ $\mu$s, $t=570$ $\mu$s]{\includegraphics[height = 5cm]{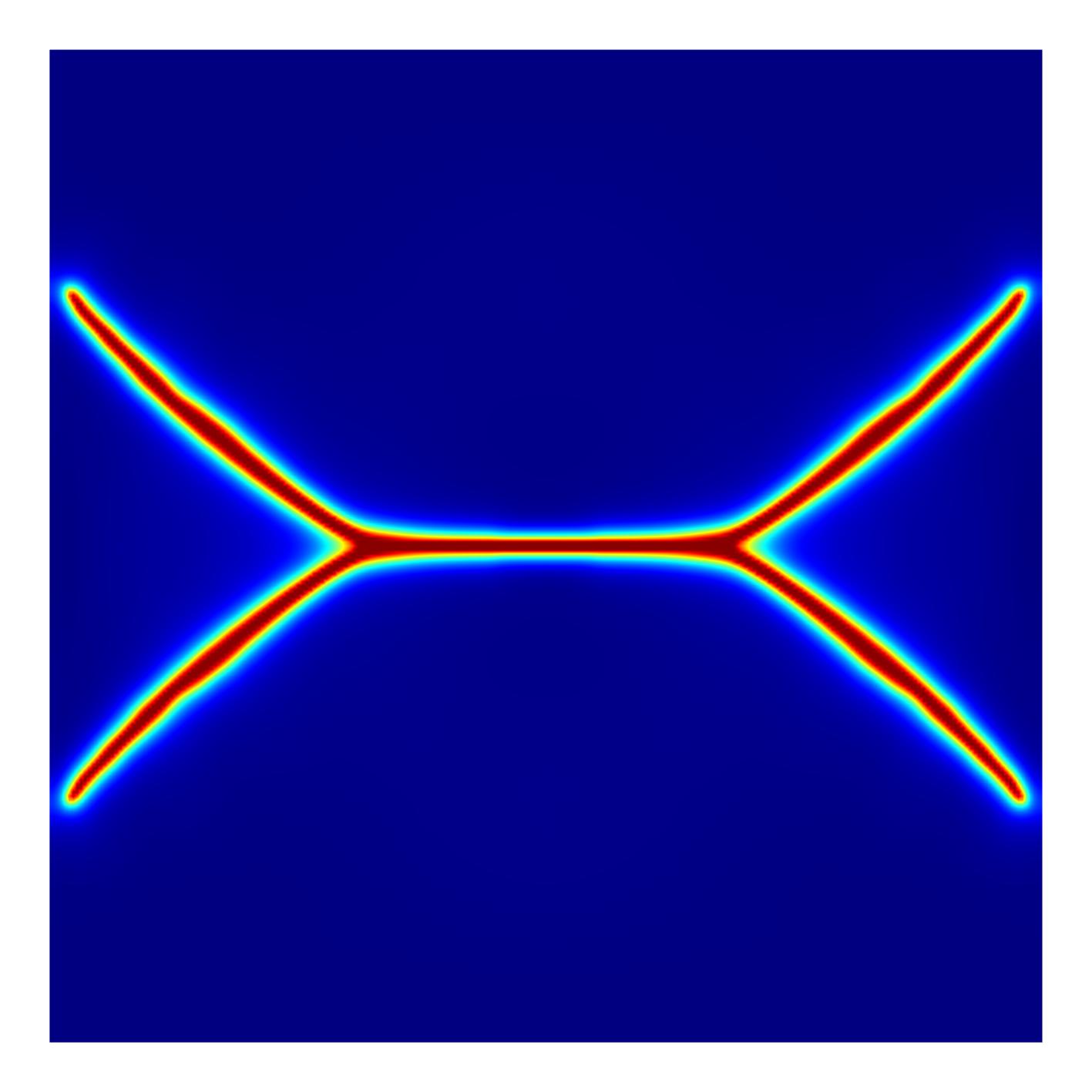}}
	\subfigure[$\Delta t=0.05$ $\mu$s, $t=585$ $\mu$s]{\includegraphics[height = 5cm]{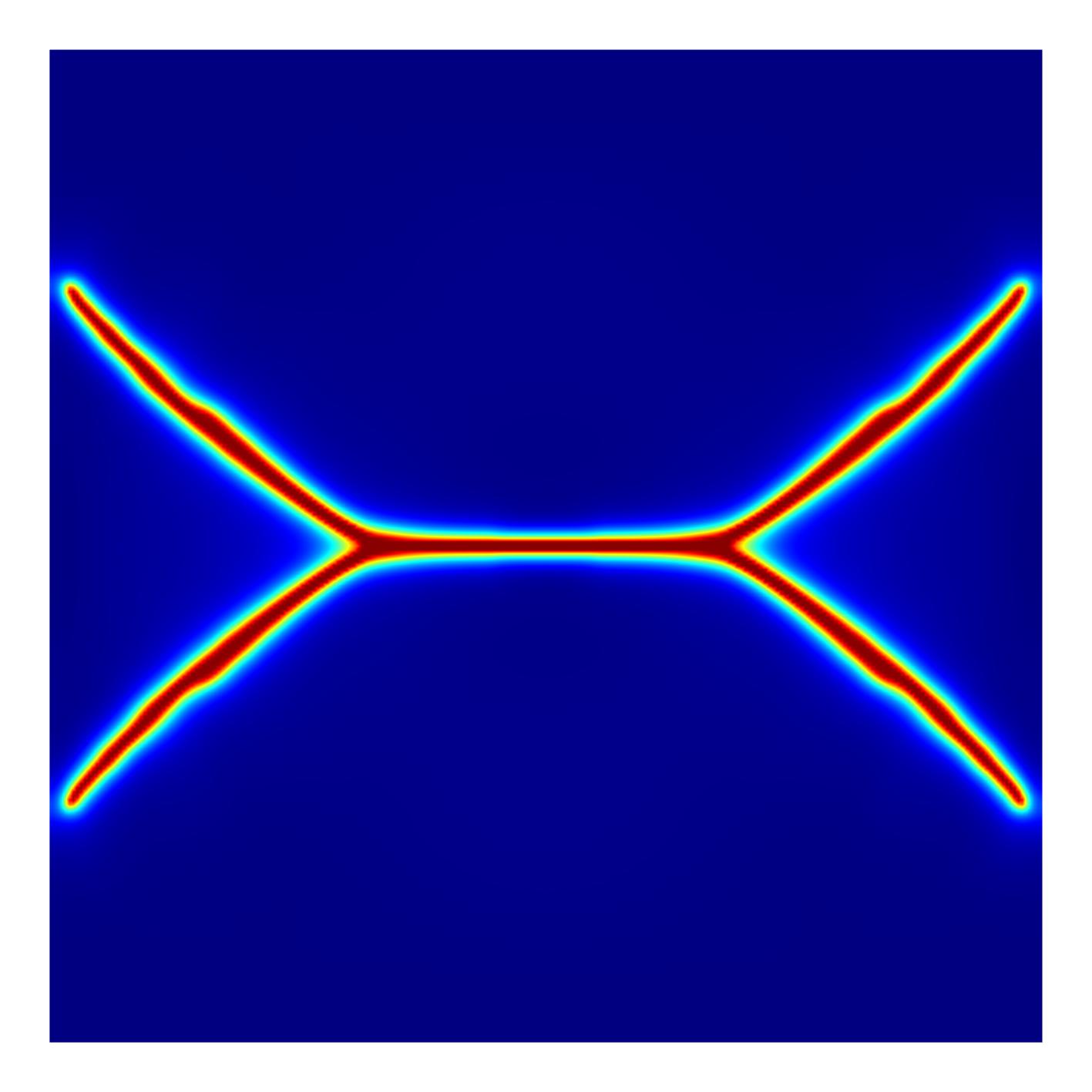}}
	\subfigure[$\Delta t=0.025$ $\mu$s, $t=560$ $\mu$s]{\includegraphics[height = 5cm]{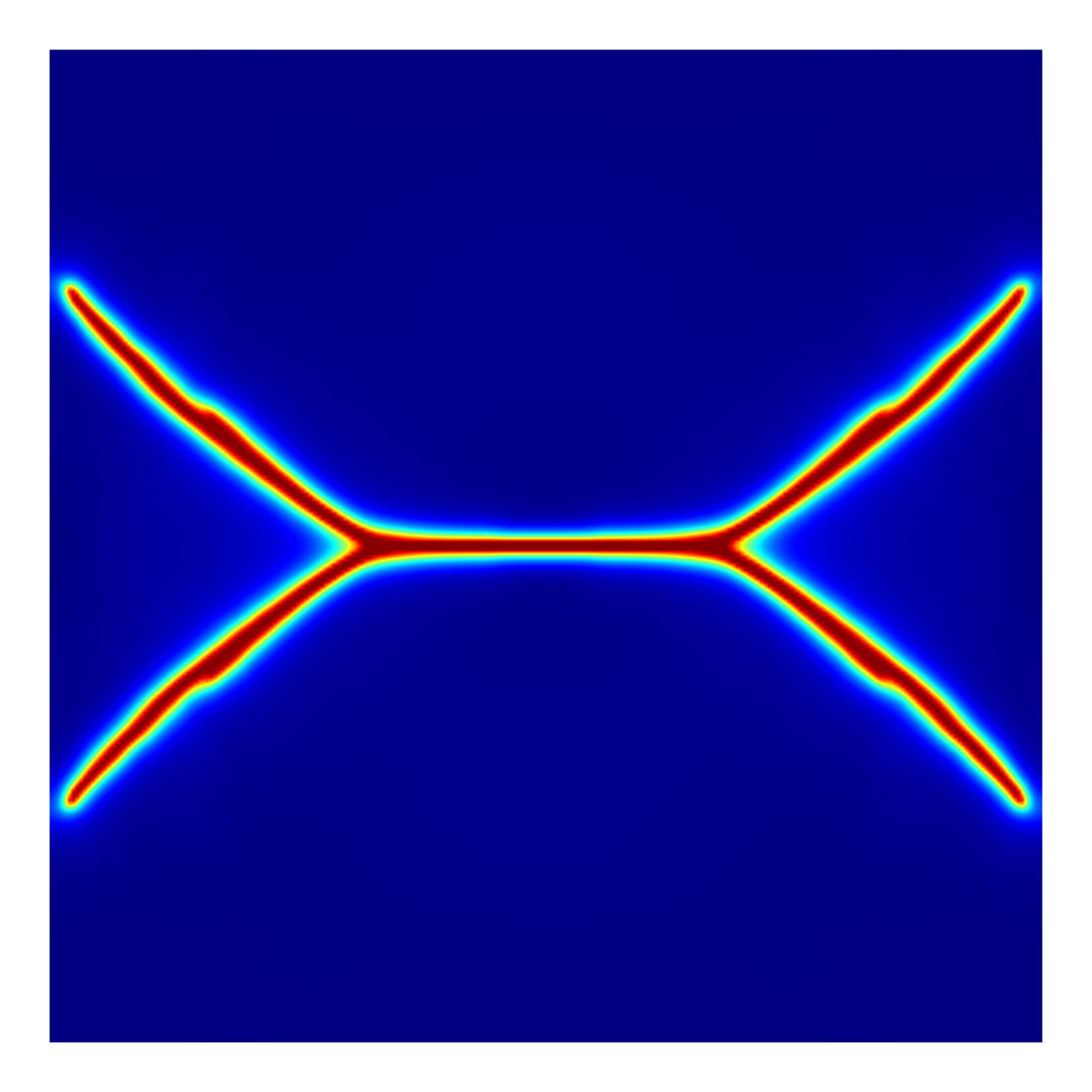}}
	\caption{Final crack patterns under different time steps}
	\label{Final crack patterns under different time steps}
	\end{figure}

\begin{figure}[htbp]
	\centering
	\includegraphics[width = 8cm]{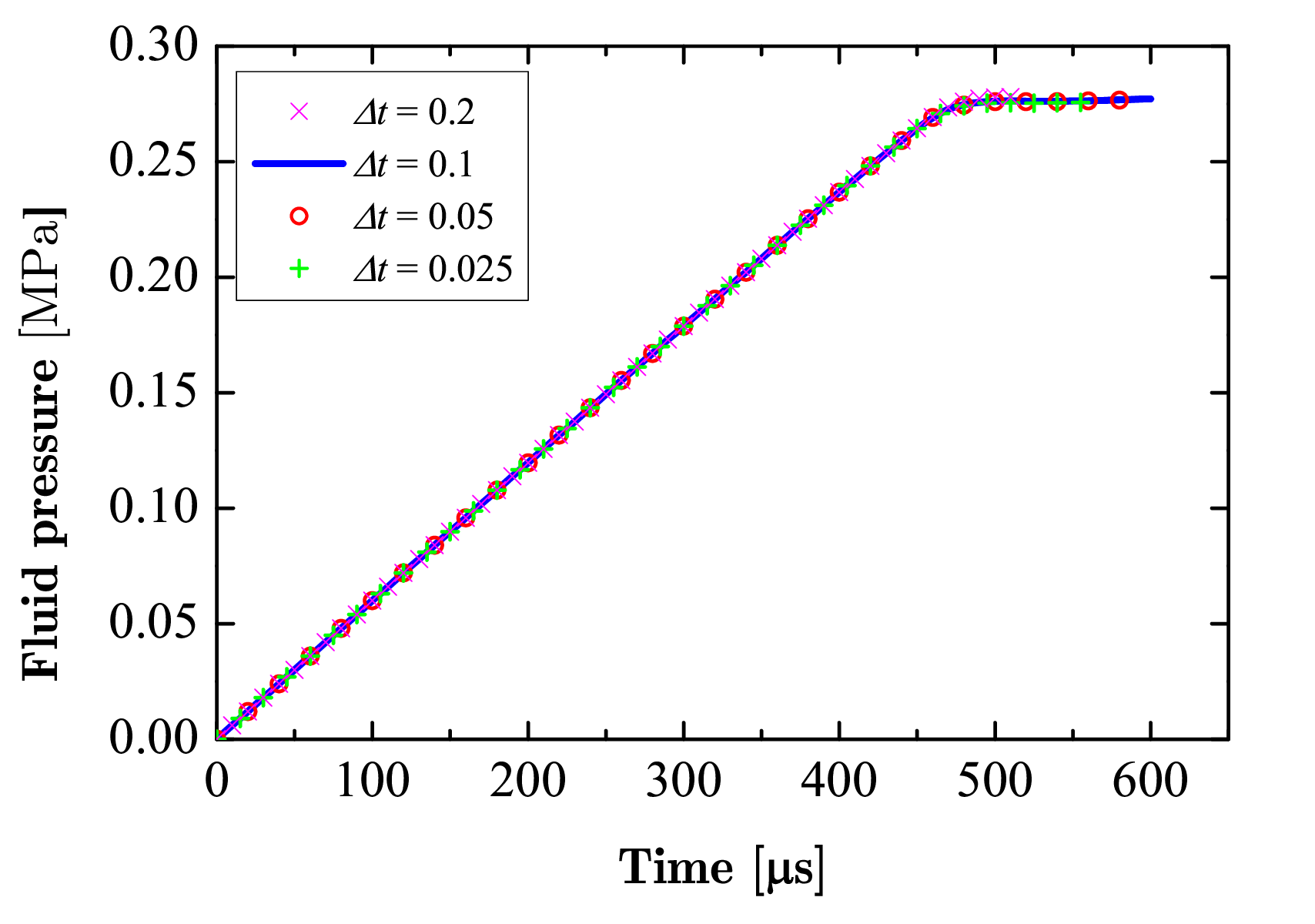}
	\caption{Fluid pressure under different time steps}
	\label{Fluid pressure under different time steps}
	\end{figure}

To test the influence of time step $\Delta t$, we change $\Delta t$ to 0.2 $\mu$s, 0.05 $\mu$s and 0.025$\mu$s and conduct the simulation again with other parameters unchanged. The final crack patterns under different time steps are shown in Fig. \ref{Final crack patterns under different time steps}. As observed, all the crack patterns are similar and the time step has little effect on the crack patterns. We also compare the fluid pressure under different time steps in Fig. \ref{Fluid pressure under different time steps}. The fluid pressure data come from the center of the initial pre-existing crack and time step is also observed to have little influence on the fluid pressure.

\begin{figure}[htbp]
	\centering
	\subfigure[$l_0=2\times10^{-3}$ m, $t=862$ $\mu$s]{\includegraphics[height = 5cm]{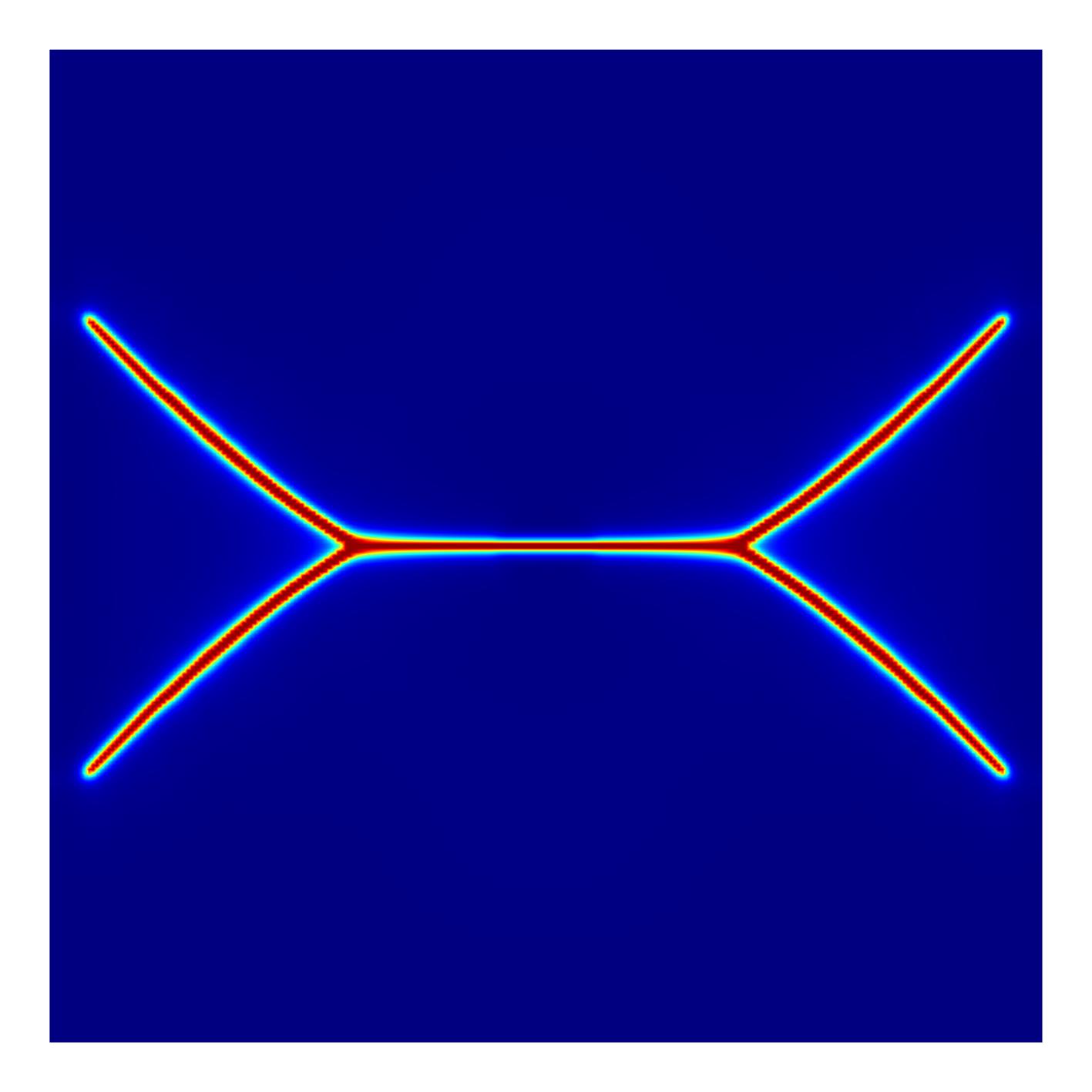}}
	\subfigure[$l_0=4\times10^{-3}$ m, $t=570$ $\mu$s]{\includegraphics[height = 5cm]{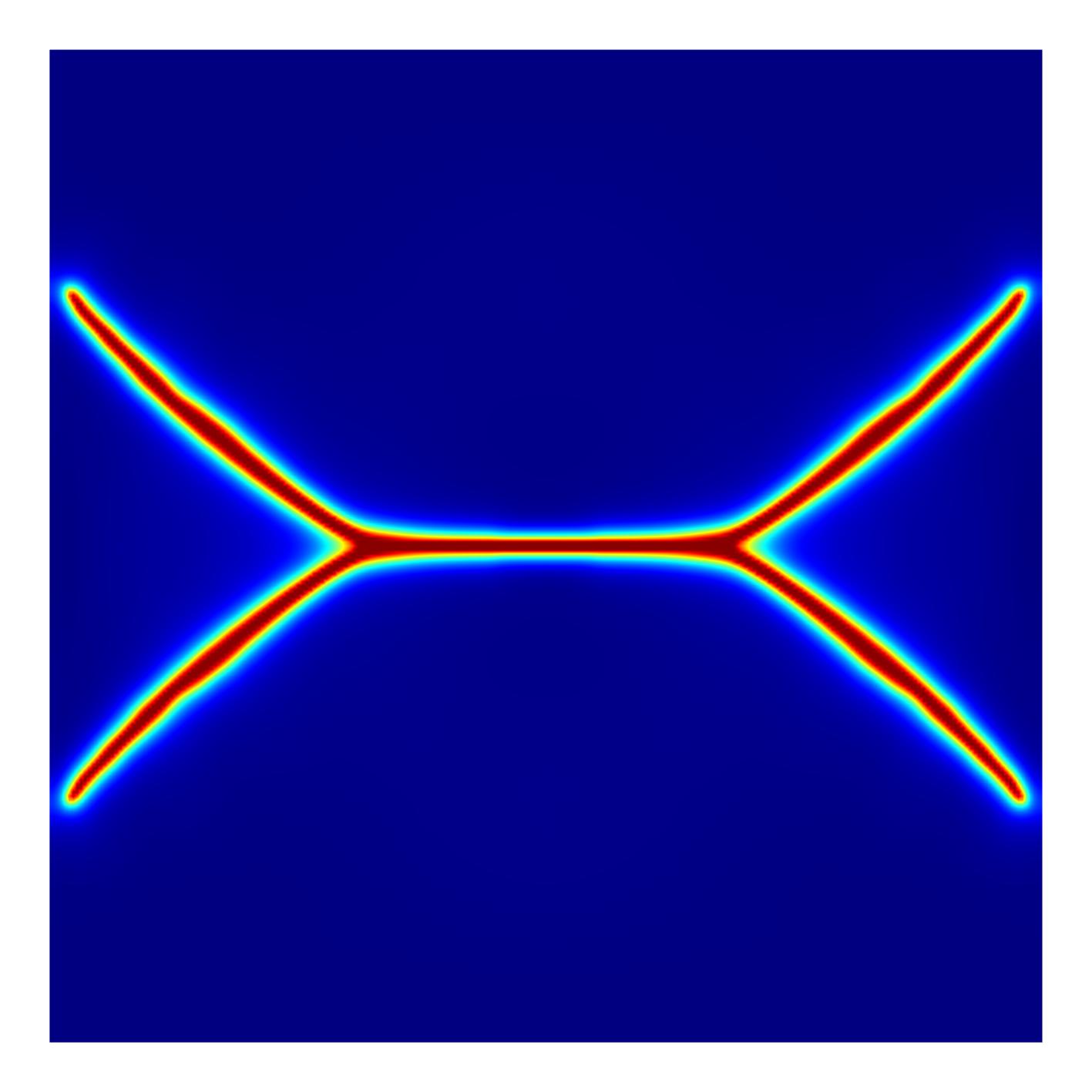}}
	\subfigure[$l_0=6\times10^{-3}$ m, $t=496$ $\mu$s]{\includegraphics[height = 5cm]{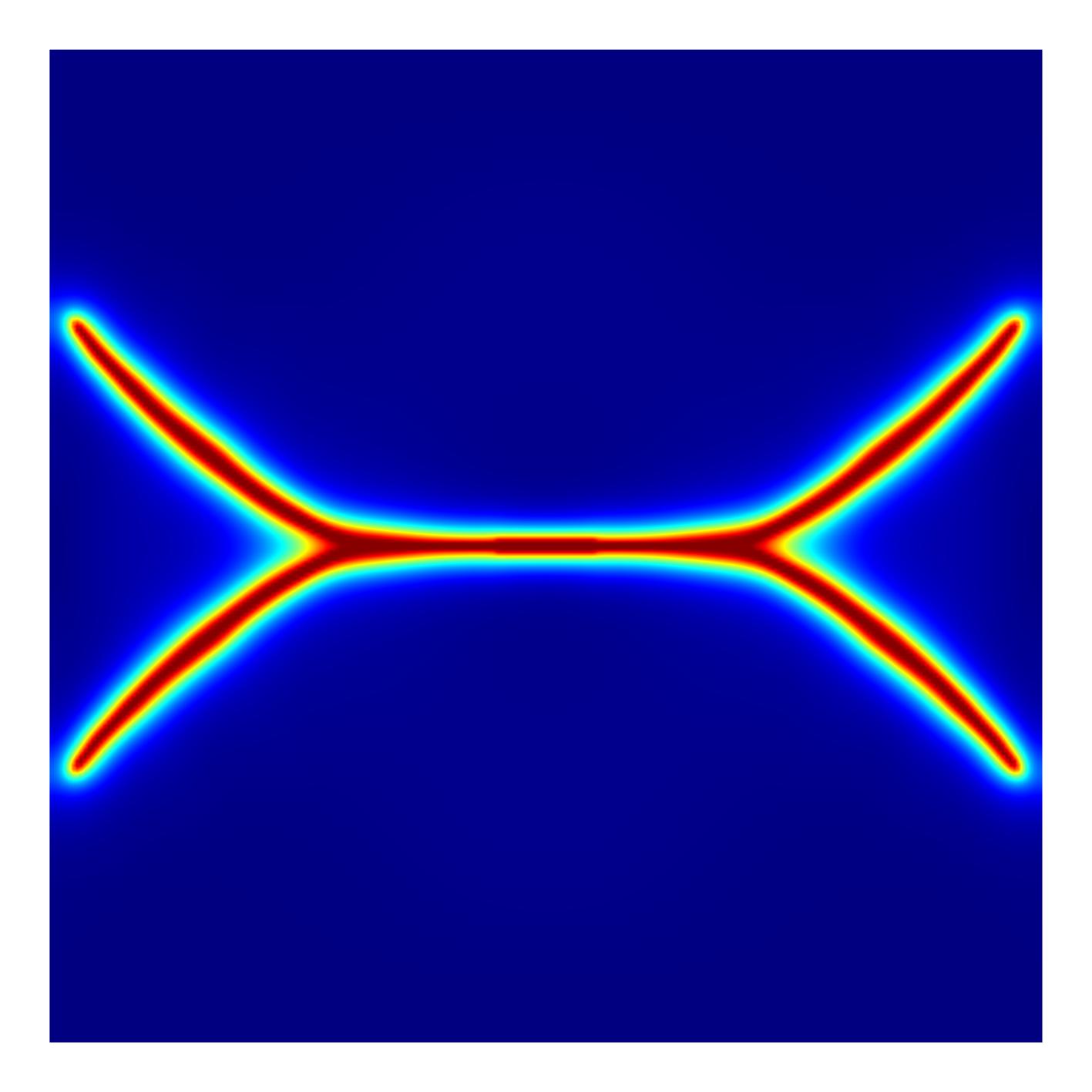}}
	\caption{Final crack patterns under different length scale parameters}
	\label{Final crack patterns under different length scale parameters}
	\end{figure}

\begin{figure}[htbp]
	\centering
	\includegraphics[width = 8cm]{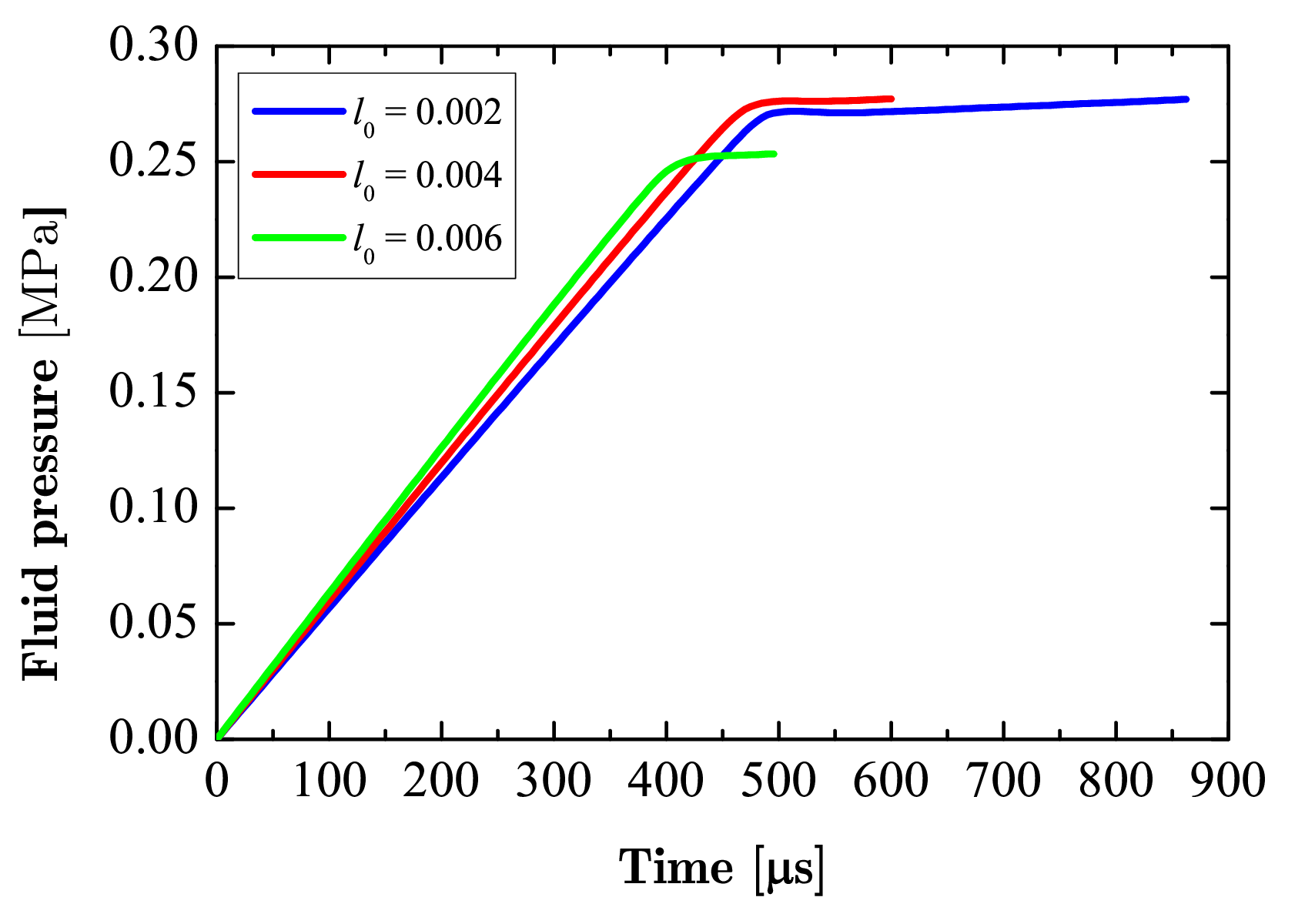}
	\caption{Fluid pressure under different length scale parameters}
	\label{Fluid pressure under different length scale parameters}
	\end{figure}

We then change the length scale $l_0$ to $2\times10^{-3}$, $4\times10^{-3}$, and $6\times10^{-3}$ m to test the influence of the length scale parameter. Likewise, the other parameters in Table \ref{Base parameters for the specimen with a pre-existing crack subjected to internal fluid injection} remain unchanged. Final crack patterns under different length scale parameters are shown in Fig. \ref{Final crack patterns under different length scale parameters}. As expected, the crack branching still occur and a larger length scale parameter $l_0$ results in a larger crack width. In addition, it will take more time for the crack to propagate to the left and right boundaries of the specimen for a smaller $l_0$. The comparison of fluid pressure under different length scale parameters is presented in Fig. \ref{Fluid pressure under different length scale parameters}. For a smaller length scale, the fluid pressure will have a longer post-peak stage.

\begin{figure}[htbp]
	\centering
	\subfigure[$G_c=1\times10^{-5}$ N/m, $t=91$ $\mu$s]{\includegraphics[height = 5cm]{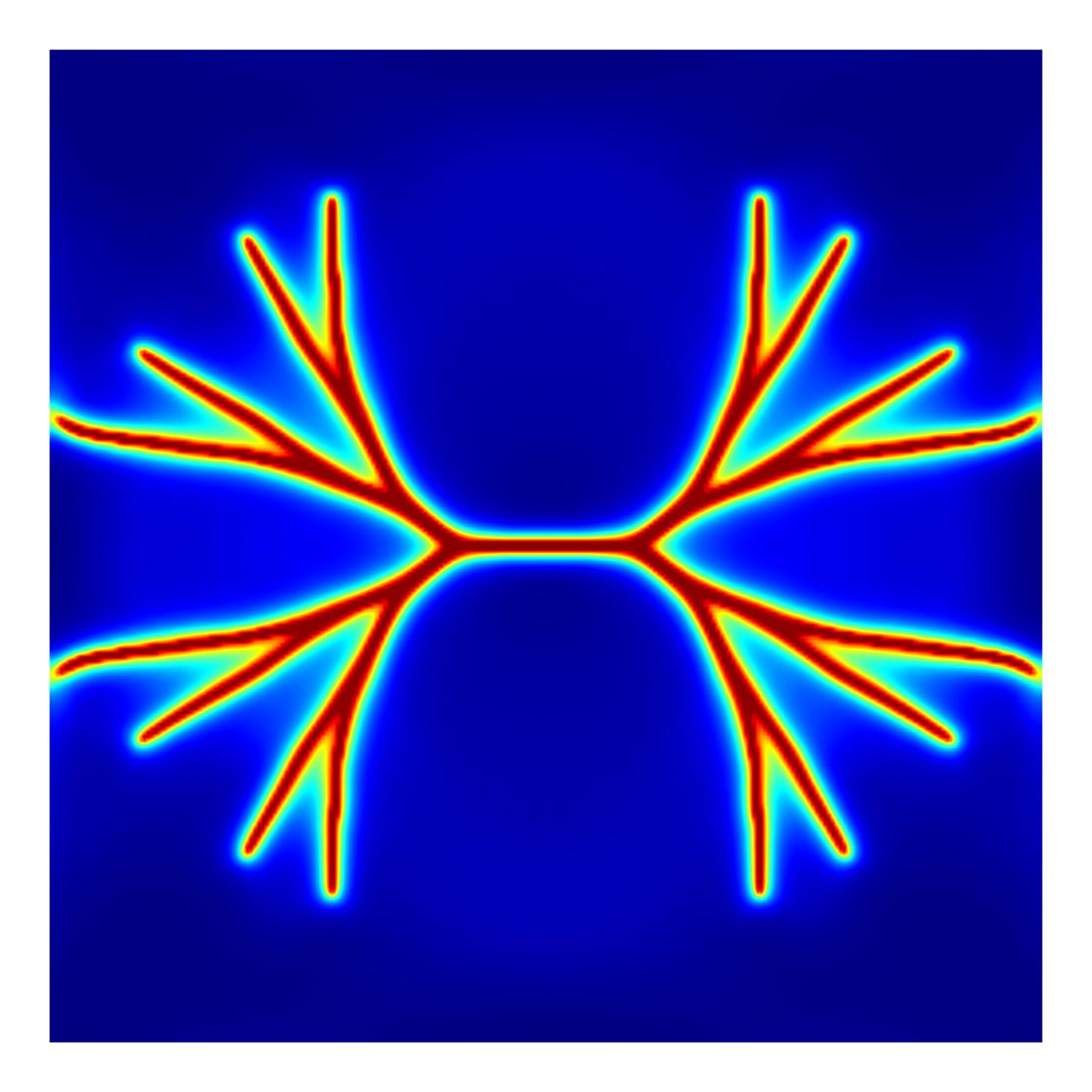}}
	\subfigure[$G_c=1\times10^{-4}$ N/m, $t=135$ $\mu$s]{\includegraphics[height = 5cm]{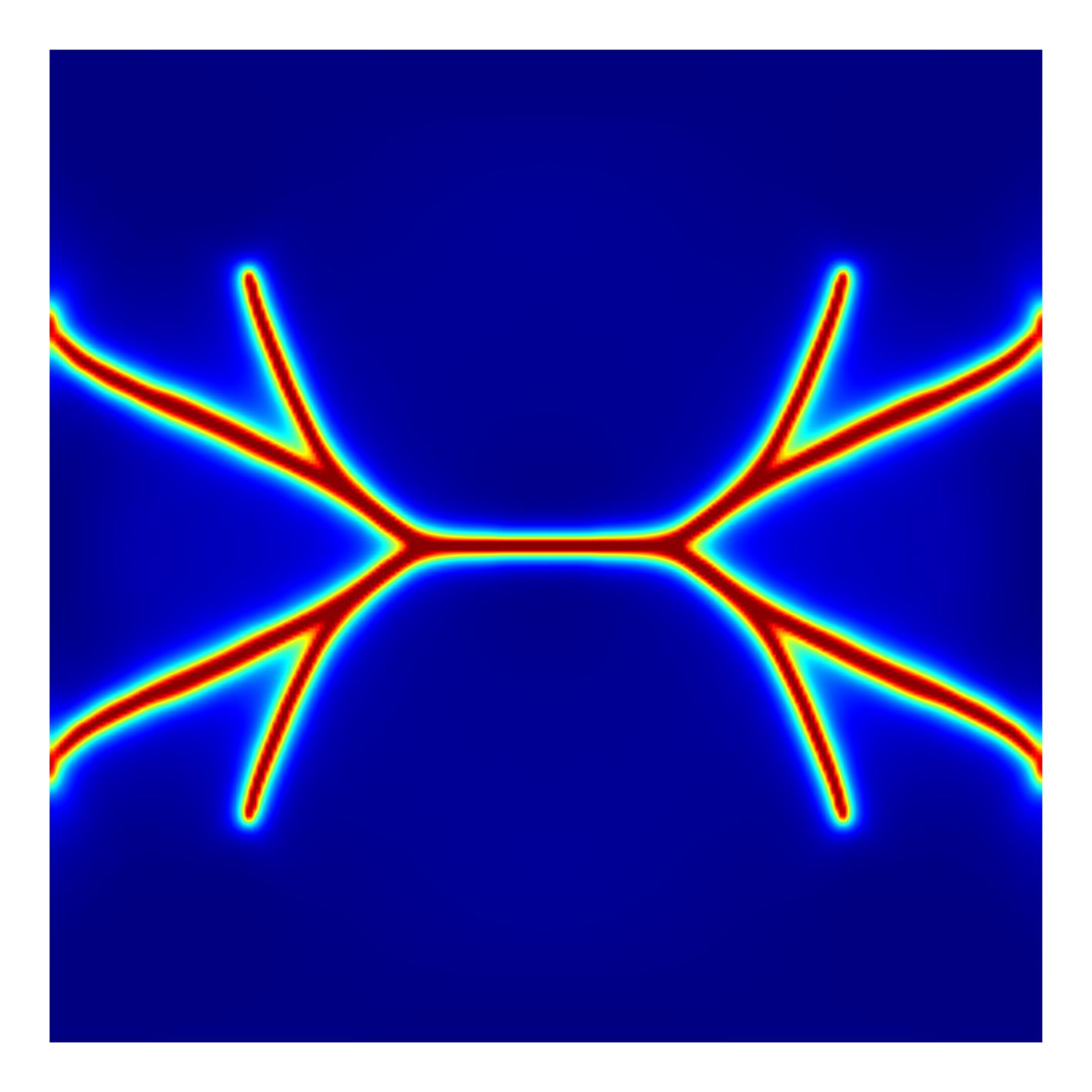}}\\
	\subfigure[$G_c=1\times10^{-3}$ N/m, $t=250$ $\mu$s]{\includegraphics[height = 5cm]{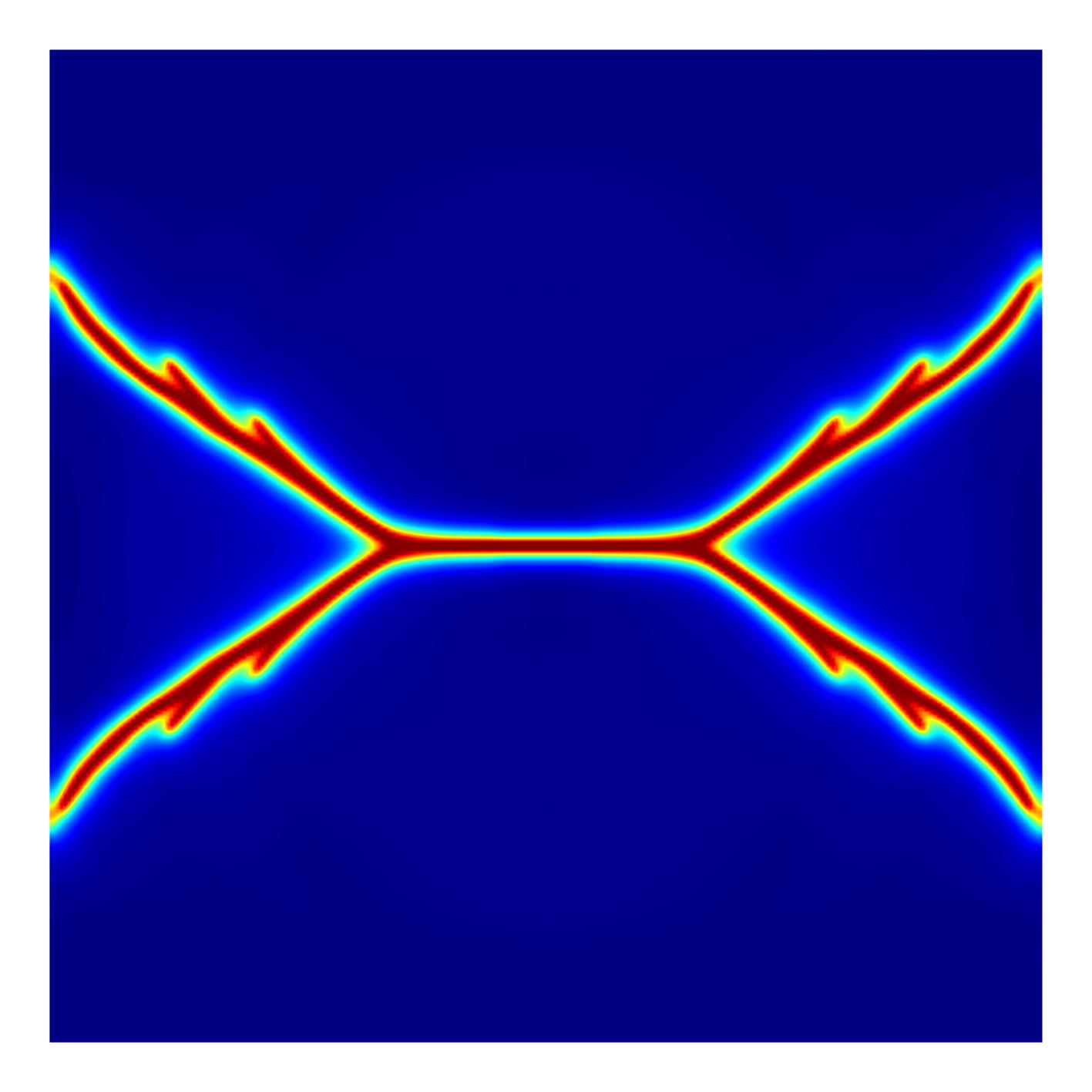}}
	\subfigure[$G_c=1\times10^{-2}$ N/m, $t=570$ $\mu$s]{\includegraphics[height = 5cm]{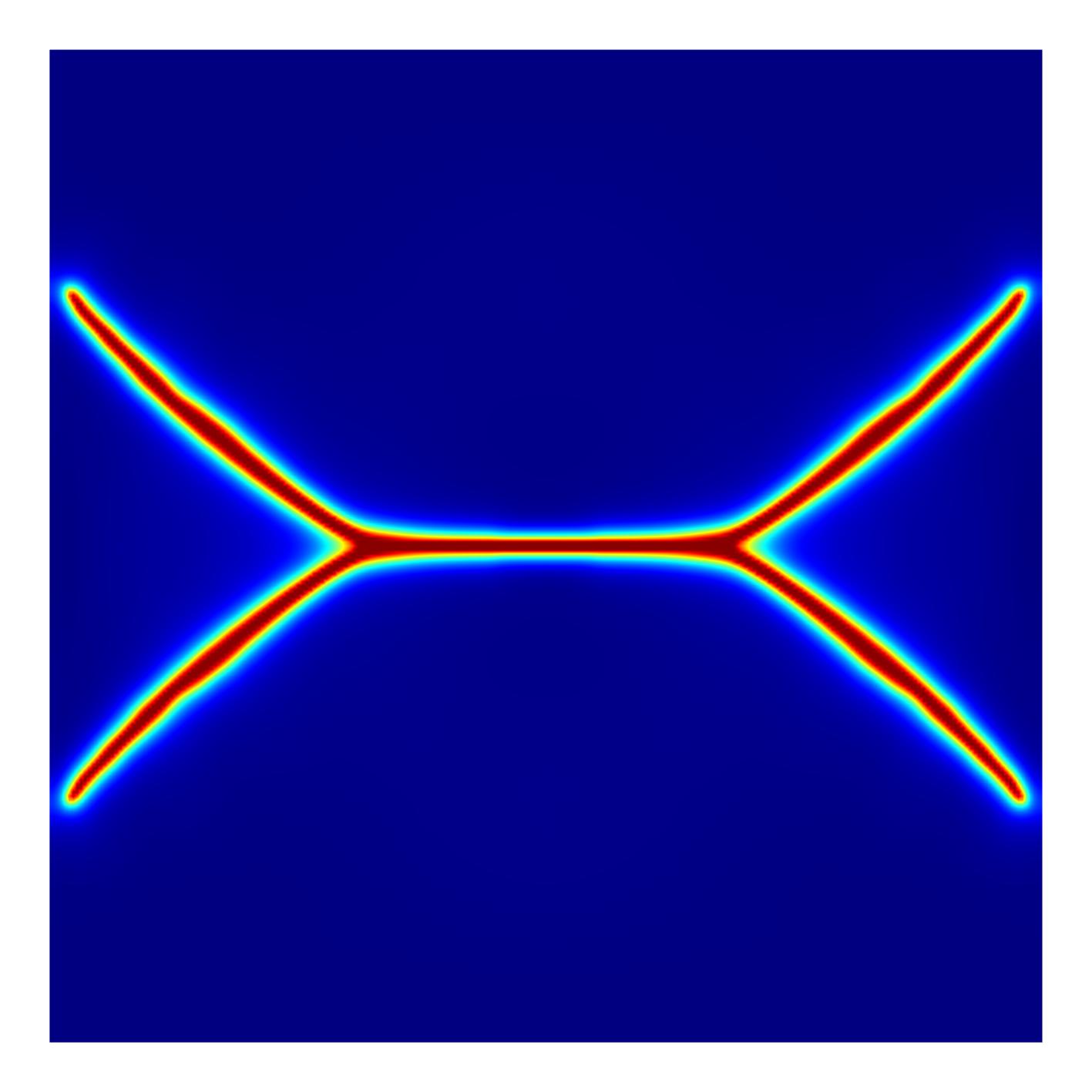}}
	\caption{Final crack patterns under different critical energy release rate $G_c$}
	\label{Final crack patterns under different critical energy release rate}
	\end{figure}

\begin{figure}[htbp]
	\centering
	\includegraphics[width = 8cm]{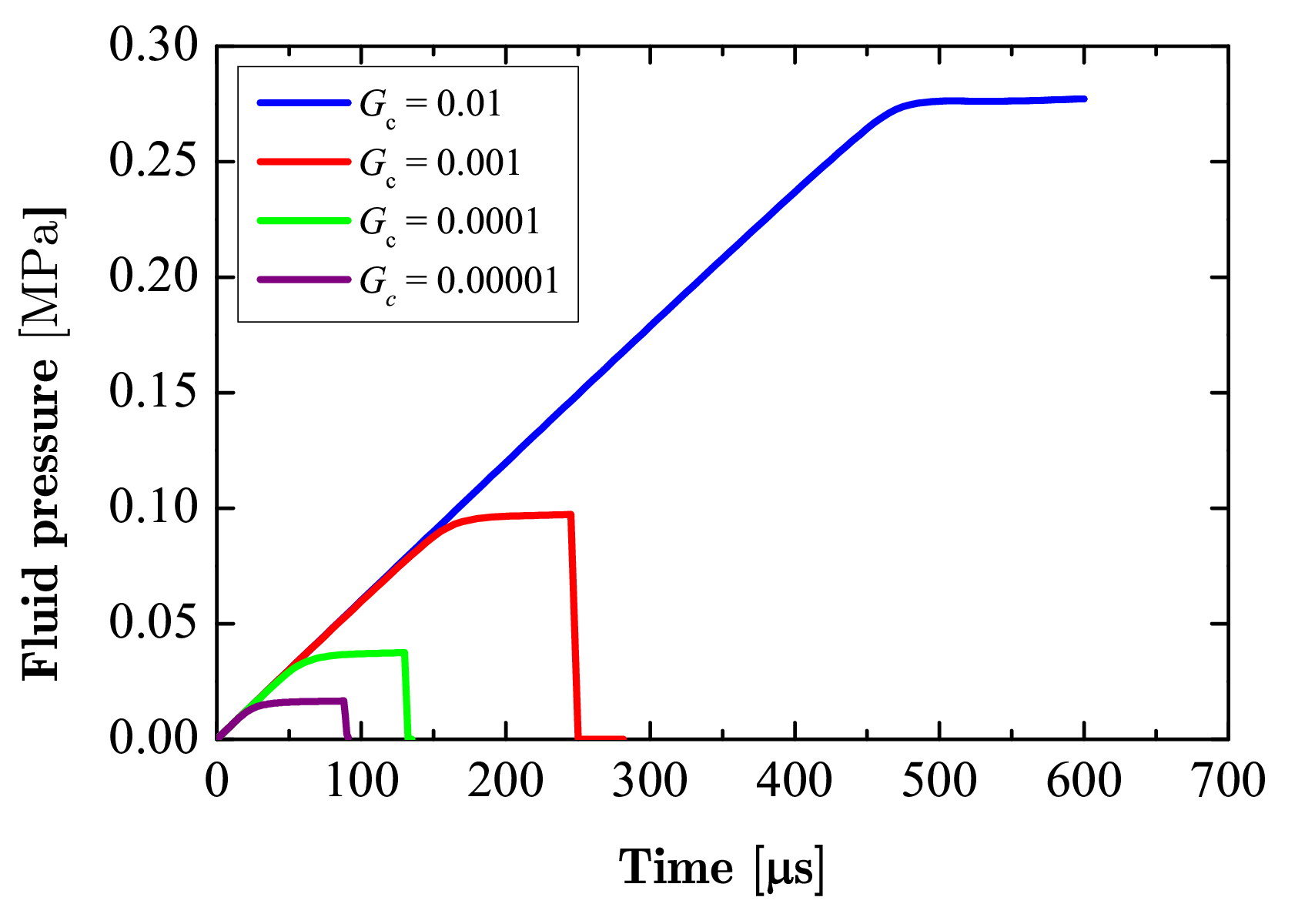}
	\caption{Fluid pressure under different critical energy release rate $G_c$}
	\label{Fluid pressure under different critical energy release rate}
	\end{figure}

In addition, we test the influence of critical energy release rate $G_c$. We fix the other parameters and change $G_c$ to $1\times10^{-5}$, $1\times10^{-4}$, $1\times10^{-3}$, and $1\times10^{-2}$ N/m, respectively. Similar simulations are conducted and the final crack patterns under different critical energy release rate are shown in Fig. \ref{Final crack patterns under different critical energy release rate}. More crack branching is observed and the cracks reach the left and right boundaries at a larger rate for a relatively smaller critical energy release rate $G_c$. Figure \ref{Fluid pressure under different critical energy release rate} shows the fluid pressure under different critical energy release rate. As observed, a larger $G_c$ will cause a larger fluid pressure.

\begin{figure}[htbp]
	\centering
	\includegraphics[width = 8cm]{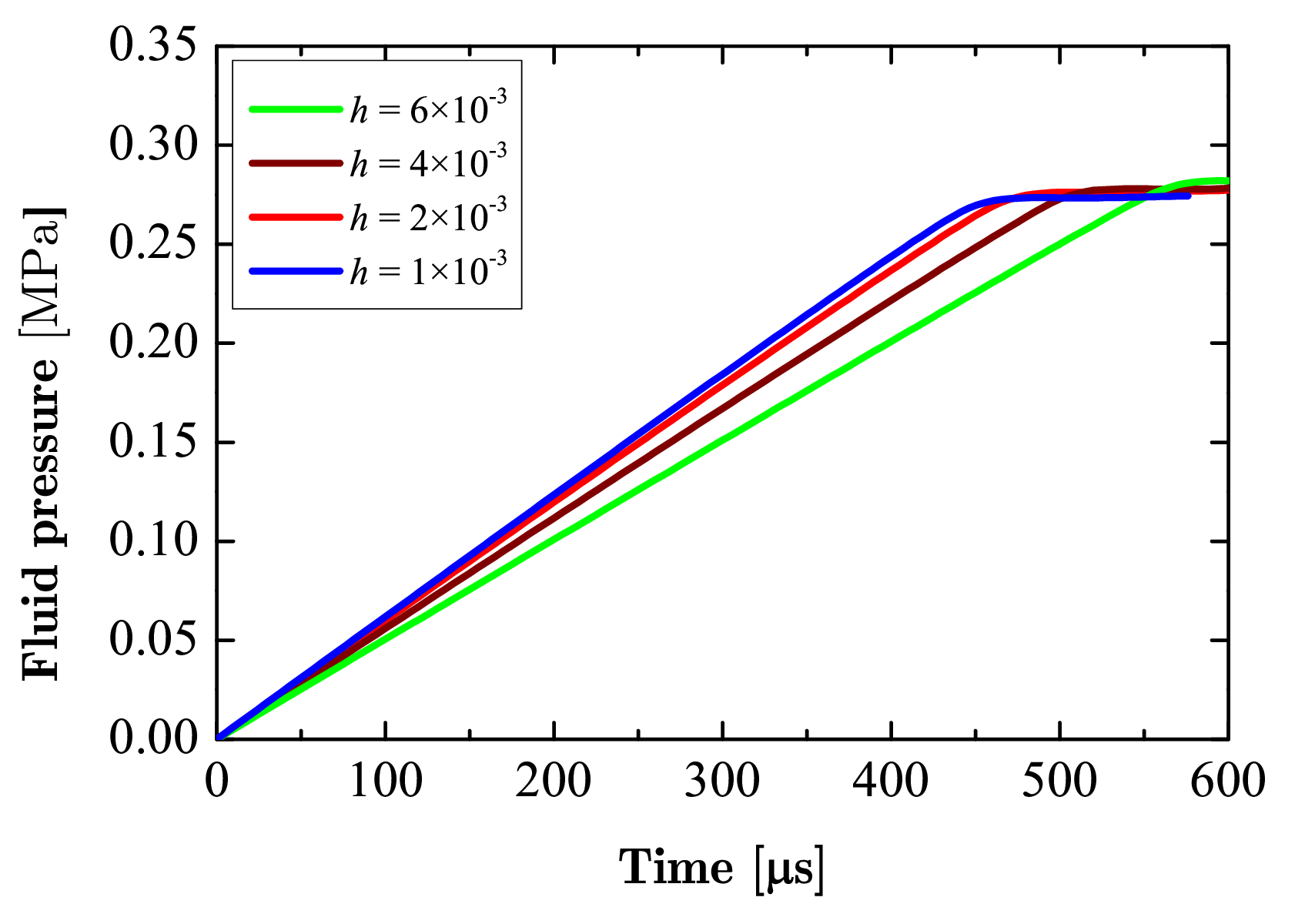}
	\caption{Fluid pressure under different mesh size $h$}
	\label{Fluid pressure under different mesh size}
	\end{figure}

\begin{figure}[htbp]
	\centering
	\subfigure[$q_{F}=$ 500 kg/(m$^3\cdot \textrm s)$, $t=9605$ $\mu$s]{\includegraphics[height = 5cm]{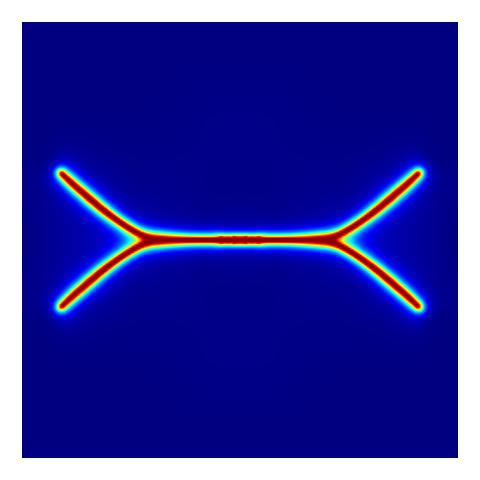}}
	\subfigure[$q_{F}=$ 1000 kg/(m$^3\cdot \textrm s)$, $t=5100$ $\mu$s]{\includegraphics[height = 5cm]{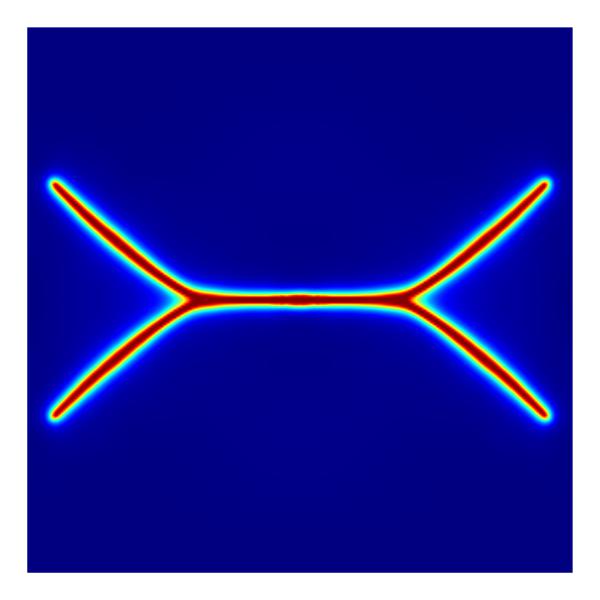}}\\
	\subfigure[$q_{F}=$ 5000 kg/(m$^3\cdot \textrm s)$, $t=1000$ $\mu$s]{\includegraphics[height = 5cm]{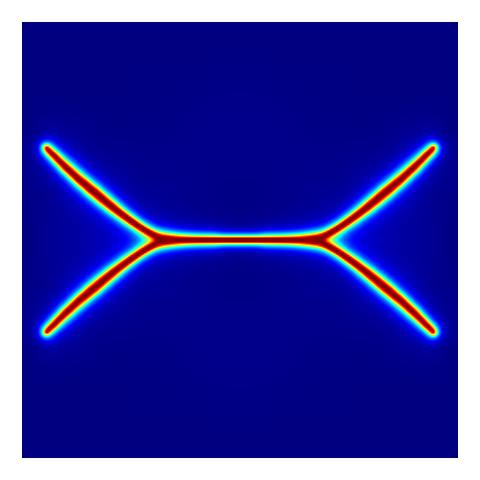}}
	\subfigure[$q_{F}=$ 10000 kg/(m$^3\cdot \textrm s)$, $t=570$ $\mu$s]{\includegraphics[height = 5cm]{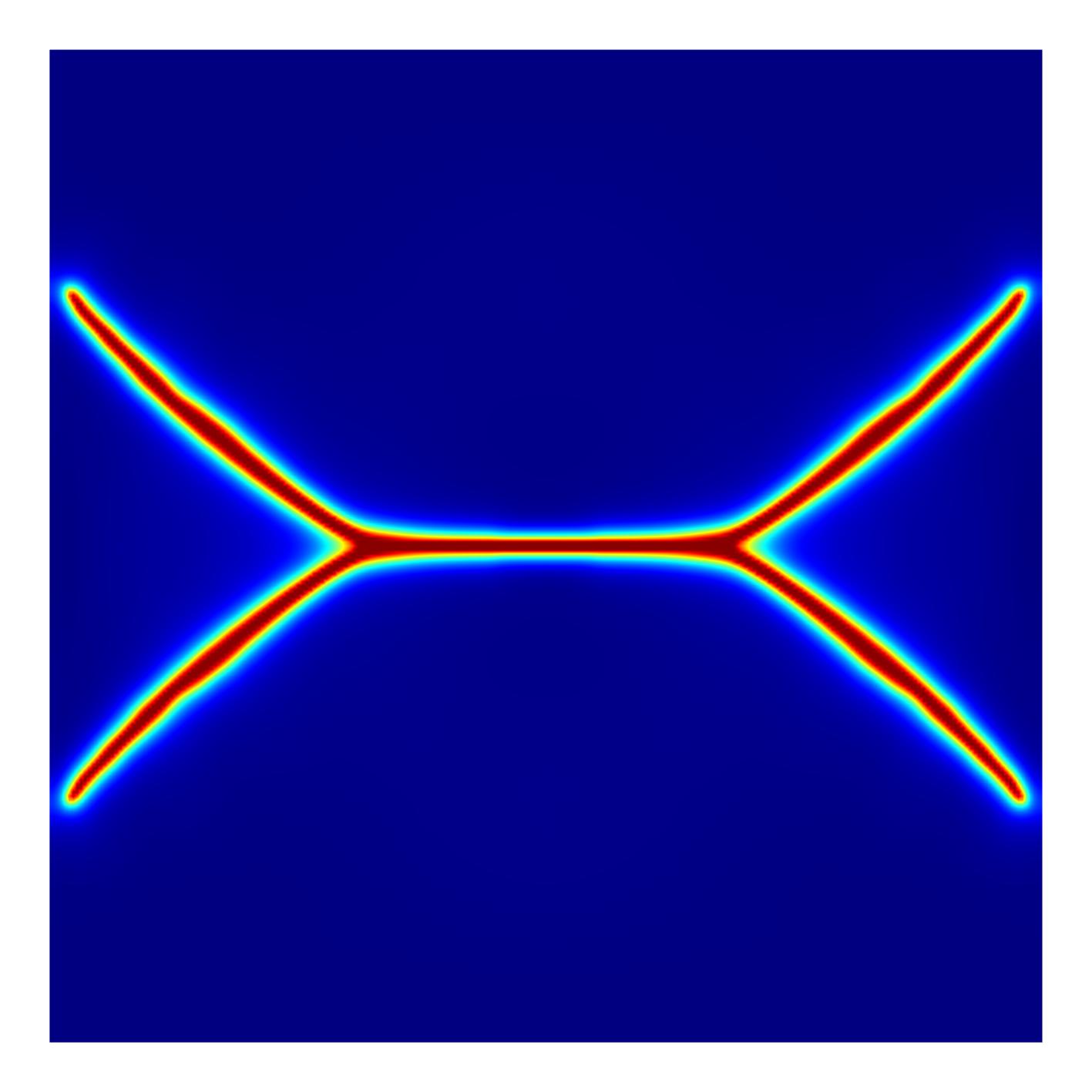}}
	\caption{Final crack patterns under different fluid source $q_F$}
	\label{Final crack patterns under different fluid source}
	\end{figure}

\begin{figure}[htbp]
	\centering
	\includegraphics[width = 8cm]{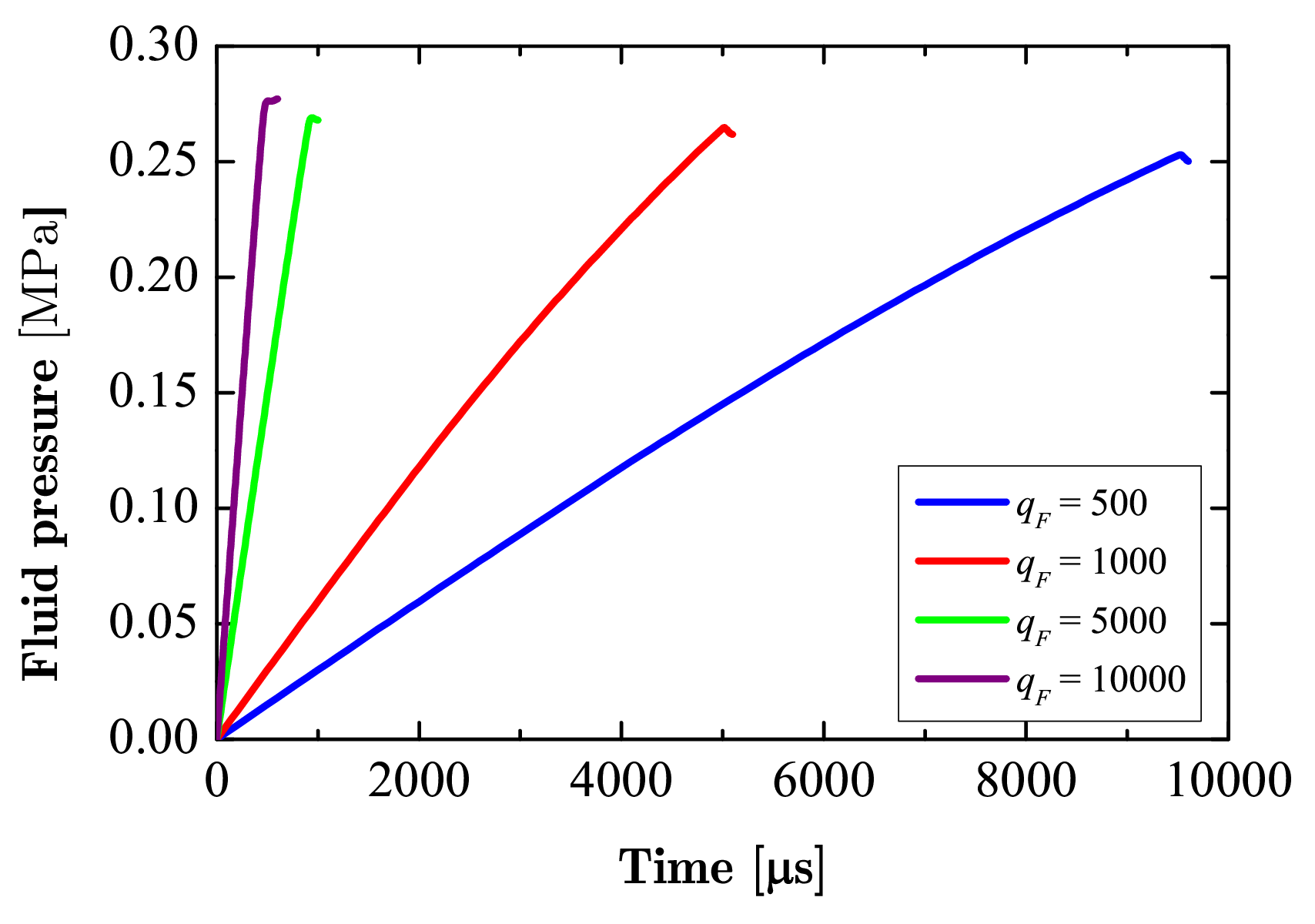}
	\caption{Fluid pressure under different fluid source $q_F$}
	\label{Fluid pressure under different fluid source}
	\end{figure}

We also change the mesh size $h$ to $1\times10^{-3}$, $4\times10^{-3}$, and $6\times10^{-3}$ m to investigate the influence of mesh refinement. It is found that the mesh size $h$ has a negligible effect on the fracture pattern unless it is too large and cannot exactly characterize the fracture path. Figure \ref{Fluid pressure under different mesh size} shows the fluid pressure under different mesh size $h$. As observed, if the mesh size decreases, the slope of the pressure curve increases. Finally, we test the influence of the fluid source term $q_F$. These $q_F$ are used in the simulation: $q_{F}= 500$, 1000, 5000, and 10000 kg/(m$^3\cdot \textrm s)$. Figure \ref{Final crack patterns under different fluid source} shows the final crack patterns in the porous domain under different $q_F$, which is observed to have little effect on the final crack pattern. Figure \ref{Fluid pressure under different fluid source} presents the fluid pressure in the middle of the pre-existing crack under different $q_F$. $q_F$ also has little effect on the maximum fluid pressure of the domain; however, a larger $q_F$ takes the shortest time to achieve the maximum fluid pressure.

\subsection{3D examples}\label{3D examples}

We now test a 3D example. A pre-existing penny-shaped cracks with a 0.04 m radius is placed at the center of a cubic specimen of 0.4 m $\times$ 0.4 m $\times$ 0.4 m. The same parameters listed in Table \ref{Base parameters for the specimen with a pre-existing crack subjected to internal fluid injection} are used except $l_0=6\times10^{-3}$ m. All the outer boundaries of the specimen are displacement-fixed and permeable with $p=0$. 6-node prism elements are used to discretize the volume and the element size is no more than $6\times10^{-3}$ m. Note that the 3D example is computed within 35 d 20 h by using 2 I5-6200U CPUs and 8GB physical RAM, although only around 4 staggered iterations are observed in one time step.

\begin{figure}[htbp]
	\centering
	\subfigure[$t=600$ $\mu$s]{\includegraphics[height = 6cm]{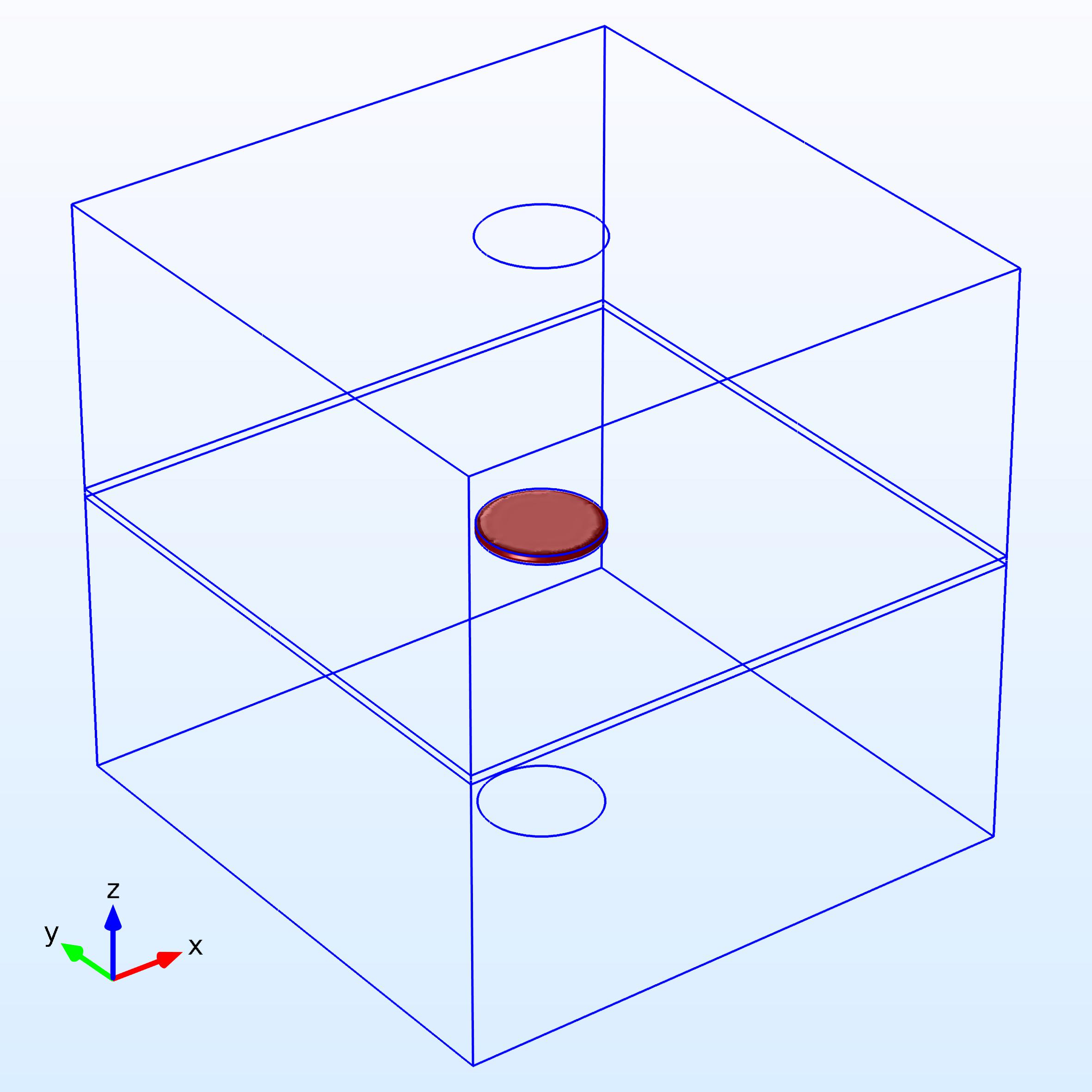}}
	\subfigure[$t=740$ $\mu$s]{\includegraphics[height = 6cm]{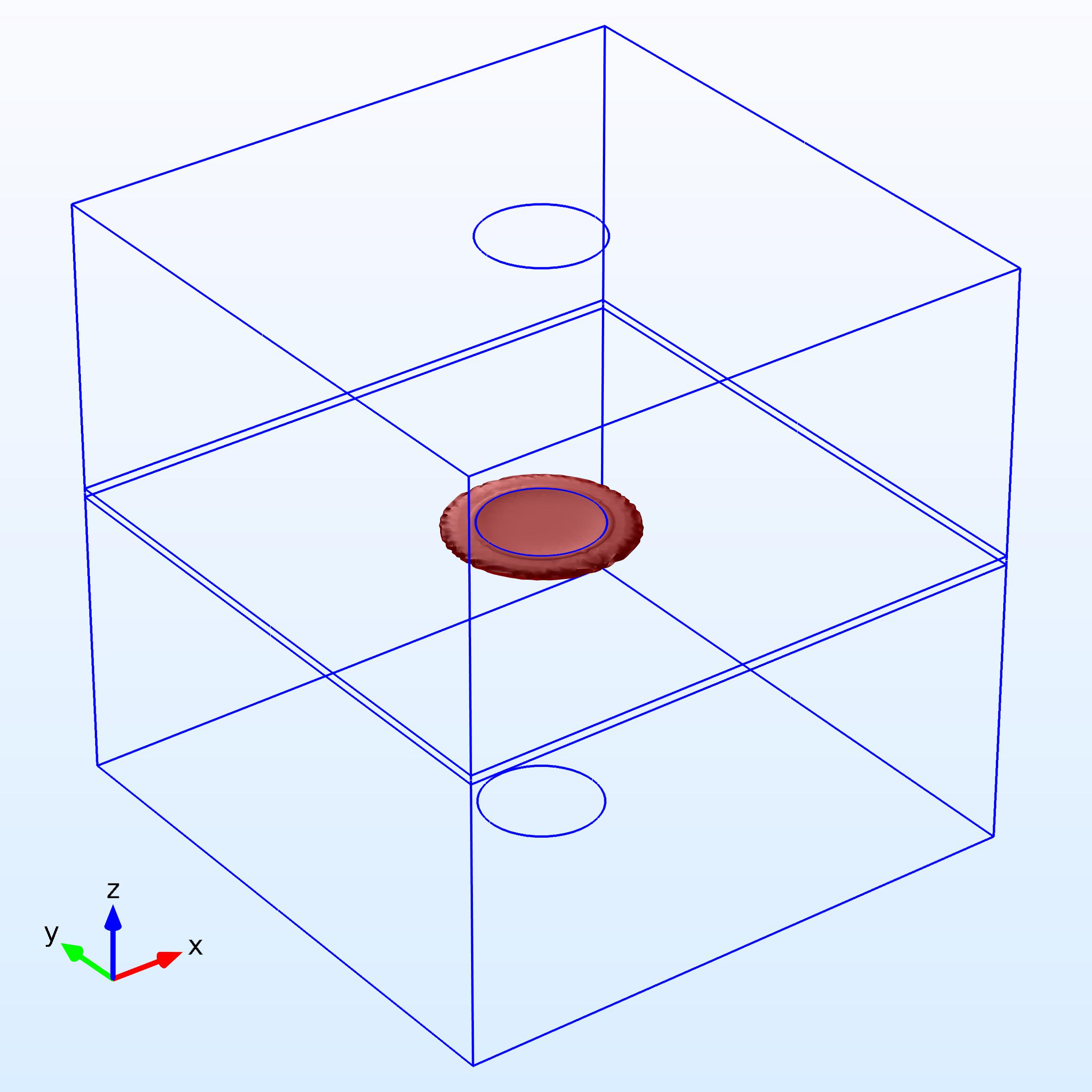}}\\
	\subfigure[$t=750$ $\mu$s]{\includegraphics[height = 6cm]{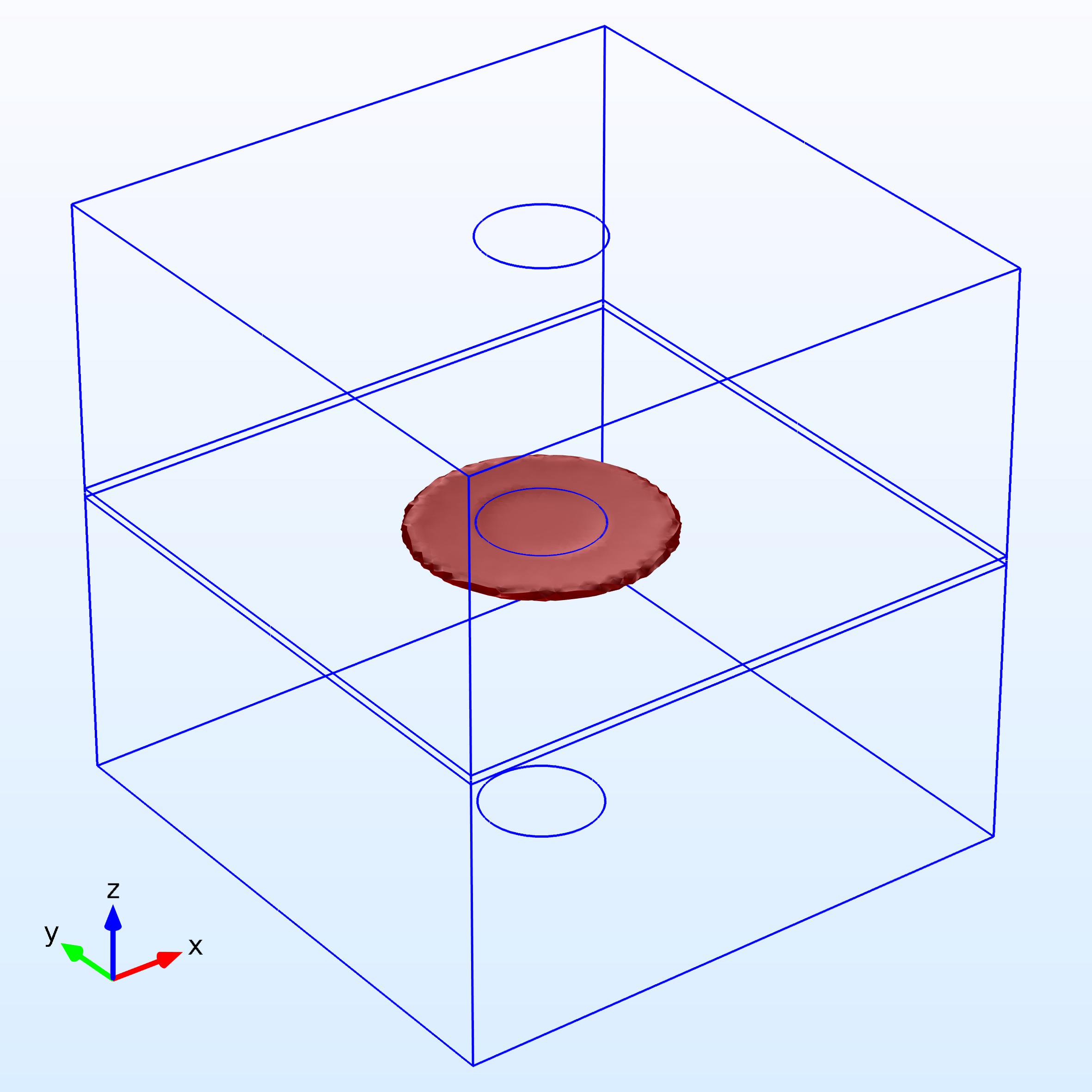}}
	\subfigure[$t=760$ $\mu$s]{\includegraphics[height = 6cm]{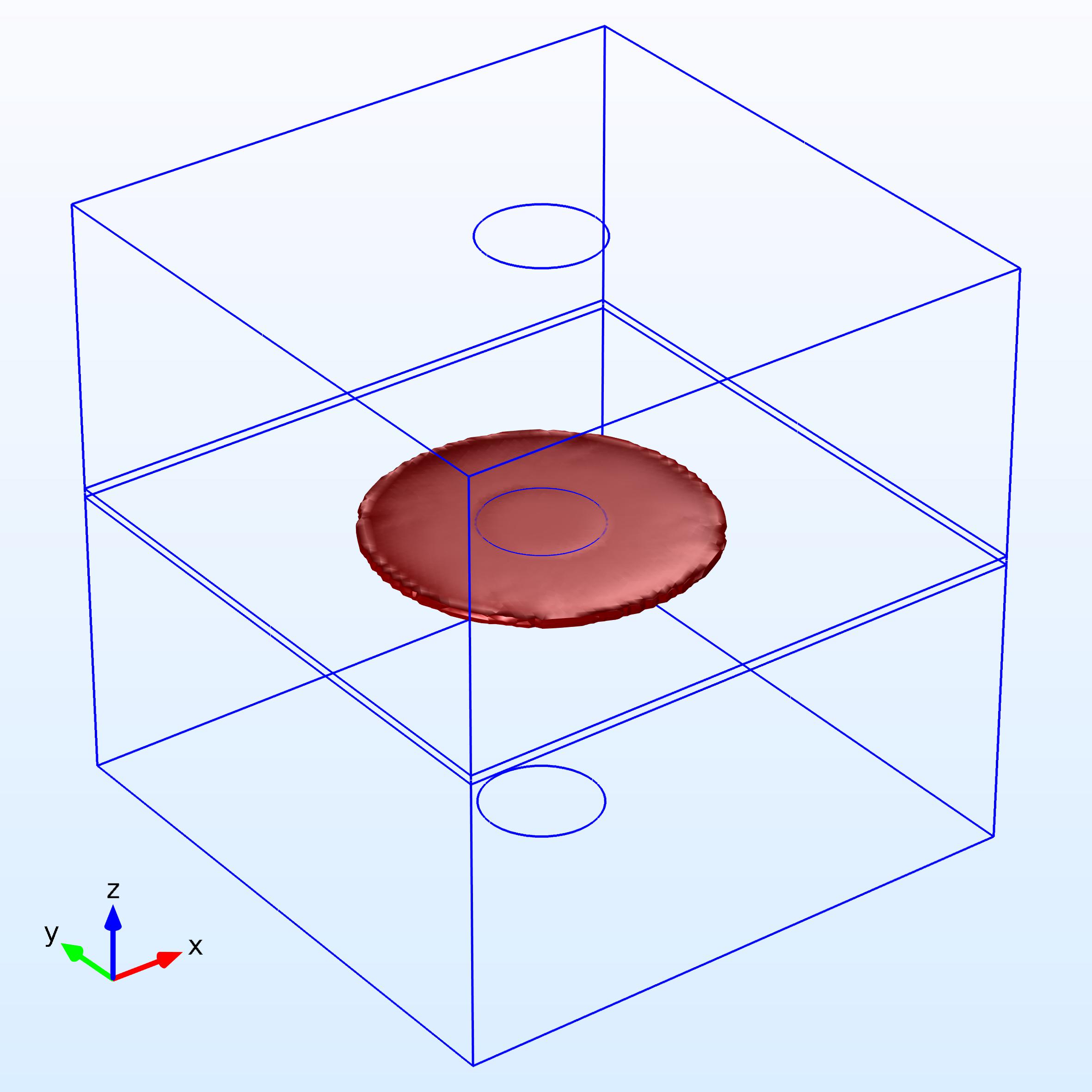}}\\
	\subfigure[$t=770$ $\mu$s]{\includegraphics[height = 6cm]{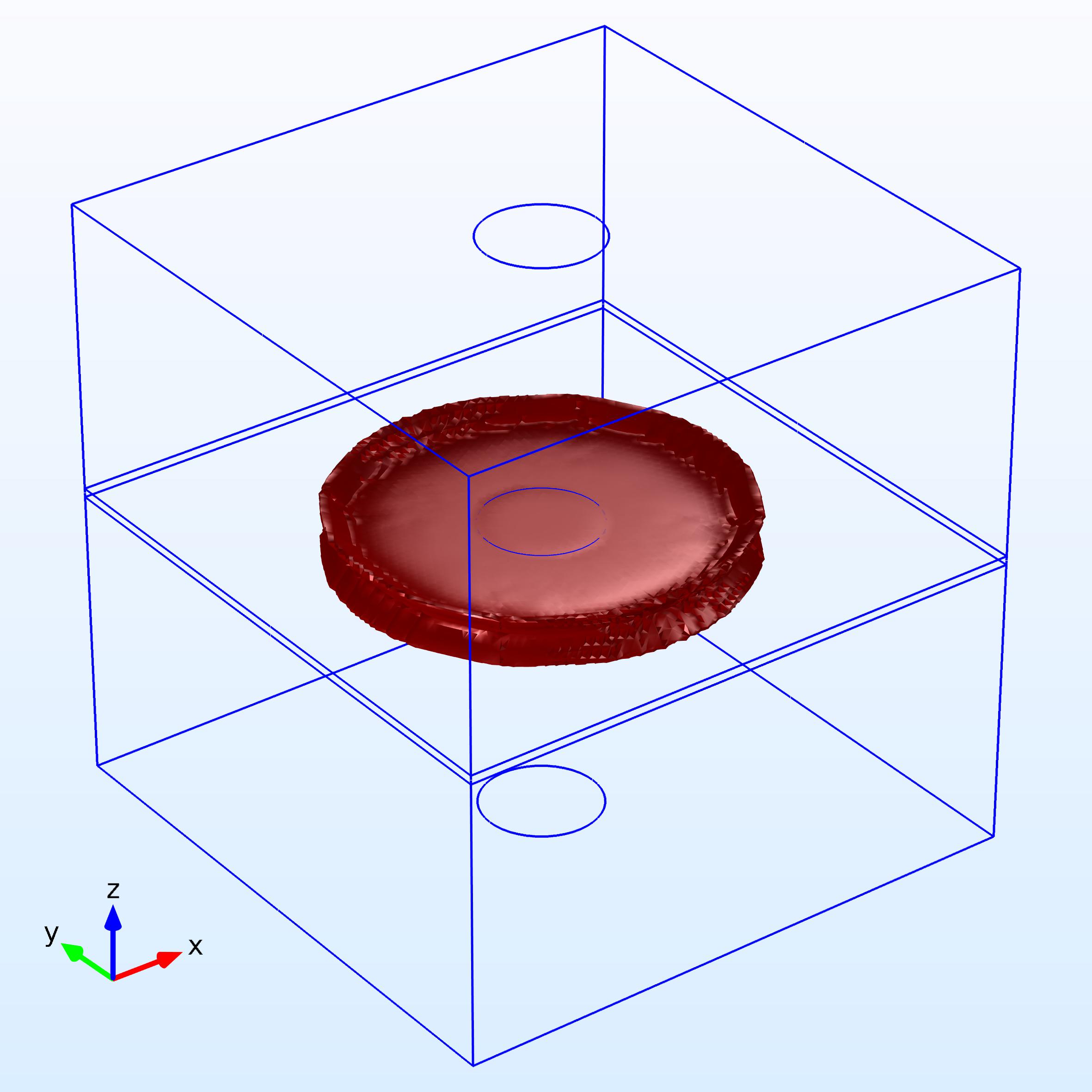}}
	\subfigure[$t=783$ $\mu$s]{\includegraphics[height = 6cm]{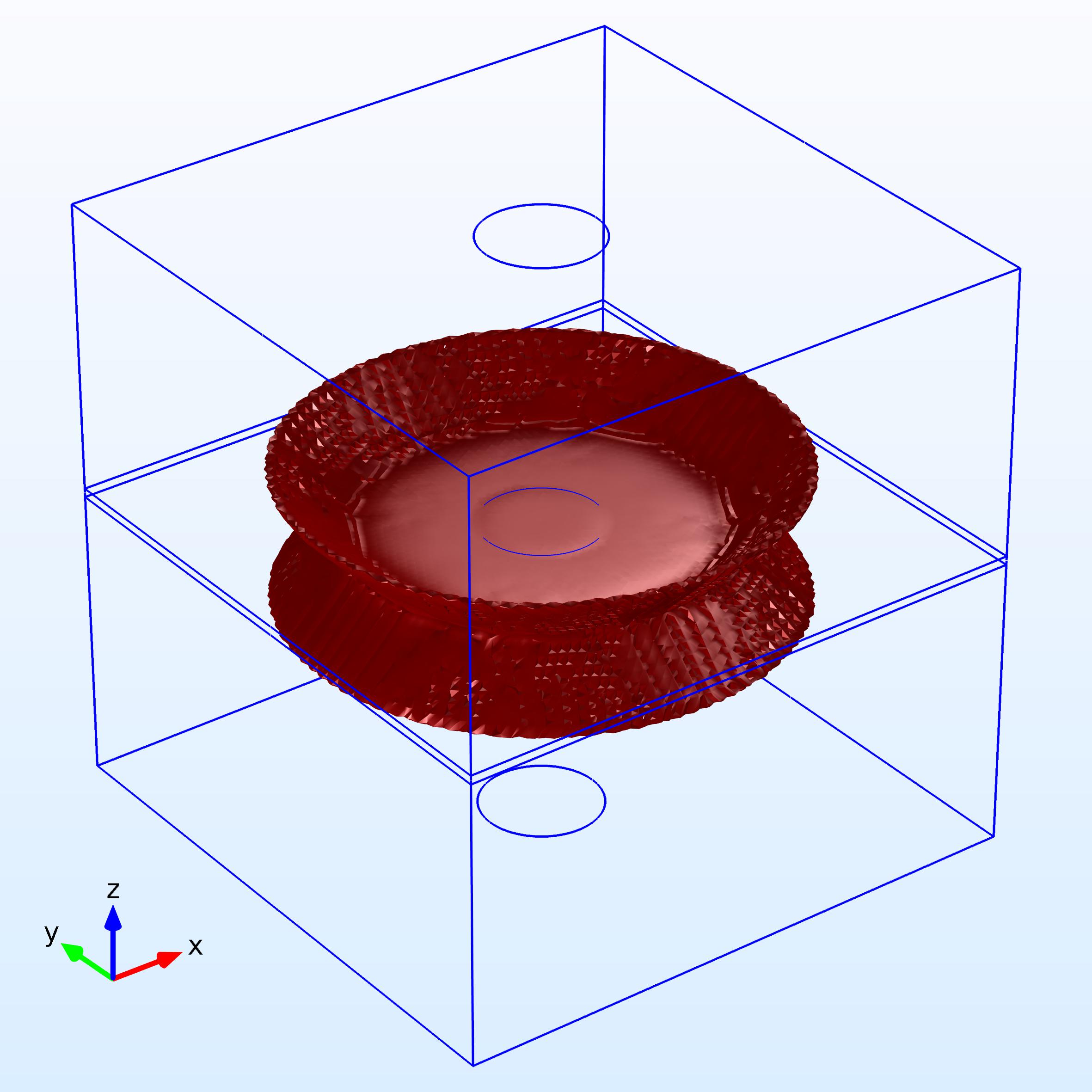}}
	\caption{Propagation of the penny-shaped crack in the 3D specimen}
	\label{Propagation of the penny-shaped crack in the 3D specimen}
	\end{figure}

Propagation patterns of the penny-shaped crack in the 3D specimen are presented in Fig. \ref{Propagation of the penny-shaped crack in the 3D specimen}. The domain with $\phi>0.95$ is displayed for the 3D crack shape. The crack starts to propagate at $t=600$ $\mu$s. The progressive crack propagation is observed when $t=600$ $\mu$s, $740$ $\mu$s, $750$ $\mu$s, and $760$ $\mu$s. The direction of the crack propagation is parallel to the x-y plane. When $t=770$ $\mu$s, the crack starts to branch. The branching crack continues to propagate towards the outer boundaries of the specimen when $t=783$ $\mu$s. 

\section{Interaction of hydraulic fracturing with natural cracks}\label{Interaction of hydraulic fracturing with natural cracks}

\begin{figure}[htbp]
\centering
\includegraphics[height = 6cm]{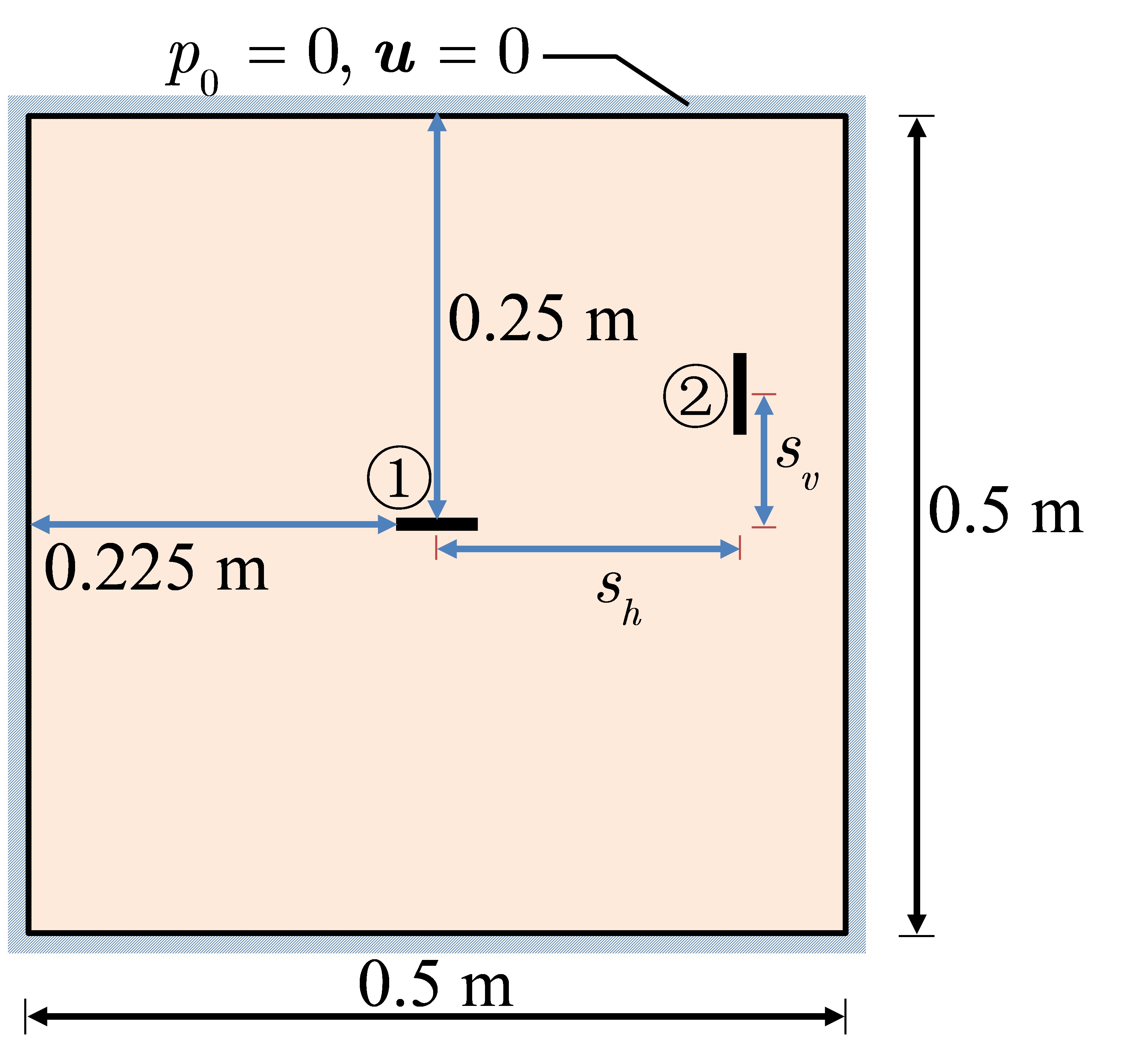}
\caption{Geometry and boundary conditions of the example for interaction between hydraulic fracturing and a natural crack}
\label{Geometry and boundary conditions of the example for interaction between hydraulic fracturing and a natural crack}
\end{figure}

For hydraulic fracturing, interaction between hydraulic and natural fractures is among the most basic and important issues. Therefore, in this section, we set a new natural crack in the calculation domain of the example of 2D dynamic crack branching to show how propagating hydraulic cracks interact with the natural crack. Geometry and boundary conditions of the example for interaction between hydraulic fracturing and the natural crack is shown in Fig. \ref{Geometry and boundary conditions of the example for interaction between hydraulic fracturing and a natural crack}.

Fluid is injected into the horizontal crack $\textcircled 1$, which is centered at (0, 0), while the vertical fracture $\textcircled 2$ is a pre-existing natural crack. The centers of the cracks $\textcircled 1$ and $\textcircled 2$ have a horizontal spacing of $s_h$ and a vertical spacing of $s_v$. In addition, the initial lengths of cracks  $\textcircled 1$ and $\textcircled 2$ are both 0.05 m. Initial history field is also established to induce the pre-existing cracks and fluid source term $q_{F}=$ 10000 kg/(m$^3\cdot \textrm s$) is set in the crack $\textcircled 1$. The parameters used for calculation are the same as those listed in Table \ref{Base parameters for the specimen with a pre-existing crack subjected to internal fluid injection}. Uniform Q4 elements are also used to discretize all the fields with the element size $h=$ $2\times 10^{-3}$ m. The time step $\Delta t$ is initially set as 0.1 $\mu$s and then changed adaptively according to the solver.

\begin{figure}[htbp]
	\centering
	\subfigure[Case 1, $t=755$ $\mu$s]{\includegraphics[height = 5cm]{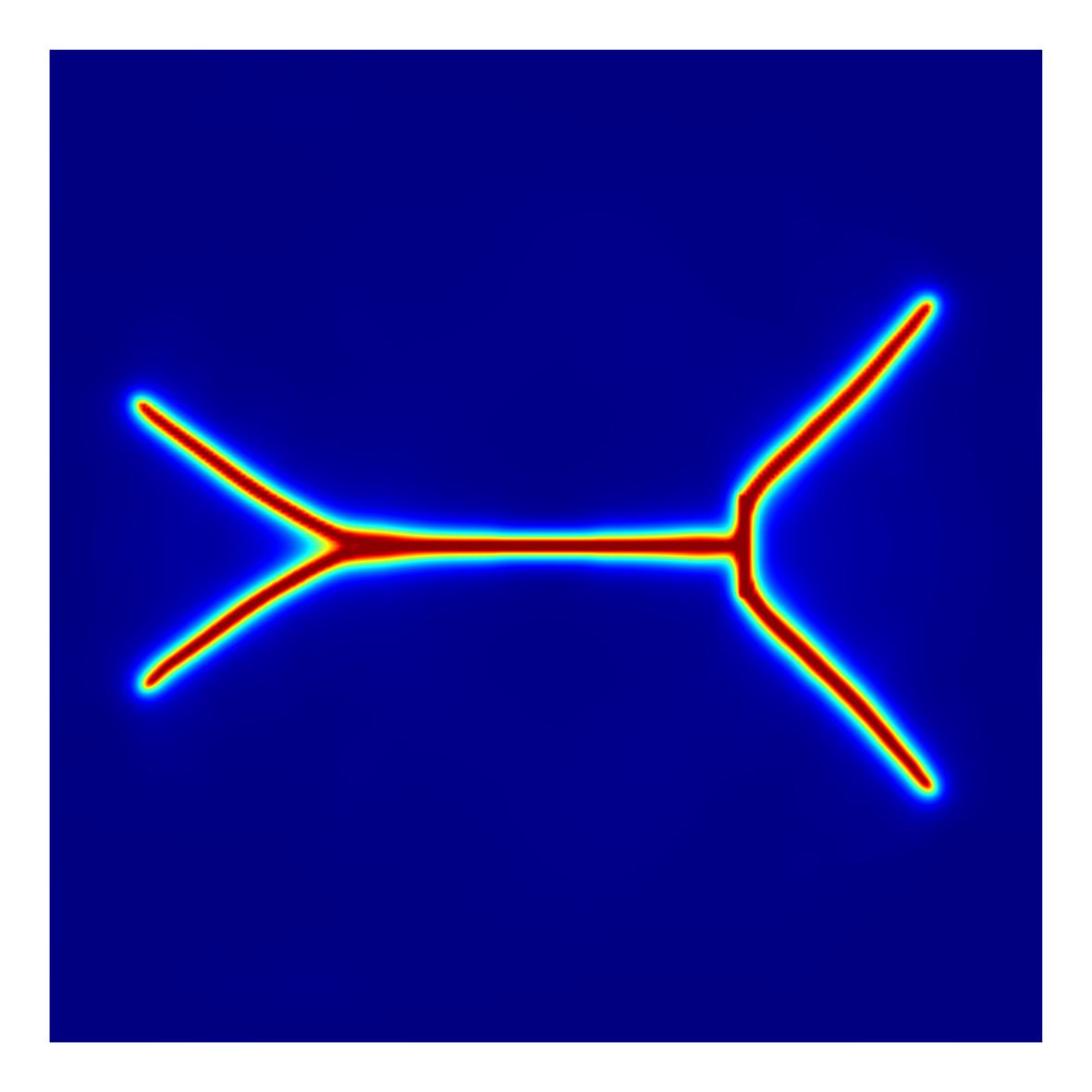}}
	\subfigure[Case 2, $t=555$ $\mu$s]{\includegraphics[height = 5cm]{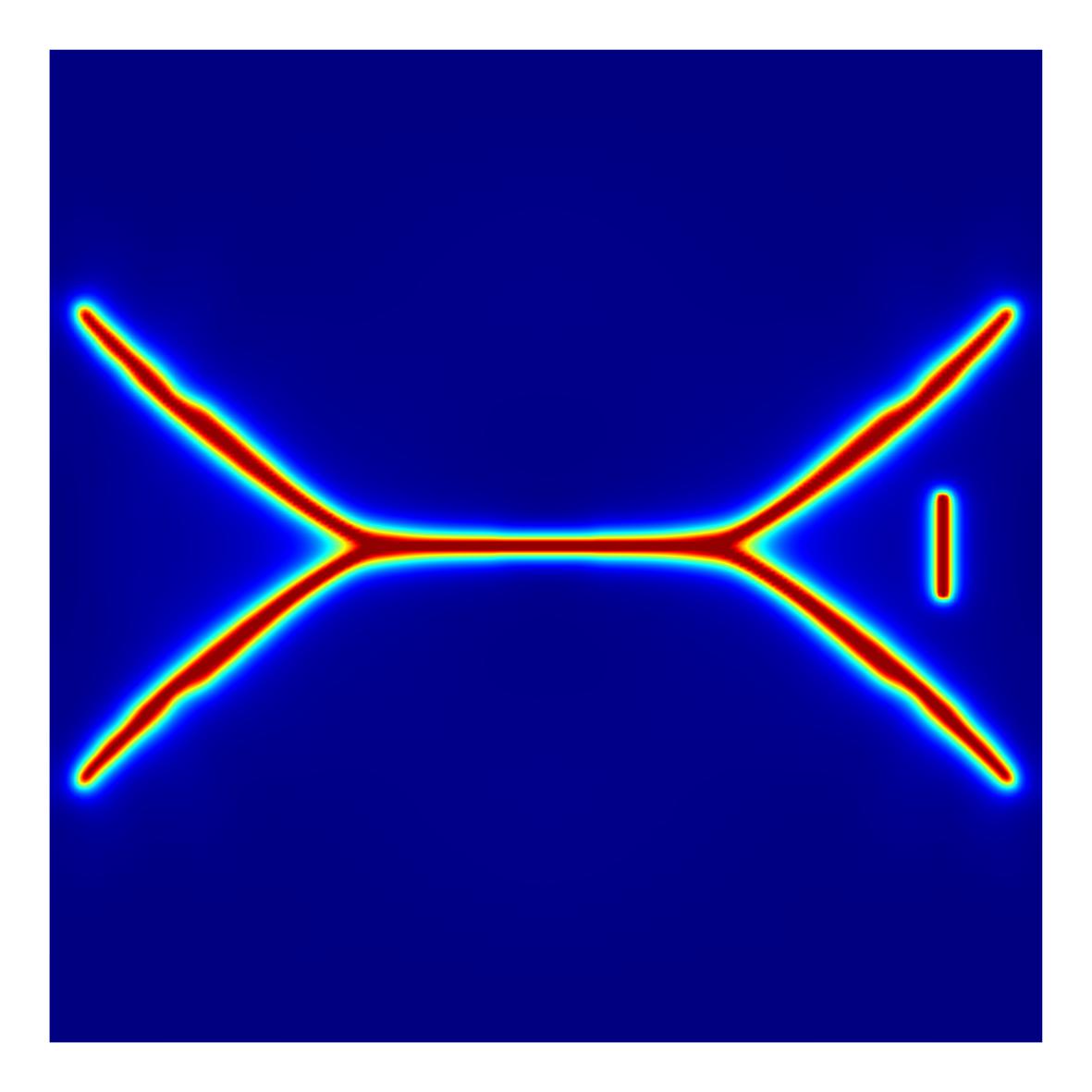}}\\
	\subfigure[Case 3, $t=539$ $\mu$s]{\includegraphics[height = 5cm]{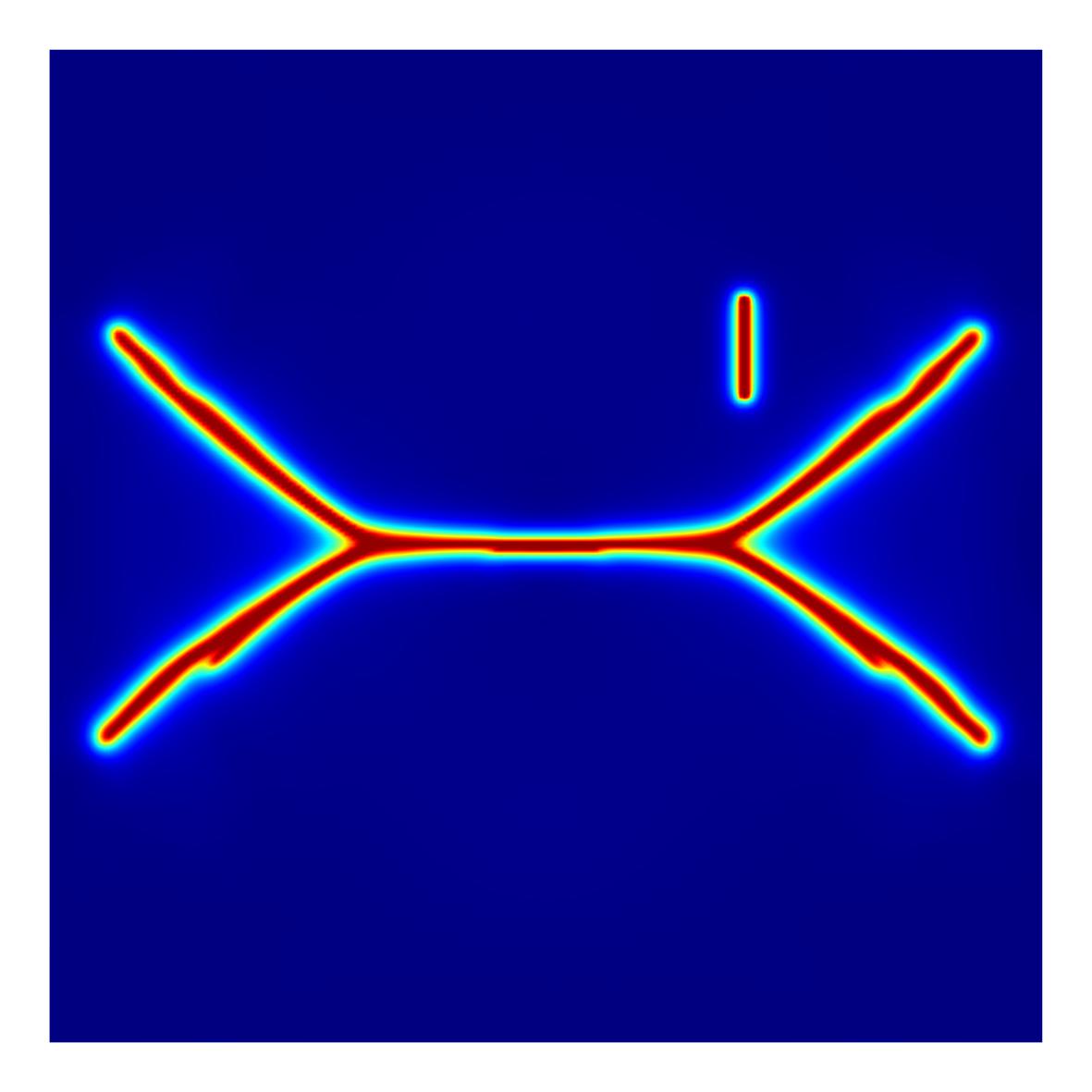}}
	\subfigure[Case 4, $t=700$ $\mu$s]{\includegraphics[height = 5cm]{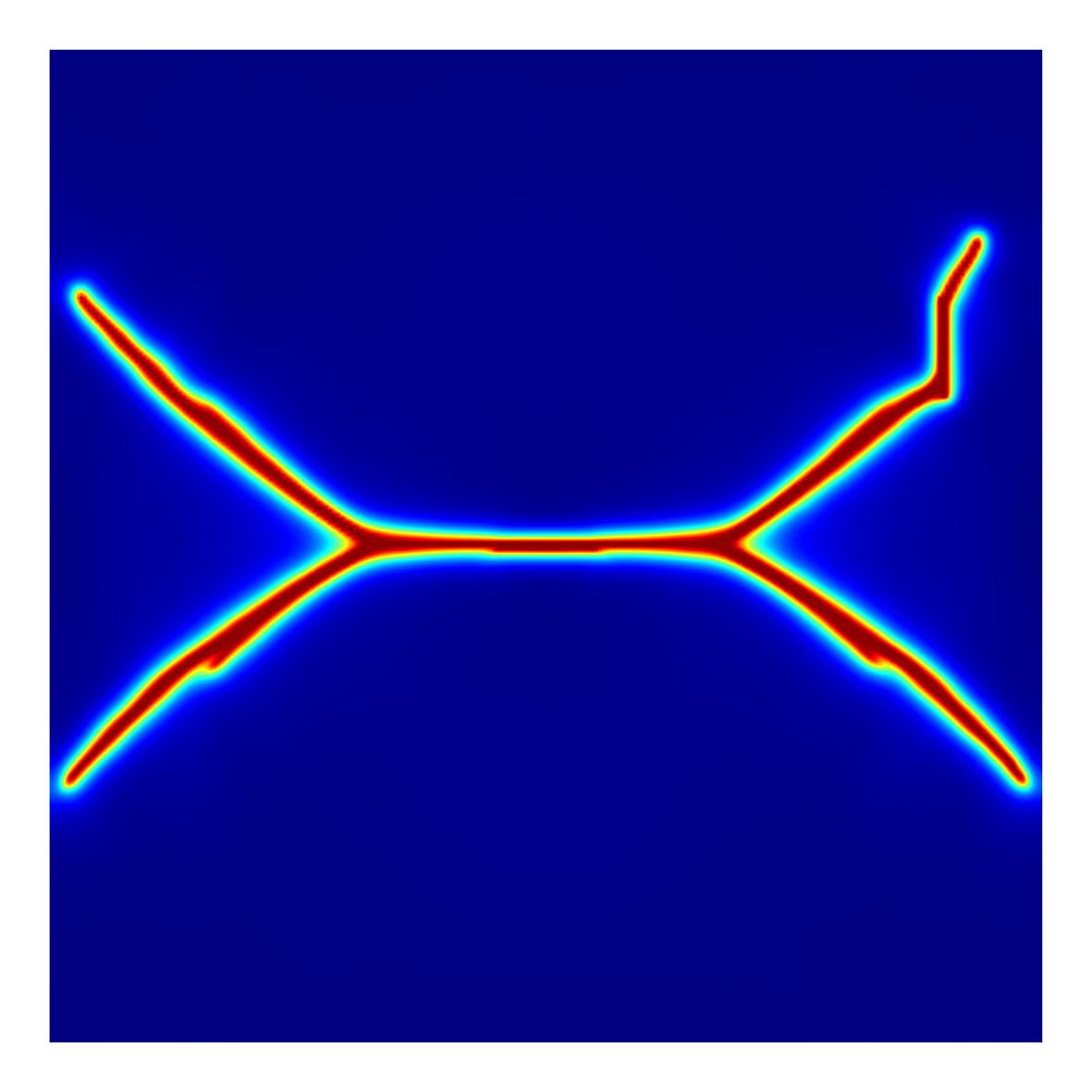}}
	\caption{Final crack patterns of the example for interaction between hydraulic fracturing and a natural crack}
	\label{Final crack patterns of the example for interaction between hydraulic fracturing and a natural crack}
	\end{figure}

For interaction between the fluid-driven and natural pre-existing cracks, we totally test four cases: Case 1 ($s_h = 0.1$ m, $s_v=0$), Case 2 ($s_h = 0.2$ m, $s_v=0$), Case 3 ($s_h = 0.1$ m, $s_v=0.1$ m), and Case 4 ($s_h = 0.2$ m, $s_v=0.1$ m). Figure \ref{Final crack patterns of the example for interaction between hydraulic fracturing and a natural crack} gives the final crack patterns of the example for interaction between hydraulic fracturing and a natural crack. As observed, crack propagation and branching in the left part of the porous domain are the similar with Fig. \ref{Progressive crack propagation in the specimen with a pre-existing crack subjected to internal fluid injection} because no natural crack exists in the left domain. However, for the right domain, different crack patterns are shown. No crack merging between the newly generated crack and pre-existing natural crack is observed for Case 2 and Case 3 because the natural crack deviates far from the expected propagation path in Fig. \ref{Progressive crack propagation in the specimen with a pre-existing crack subjected to internal fluid injection}. Meanwhile, for Case 1 and Case 4, the natural crack is captured by the newly generated crack. Subsequently, cracks continue to propagate from the tips of the natural crack after the coalescence of the hydraulic and natural cracks. Note that when the hydraulic and natural fractures join, no special treatment is required and the calculation is automatically completed because the phase-field modeling of fluid-driven cracks belongs to the class of continuous approaches to fracture and the fluid mass balance is ensured by equation (34).

\begin{figure}[htbp]
	\centering
	\subfigure[Case 1, $t=755$ $\mu$s]{\includegraphics[height = 5cm]{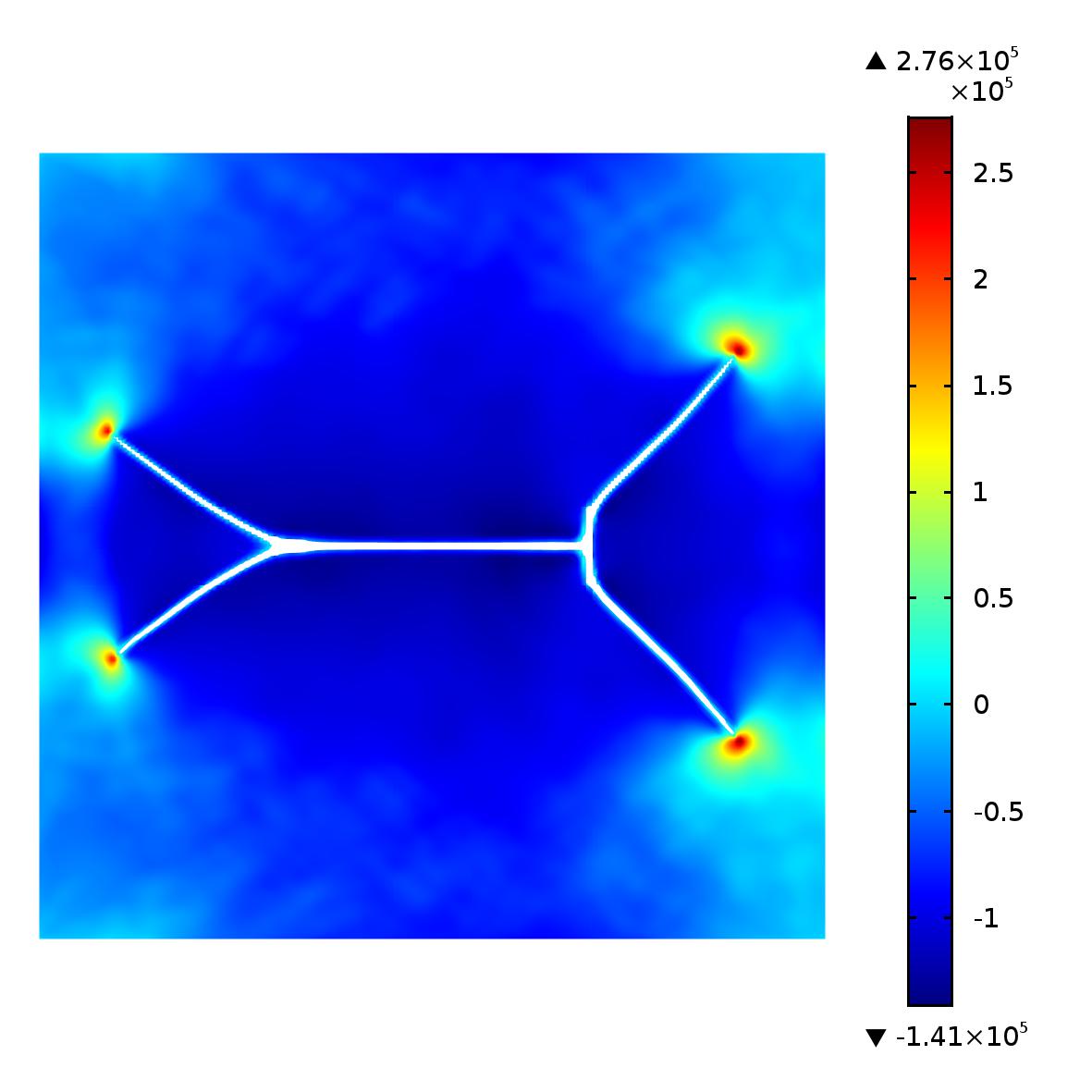}}
	\subfigure[Case 2, $t=555$ $\mu$s]{\includegraphics[height = 5cm]{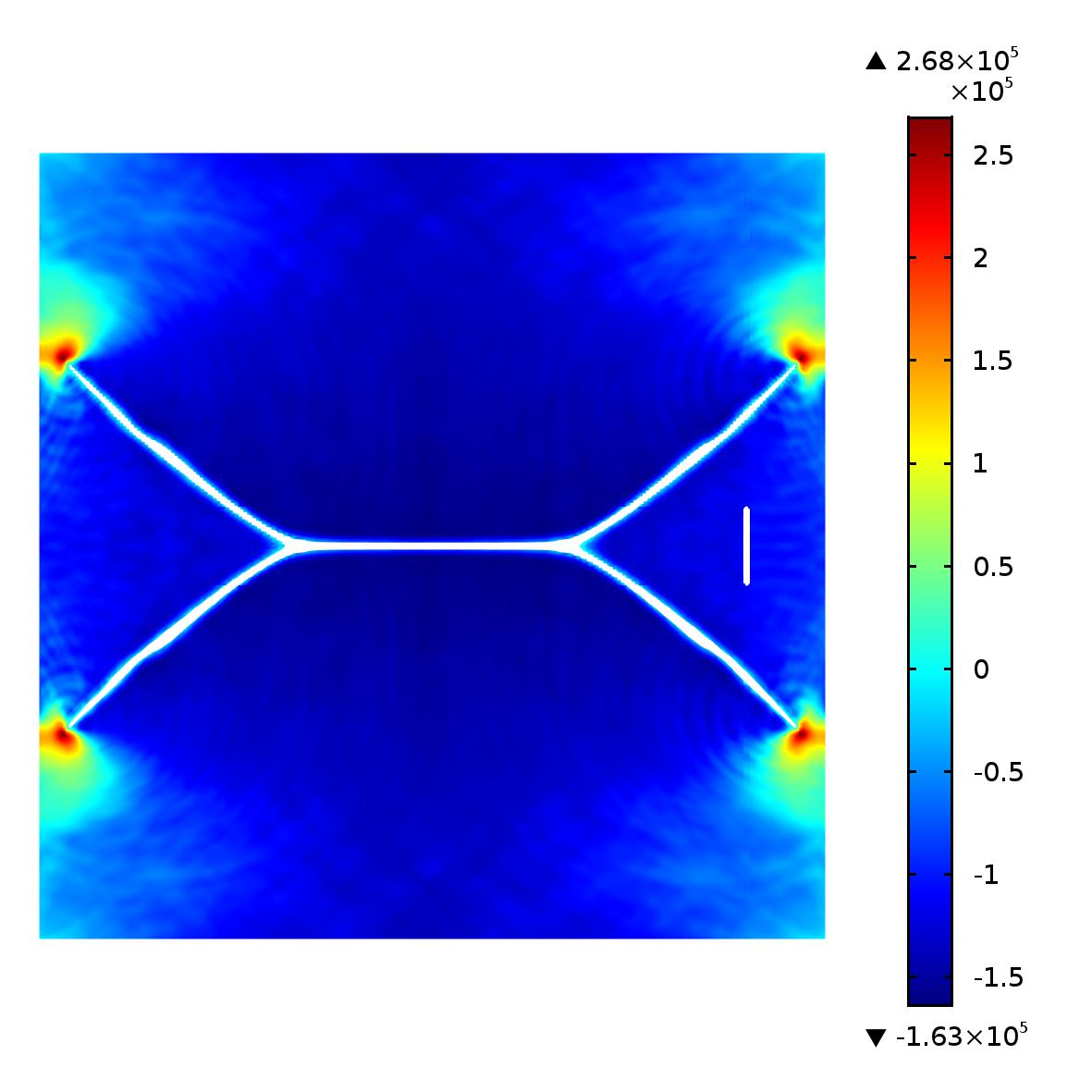}}\\
	\subfigure[Case 3, $t=539$ $\mu$s]{\includegraphics[height = 5cm]{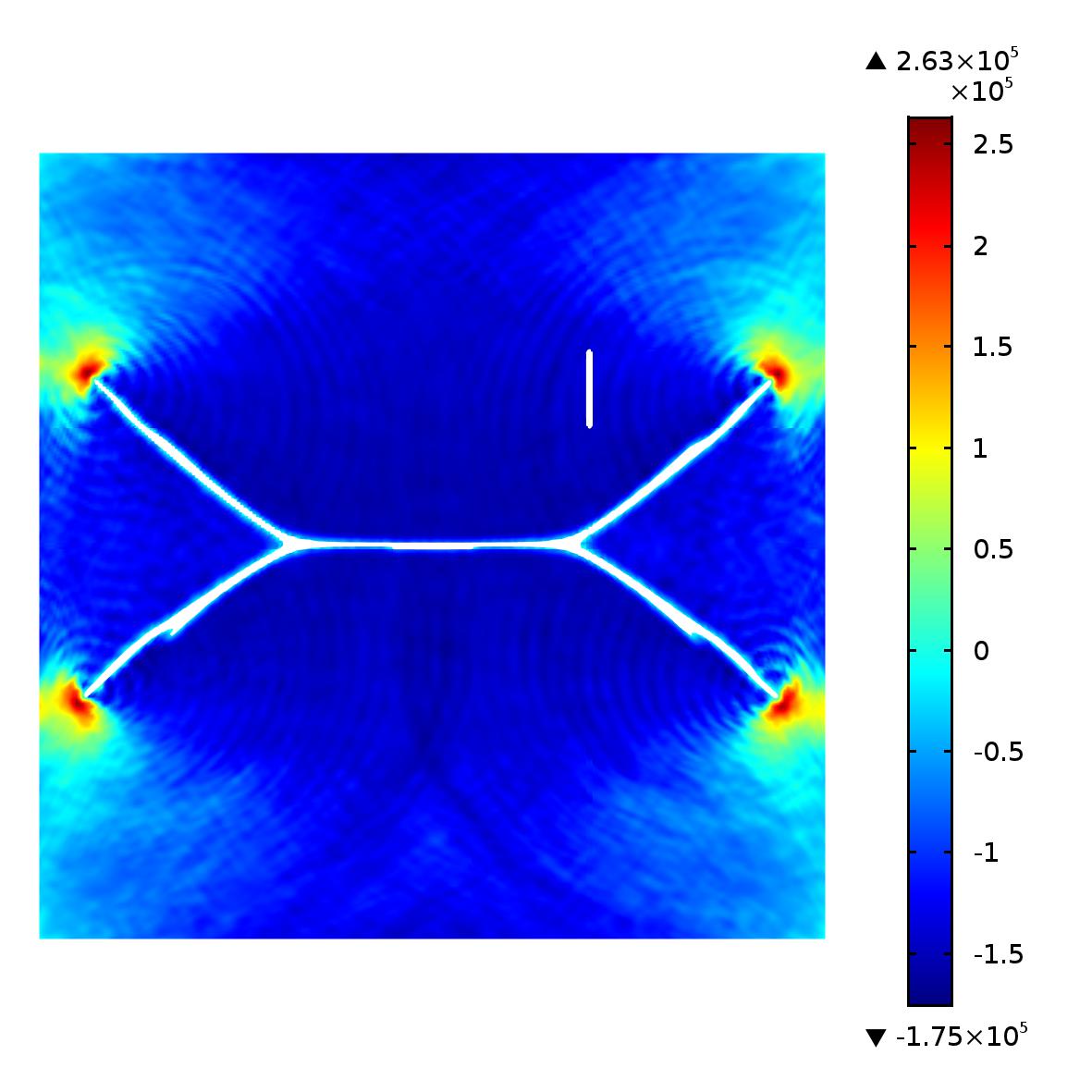}}
	\subfigure[Case 4, $t=700$ $\mu$s]{\includegraphics[height = 5cm]{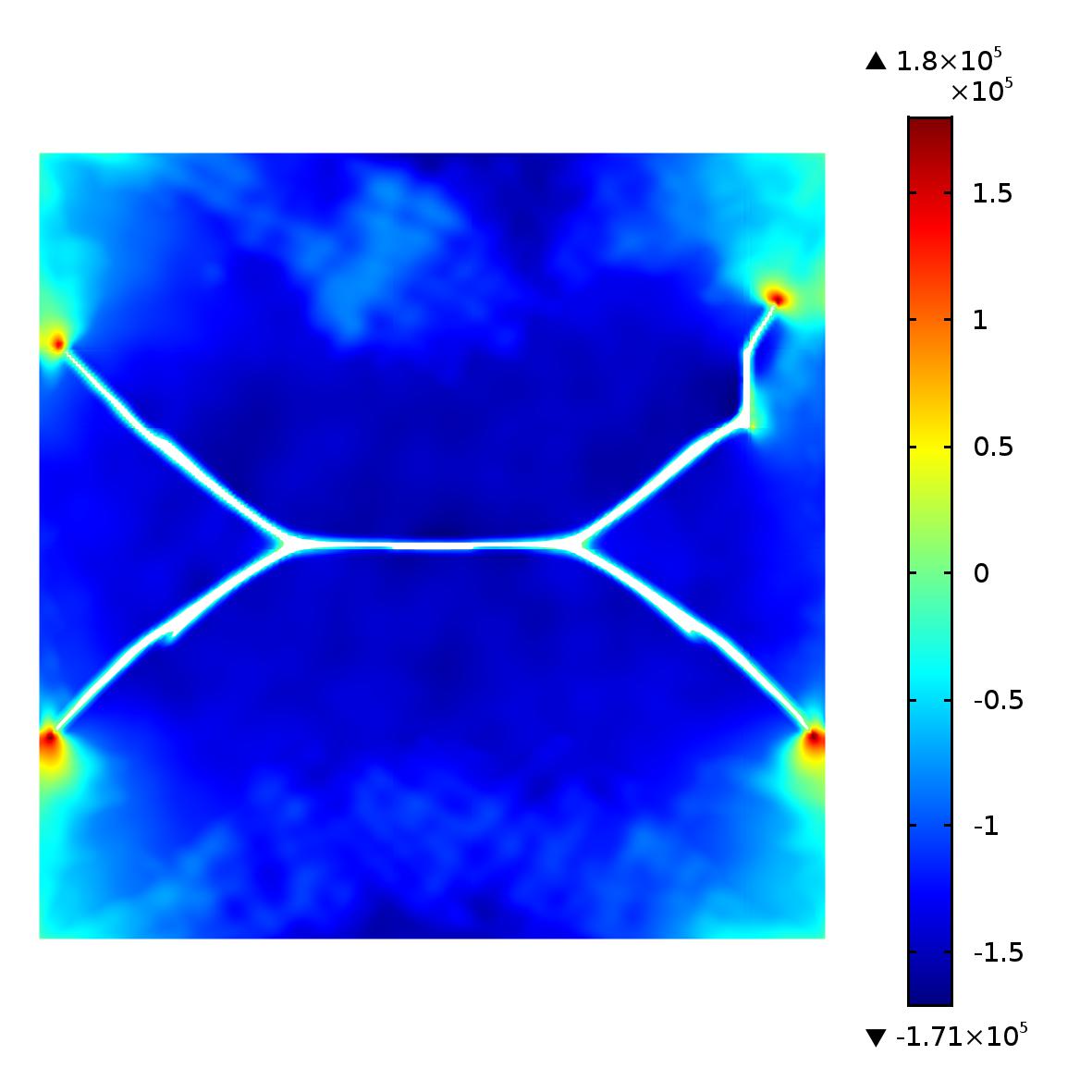}}
	\caption{First principal stress of the example for interaction between hydraulic fracturing and a natural crack (unit: Pa)}
	\label{First principal stress of the example for interaction between hydraulic fracturing and a natural crack}
	\end{figure}

\begin{figure}[htbp]
	\centering
	\includegraphics[width = 8cm]{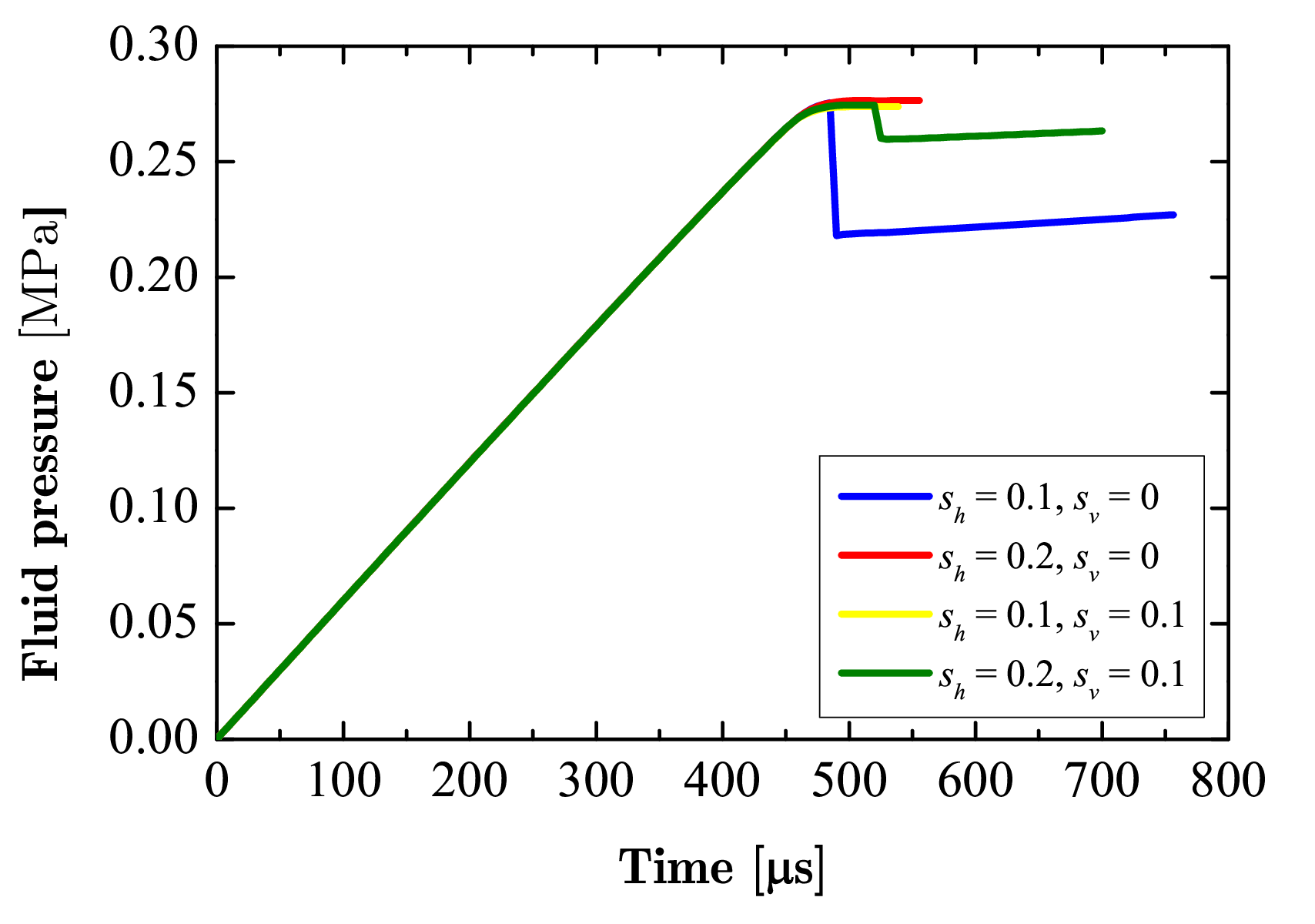}
	\caption{Fluid pressure of the example for interaction between hydraulic fracturing and a natural crack}
	\label{Fluid pressure of the example for interaction between hydraulic fracturing and a natural crack}
	\end{figure}

Figure \ref{First principal stress of the example for interaction between hydraulic fracturing and a natural crack} shows the first principal stress in the porous media. As observed, stress concentration occurs around the crack tips. In addition, Fig. \ref{Fluid pressure of the example for interaction between hydraulic fracturing and a natural crack} shows the fluid pressure at the center of crack $\textcircled 1$ for different cases. If the natural crack is not captured by the hydraulic crack, the pressure curves are similar. However, the pressure drops suddenly when the hydraulic crack joins the natural crack if the natural crack can be captured by the hydraulic crack.

\section {Conclusions}\label{Conclusions}

This paper proposes a phase field approach of modeling dynamic fracture propagation in poroelastic media. The proposed approach can be regarded as an extension of a previously proposed phase field method for only single-phasic solids to the porous media. The hydro-mechanical coupling in the porous medium is based on the classical Biot poroelasticity theory and the dynamic crack propagation is controlled by the evolution of phase field. In addition, the elastic energy drives the fracture propagation and the phase field is reconstructed as indicator functions for transiting fluid property from the intact medium to the fully broken one.

A commercial software namely, COMSOL Multiphysics is employed to implement the proposed approach. A staggered scheme is used to solve the displacement, pressure, and phase field independently. Three examples is first presented to verify the feasibility and accuracy of the proposed approach and the presented results agree well with existing analytical results. Then, some other 2D and 3D examples are presented to show dynamic crack branching and its interaction with pre-existing natural fractures. Because of the succinct implementation, the proposed approach is also suitable for researchers seeking for a quick implementation of phase field modeling of dynamic fractures. In future work, the proposed approach can be extended to predict dynamic fractures in inelastic, partially saturated, or heterogeneous porous media. 

\section*{Acknowledgement}
The financial support provided by the Sino-German (CSC-DAAD) Postdoc Scholarship Program 2016, the Natural Science Foundation of China (51474157), and RISE-project BESTOFRAC (734370) is gratefully acknowledged.

\bibliography{references}

\end{document}